\documentclass[11pt,a4paper]{article}

\usepackage{soul}
\usepackage{tikz}
\usepackage{amsthm}
\usepackage{amsmath}
\usepackage{amssymb}
\usepackage{mathrsfs}

\usepackage{geometry}
\usepackage{graphicx}
\usepackage{amsfonts}
\usepackage{epstopdf}
\usepackage{placeins}
\usepackage{enumerate}
\usepackage{enumitem}
\usepackage{subcaption}
\usepackage{microtype}
\usepackage[hidelinks]{hyperref} % avoid too many boxes with colour
\usepackage[normalem]{ulem} % enables \sout, avoids changing \emph 

\usetikzlibrary{calc,arrows.meta, positioning}

\numberwithin{equation}{section}

\allowdisplaybreaks

\newcommand{\Rmnum}[1]{\uppercase\expandafter{\romannumeral#1}} % Uppercase roman number

\def\Xint#1{\mathchoice
{\XXint\displaystyle\textstyle{#1}}%
{\XXint\textstyle\scriptstyle{#1}}%
{\XXint\scriptstyle\scriptscriptstyle{#1}}%
{\XXint\scriptscriptstyle\scriptscriptstyle{#1}}%
\!\int}
\def\XXint#1#2#3{{\setbox0=\hbox{$#1{#2#3}{\int}$ }
\vcenter{\hbox{$#2#3$ }}\kern-.6\wd0}}

\def\dashint{\Xint-}

\theoremstyle{plain}
\newtheorem{theorem}{Theorem}[section]
\newtheorem{proposition}[theorem]{Proposition}
\newtheorem{lemma}[theorem]{Lemma}
\newtheorem{corollary}[theorem]{Corollary}
\newtheorem{definition}[theorem]{Definition}

\theoremstyle{definition}
\newtheorem{remark}[theorem]{Remark}

\renewcommand{\thefootnote}{}

\AtEndDocument{
\bigskip{\footnotesize%
   \textsc{Department of Mathematics, Aarhus University, 8000 Aarhus C, Denmark} \par
  \textit{E-mail address}: \texttt{yangmengqh@gmail.com}
}}

\begin{document}

\title{\large{$p$-Poincar\'e inequalities and cutoff Sobolev inequalities on metric measure spaces}}
\author{Meng Yang}
\date{}

\maketitle

\abstract{For $p>1$, we introduce the cutoff Sobolev inequality on general metric measure spaces, and prove that there exists a metric measure space endowed with a $p$-energy that satisfies the chain condition, the volume regular condition with respect to a doubling scaling function $\Phi$, and that both the Poincar\'e inequality and the the cutoff Sobolev inequality with respect to a doubling scaling function $\Psi$ hold \emph{if and only if}
$$\frac{1}{C}\left(\frac{R}{r}\right)^p\le\frac{\Psi(R)}{\Psi(r)}\le C\left(\frac{R}{r}\right)^{p-1}\frac{\Phi(R)}{\Phi(r)}\text{ for any }r\le R.$$
In particular, given any pair of doubling functions $\Phi$ and $\Psi$ satisfying the above inequality, we construct a metric measure space endowed with a $p$-energy on which all the above conditions are satisfied. As a direct corollary, we prove that there exists a metric measure space which is $d_h$-Ahlfors regular and has $p$-walk dimension $\beta_p$ \emph{if and only if}
$$p\le\beta_p\le d_h+(p-1).$$
Our proof builds on the Laakso-type space theory, which was recently developed by Murugan [\textit{Ann. Probab.}, to appear].}

\footnote{\textsl{Date}: \today}
\footnote{The author is very grateful to Fabrice Baudoin for stimulating discussions, and to Ryosuke Shimizu for pointing out the recent paper \cite{Sas26} on (canonical) $p$-energy measures and for valuable communications on nonlinear potential theory.}

\renewcommand{\thefootnote}{\arabic{footnote}}
\setcounter{footnote}{0}

\section{Introduction}

On a complete non-compact Riemannian manifold, it was independently discovered by Grigor’yan \cite{Gri92} and Saloff-Coste \cite{Sal92,Sal95} that the following two-sided Gaussian estimate of the heat kernel
$$\frac{C_1}{V(x,\sqrt{t})}\exp\left(-C_2\frac{d(x,y)^2}{t}\right)\le p_t(x,y)\le\frac{C_3}{V(x,\sqrt{t})}\exp\left(-C_4\frac{d(x,y)^2}{t}\right),$$
where $d(x,y)$ is the geodesic distance, $V(x,r)=m(B(x,r))$ is the Riemannian measure of the open ball $B(x,r)$, is equivalent to the conjunction of the volume doubling condition and the following Poincar\'e inequality
$$\int_{B(x,r)}|f-f_{B(x,r)}|^2\mathrm{d} m\le Cr^2\int_{B(x,2r)}|\nabla f|^2\mathrm{d} m,$$
where $f_{B(x,r)}$ is the mean value of $f$ over $B(x,r)$. However, on many fractals, including the Sierpi\'nski gasket and the Sierpi\'nski carpet, the heat kernel has the following two-sided sub-Gaussian estimate
\begin{align*}\label{eq_HKbeta}\tag*{HK($\beta$)}
\frac{C_1}{V(x,t^{1/\beta})}\exp\left(-C_2\left(\frac{d(x,y)}{t^{1/\beta}}\right)^{\frac{\beta}{\beta-1}}\right)\le p_t(x,y)\le\frac{C_3}{V(x,t^{1/\beta})}\exp\left(-C_4\left(\frac{d(x,y)}{t^{1/\beta}}\right)^{\frac{\beta}{\beta-1}}\right),\nonumber
\end{align*}
where $\beta$ is a new parameter called the walk dimension, which is always strictly greater than 2 on fractals. For example, $\beta=\log5/\log2$ on the Sierpi\'nski gasket (see \cite{BP88,Kig89}), $\beta\approx2.09697$ on the Sierpi\'nski carpet (see \cite{BB89,BB90,BBS90,BB92,KZ92,HKKZ00}). On general metric measure spaces that may exhibit different scaling behavior, such as fractal-like domains and fractal-like cable systems (see \cite{BB00,BBK06}), the heat kernel would display scale-dependent behavior as follows.
\begin{equation*}\label{eq_HKPsi}\tag*{HK($\Psi$)}
\frac{C_1}{V\left(x,\Psi^{-1}(t)\right)}\exp\left(-\Upsilon\left(C_2d(x,y),t\right)\right)\le p_t(x,y)\le\frac{C_3}{V\left(x,\Psi^{-1}(t)\right)}\exp\left(-\Upsilon\left(C_4d(x,y),t\right)\right),
\end{equation*}
where
$$\Upsilon(R,t)=\sup_{s\in(0,+\infty)}\left(\frac{R}{s}-\frac{t}{\Psi(s)}\right),$$
here $\Psi$ is a scaling function, which is also assumed to be doubling. \ref{eq_HKbeta} corresponds to the homogeneous case $\Psi(r)=r^\beta$. It was proved in \cite{BB04,BBK06,GHL15} that on a general metric measure space, under the volume doubling condition, \ref{eq_HKPsi} is equivalent to the conjunction of the following Poincar\'e inequality
$$\int_{B(x,r)}|f-f_{B(x,r)}|^2\mathrm{d} m\le C\Psi(r)\int_{B(x,2r)}\mathrm{d}\Gamma(f,f),$$
where $\Gamma(f,f)$ is the energy measure of $f$ (see \cite{FOT11} for more details), and the cutoff Sobolev inequality or the generalized capacity condition, with the scaling function $\Psi$. It was proved in \cite{Mur24a} that the cutoff Sobolev inequality can be replaced by the following simplified version. There exist $C_1, C_2>0$, $A>1$, such that for any ball $B(x,r)$, there exists a cutoff function $\phi$ for $B(x,r)\subseteq B(x,Ar)$ such that for any $f$, we have
$$\int_{B(x,Ar)}|\widetilde{f}|^2\mathrm{d}\Gamma(\phi,\phi)\le C_1\int_{B(x,Ar)}\mathrm{d}\Gamma(f,f)+\frac{C_2}{\Psi(r)}\int_{B(x,Ar)}|f|^2\mathrm{d} m,$$
where $\widetilde{f}$ is a quasi-continuous modification of $f$.

However, \ref{eq_HKPsi} would impose certain constraints on $\Psi$. In the homogeneous case, assume that the space is $d_h$-Ahlfors regular for some $d_h>0$, it was proved in \cite[Corollary 4.9]{GHL03} that under the chain condition\footnote{It was recently proved in \cite[Corollary 1.8]{Mur20} that the chain condition can be dropped.}, \ref{eq_HKbeta} implies
$$2\le \beta\le d_h+1,$$
see also \cite[Theorem 1]{Bar04} for the infinite graph case. In the general case, assume that the space satisfies the volume regular condition \hypertarget{eq_VPhi_intro}{\textrm{V}($\Phi$)} with a scaling function $\Phi$, that is, $V(x,r)\asymp\Phi(r)$. Here $\Phi$ is also assumed to be doubling. It was recently proved in \cite{Mur24} that \ref{eq_HKPsi} implies
$$\frac{1}{C}\left(\frac{R}{r}\right)^2\le \frac{\Psi(R)}{\Psi(r)}\le C\frac{R\Phi(R)}{r\Phi(r)}\text{ for any }r\le R.$$

A natural question is whether the converse also holds. In the homogeneous case, for any $d_h, \beta$ satisfying $2\le \beta\le d_h+1$, Barlow \cite{Bar04} constructed a $d_h$-Ahlfors regular space on which \ref{eq_HKbeta} holds. He used a construction similar to the Laakso space in \cite{Laa00}, adapted to the graph setting. In the general case, for any pair of doubling functions $\Phi,\Psi$ satisfying the above inequality, Murugan \cite{Mur24} recently constructed a metric measure space on which both \hyperlink{eq_VPhi_intro}{\textrm{V}($\Phi$)} and \ref{eq_HKPsi} hold. He developed a general Laakso-type space theory, and proved that the Poincar\'e inequality and the cutoff Sobolev inequality, whose conjunction is equivalent to \ref{eq_HKPsi}, holds on this space.

The preceding results are formulated within the Dirichlet form framework, which generalizes the classical Dirichlet integral $\int_{\mathbb{R}^d}|\nabla f(x)|^2\mathrm{d} x$ in $\mathbb{R}^d$. For general $p>1$, extending the classical $p$-energy $\int_{\mathbb{R}^d}|\nabla f(x)|^p\mathrm{d} x$ in $\mathbb{R}^d$, as initiated by \cite{HPS04}, the study of $p$-energy on fractals and general metric measure spaces has been recently advanced considerably, see \cite{CGQ22,Shi24,BC23,MS25,Kig21book,CGYZ26,AB25,AES25a}. In this setting, a new parameter $\beta_p$, called the $p$-walk dimension, naturally arises in connection with a $p$-energy. Notably, $\beta_2$ coincides with $\beta$ in \ref{eq_HKbeta}. However, for general $p>1$, there is no known analogue of the heat kernel estimate, instead, the $p$-walk dimension can be characterized in terms of the following Poincar\'e inequality
\begin{align*}\label{eq_PIbeta_intro}\tag*{\textrm{PI}($\beta_p$)}
\int_{B(x,r)}|f-f_{B(x,r)}|^p\mathrm{d} m\le Cr^{\beta_p}\int_{B(x,2r)}\mathrm{d}\Gamma(f),
\end{align*}
where $\Gamma(f)$ is the $p$-energy measure of $f$, and the following cutoff Sobolev inequality
\begin{align*}\label{eq_CSbeta_intro}\tag*{\textrm{CS}($\beta_p$)}
\int_{B(x,Ar)}|\widetilde{f}|^p\mathrm{d}\Gamma(\phi)\le C_1\int_{B(x,Ar)}\mathrm{d}\Gamma(f)+\frac{C_2}{r^{\beta_p}}\int_{B(x,Ar)}|f|^p\mathrm{d} m,
\end{align*}
where $\phi$ is a cutoff function for $B(x,r)\subseteq B(x,Ar)$. By taking $f\equiv1$ in $B(x,Ar)$, it is easy to see that \ref{eq_CSbeta_intro} implies the following capacity upper bound
\begin{align*}\label{eq_ucapbeta_intro}\tag*{\textrm{cap}$(\beta_p)_\le$}
\mathrm{cap}(B(x,r),X\backslash B(x,Ar))\le C_2\frac{V(x,Ar)}{r^{\beta_p}},
\end{align*}
where $\mathrm{cap}(\cdot,\cdot)$ is the capacity associated with a $p$-energy. On general metric measure spaces, scale-dependent formulations of the Poincar\'e inequaity \hypertarget{eq_PIPsi_intro}{\textrm{PI}$(\Psi)$}, the cutoff Sobolev inequality \hypertarget{eq_CSPsi_intro}{\textrm{CS}$(\Psi)$}, and the capacity upper bound \hypertarget{eq_capPsi_intro}{\textrm{cap}$(\Psi)_\le$} can be considered. Similar to the case $p=2$, one can ask the following question.

\vspace{5pt}
\noindent{\textbf{Question 1:}} In the homogeneous case, assume that the space is $d_h$-Ahlfors regular for some $d_h>0$, and that both \ref{eq_PIbeta_intro} and \ref{eq_CSbeta_intro} hold. What are the constraints on $d_h$ and $\beta_p$? In the general case, assume that \hyperlink{eq_VPhi_intro}{\textrm{V}($\Phi$)}, \hyperlink{eq_PIPsi_intro}{\textrm{PI}$(\Psi)$} and \hyperlink{eq_CSPsi_intro}{\textrm{CS}$(\Psi)$} hold. What are the constraints on $\Phi$ and $\Psi$?
\vspace{5pt}

In the homogeneous case, it was proved in \cite{Bau24,Shi24a} that under the chain condition, the conjunction of \ref{eq_PIbeta_intro} and \ref{eq_ucapbeta_intro} implies
$$p\le\beta_p\le d_h+(p-1).$$
The initial purpose of this paper is to consider the general case. We will prove that under the chain condition, the conjunction of \hyperlink{eq_VPhi_intro}{\textrm{V}($\Phi$)}, \hyperlink{eq_PIPsi_intro}{\textrm{PI}$(\Psi)$} and \hyperlink{eq_capPsi_intro}{\textrm{cap}$(\Psi)_\le$} implies
$$\frac{1}{C}\left(\frac{R}{r}\right)^p\le\frac{\Psi(R)}{\Psi(r)}\le C\left(\frac{R}{r}\right)^{p-1}\frac{\Phi(R)}{\Phi(r)}\text{ for any }r\le R,$$
which is Proposition \ref{prop_PhiPsi_bound}. It is also natural to ask whether the converse holds as follows.

\vspace{5pt}
\noindent {\textbf{Question 2:}} For any pair of doubling functions $\Phi,\Psi$ satisfying the above inequality, does there exist a metric measure space endowed with a $p$-energy, on which \hyperlink{eq_VPhi_intro}{\textrm{V}($\Phi$)}, \hyperlink{eq_PIPsi_intro}{\textrm{PI}$(\Psi)$} and \hyperlink{eq_CSPsi_intro}{\textrm{CS}$(\Psi)$} hold?
\vspace{5pt}

The main purpose of this paper is to give a positive answer to this question, which is Theorem \ref{thm_main}. Very recently, Eriksson-Bique \cite{Eri26} proved that, under the volume doubling condition and \hyperlink{eq_PIPsi_intro}{\textrm{PI}$(\Psi)$}, \hyperlink{eq_CSPsi_intro}{\textrm{CS}$(\Psi)$} and \hyperlink{eq_capPsi_intro}{\textrm{cap}$(\Psi)_\le$} are equivalent. This result resolves the long-standing resistance conjecture (see \cite{GHL14}). See also Murugan \cite{Mur26} for the resolution of the resistance conjecture in the non-local setting (see \cite{GHH18}). However, in the present paper, we will establish \hyperlink{eq_CSPsi_intro}{\textrm{CS}$(\Psi)$} \emph{directly}, without relying on these results.

As an application of our results, we will consider the critical exponent $\alpha_p$ of Besov spaces $B^{p,\alpha}$ for $p>1$, which was recently studied in \cite{ABCRST3,Bau24}. We will provide a characterization of the range in which the critical exponent $\alpha_p$ is attainable.

To conclude the introduction, let us give a brief overview of our solution to Question 2 based on the Laakso-type space theory developed in \cite{Mur24}.

Firstly, we use an $\mathbb{R}$-tree $\mathcal{T}(\mathbf{b})$ introduced in \cite{Mur24}, where $\mathbf{b}$ is a function that determines the branching numbers of the tree and is thus referred to as a branching function. Owing to the tree property, suitable differential structure can be introduced on $\mathcal{T}(\mathbf{b})$, from which we can construct a $p$-energy, obtain the Poincar\'e inequality and the capacity upper bound in a relatively straightforward manner. See also \cite{BC23,BC24,CGYZ26} for the Vicsek set case. However, in contrast to the case $p=2$, which was considered in \cite{Mur24} and for which numerous equivalent characterizations of heat kernel estimates are available, we will prove the cutoff Sobolev inequality \emph{directly}.

Secondly, an ultrametric space $\mathcal{U}(\mathbf{g})$ is introduced to define a product space $\mathcal{P}(\mathbf{g},\mathbf{b})=\mathcal{U}(\mathbf{g})\times\mathcal{T}(\mathbf{b})$. The Laakso-type space $\mathcal{L}(\mathbf{g},\mathbf{b})$ is then obtained as the quotient of $\mathcal{P}(\mathbf{g},\mathbf{b})$ with respect to an suitable equivalence relation, or roughly speaking, by gluing different copies of $\mathcal{T}(\mathbf{b})$ along ``carefully" chosen ``wormholes" determined by the function $\mathbf{g}$. For this reason, $\mathbf{g}$ is referred to as a gluing function. We will introduce a metric $\mathbf{d}_{\mathcal{L}(\mathbf{g},\mathbf{b})}$, which is slightly different from the one used in \cite{Mur24}. We will show that $(\mathcal{L}(\mathbf{g},\mathbf{b}),\mathbf{d}_{\mathcal{L}(\mathbf{g},\mathbf{b})})$ is a geodesic space and give a characterization of its geodesics, a topic not addressed in \cite{Mur24}. Let $\mathbf{m}_{\mathcal{L}(\mathbf{g},\mathbf{b})}$ be the pushforward of the product measure on $\mathcal{P}(\mathbf{g},\mathbf{b})$, then $(\mathcal{L}(\mathbf{g},\mathbf{b}),\mathbf{d}_{\mathcal{L}(\mathbf{g},\mathbf{b})},\mathbf{m}_{\mathcal{L}(\mathbf{g},\mathbf{b})})$ is a metric measure space on which \hyperlink{eq_VPhi_intro}{\textrm{V}($\Phi$)} holds, where $\Phi=V_\mathbf{g} V_\mathbf{b}$ and $V_\mathbf{g}, V_\mathbf{b}$ are determined by $\mathbf{g}, \mathbf{b}$, respectively.

Thirdly, due to the ultrametric property of $\mathcal{U}(\mathbf{g})$ and the ``carefully" chosen ``wormholes" in $\mathcal{T}(\mathbf{b})$, the differential structure on $\mathcal{T}(\mathbf{b})$ can be transferred to $\mathcal{L}(\mathbf{g},\mathbf{b})$, from which we construct a $p$-energy on $\mathcal{L}(\mathbf{g},\mathbf{b})$. We will prove the Poincar\'e inequality \hyperlink{eq_PIPsi_intro}{\textrm{PI}$(\Psi)$} using the technique of the pencil of curves, following an argument similar to that in \cite{Mur24}, adapted to the $p$-energy setting, and the cutoff Sobolev inequality \hyperlink{eq_CSPsi_intro}{\textrm{CS}$(\Psi)$}, the capacity upper bound \hyperlink{eq_capPsi_intro}{\textrm{cap}$(\Psi)_\le$} using the corresponding result on $\mathcal{T}(\mathbf{b})$, here $\Psi$ is given by $\Psi(r)=r^{p-1}V_\mathbf{b}(r)$. To complete the proof, it remains to find $\mathbf{g}, \mathbf{b}$ corresponding to given $\Phi, \Psi$.

Throughout this paper, $p>1$ is fixed. The letters $C,C_1,C_2,C_A,C_B$ will always refer to some positive constants and may change at each occurrence. The sign $\asymp$ means that the ratio of the two sides is bounded from above and below by positive constants. The sign $\lesssim$ ($\gtrsim$) means that the LHS is bounded by positive constant times the RHS from above (below). We will use the notation $\lfloor x\rfloor$ ($\lceil x\rceil$) to denote the largest integer less than or equal to (the smallest integer greater than or equal to) $x\in \mathbb{R}$. We use $\#A$ to denote the cardinality of a set $A$.

\section{Statement of main results}

We say that a function $\Phi:[0,+\infty)\to[0,+\infty)$ is doubling if $\Phi$ is a homeomorphism, which implies that $\Phi$ is strictly increasing continuous and $\Phi(0)=0$, and there exists $C_\Phi>1$, called a doubling constant of $\Phi$, such that $\Phi(2r)\le C_\Phi\Phi(r)$ for any $r>0$. It follows directly that for any $R,r>0$ with $r\le R$, we have $\frac{\Phi(R)}{\Phi(r)}\le C_{\Phi}\left(\frac{R}{r}\right)^{\log_2C_\Phi}$. Throughout this paper, we always assume that $\Phi,\Psi$ are two doubling functions with doubling constants $C_\Phi,C_\Psi$, respectively. Moreover, we assume that there exist $\beta^{(1)}_\Psi,\beta^{(2)}_\Psi>0$ with $\beta^{(1)}_\Psi\le\beta^{(2)}_\Psi$ such that
\begin{equation}\label{eq_beta12}
\frac{1}{C_\Psi}\left(\frac{R}{r}\right)^{\beta^{(1)}_\Psi}\le \frac{\Psi(R)}{\Psi(r)}\le C_\Psi \left(\frac{R}{r}\right)^{\beta^{(2)}_\Psi}\text{ for any }r\le R.
\end{equation}
Indeed, we can take $\beta^{(2)}_\Psi=\log_2C_\Psi$.

Let $(X,d,m)$ be a metric measure space, that is, $(X,d)$ is a locally compact separable metric space and $m$ is a positive Radon measure on $X$ with full support. For any $x\in X$ and any $r>0$, denote $B(x,r)=\{y\in X:d(x,y)<r\}$, $\overline{B}(x,r)=\{y\in X:d(x,y)\le r\}$, and $V(x,r)=m(B(x,r))$. If $B=B(x,r)$, then denote $\delta B=B(x,\delta r)$ for any $\delta>0$. Let $\mathrm{diam}(X)=\sup\{d(x,y):x,y\in X\}$ be the diameter of $(X,d)$. Let $\mathcal{B}(X)$ be the family of all Borel measurable subsets of $X$. Let $C(X)$ be the family of all continuous functions on $X$. Let $C_c(X)$ be the family of all continuous functions on $X$ with compact support. Denote $\dashint_A=\frac{1}{m(A)}\int_A$ and $u_A=\dashint_Au\mathrm{d} m$ for any measurable set $A$ with $m(A)\in(0,+\infty)$ and any function $u$ such that the integral $\int_Au\mathrm{d} m$ is well-defined.

We say that the chain condition \ref{eq_CC} holds if there exists $C_{cc}>0$ such that for any $x,y\in X$ and any positive integer $n$, there exists a sequence $\{x_k:0\le k\le n\}$ of points in $X$ with $x_0=x$ and $x_n=y$ such that
\begin{equation*}\label{eq_CC}\tag*{CC}
d(x_k,x_{k-1})\le C_{cc} \frac{d(x,y)}{n}\text{ for any }k=1,\ldots,n.
\end{equation*}
A simple but useful consequence of \ref{eq_CC} is that for any open interval $(a,b)\subseteq(0,\mathrm{diam}(X))$, there exist $x,y\in X$ such that $d(x,y)\in(a,b)$.

We say that the volume doubling condition \ref{eq_VD} holds if there exists $C_{VD}>0$ such that
\begin{equation*}\label{eq_VD}\tag*{VD}
V(x,2r)\le C_{VD}V(x,r)\text{ for any }x\in X,r>0.
\end{equation*}

We say that the volume regular condition \ref{eq_VPhi} holds if there exists $C_{VR}>0$ such that
\begin{equation*}\label{eq_VPhi}\tag*{V($\Phi$)}
\frac{1}{C_{VR}}\Phi(r)\le V(x,r)\le C_{VR}\Phi(r)\text{ for any }x\in X,r\in(0,\mathrm{diam}(X)).
\end{equation*}
For $d_h>0$, we say that the Ahlfors regular condition \hypertarget{eq_Valpha}{V($d_h$)} holds if \ref{eq_VPhi} holds with $\Phi:r\mapsto r^{d_h}$.

We say that $(\mathcal{E},\mathcal{F})$ is a $p$-energy on $(X,d,m)$ if $\mathcal{F}$ is a dense subspace of $L^p(X;m)$ and $\mathcal{E}:\mathcal{F}\to[0,+\infty) $ satisfies the following conditions.

\begin{enumerate}[label=(\arabic*)]
\item $\mathcal{E}^{1/p}$ is a semi-norm on $\mathcal{F}$, that is, for any $f,g\in\mathcal{F}$ and any $c\in \mathbb{R}$, we have $\mathcal{E}(f)\ge0$, $\mathcal{E}(cf)^{1/p}=|c|\mathcal{E}(f)^{1/p}$, and $\mathcal{E}(f+g)^{1/p}\le\mathcal{E}(f)^{1/p}+\mathcal{E}(g)^{1/p}$.
\item (Closed property) $(\mathcal{F},\mathcal{E}(\cdot)^{1/p}+\lVert \cdot\rVert_{L^p(X;m)})$ is a Banach space.
\item (Markovian property) For any $\varphi\in C(\mathbb{R})$ with $\varphi(0)=0$ and $|\varphi(t)-\varphi(s)|\le|t-s|$ for any $t,s\in \mathbb{R}$, for any $f\in\mathcal{F}$, we have $\varphi(f)\in\mathcal{F}$ and $\mathcal{E}(\varphi(f))\le\mathcal{E}(f)$.
\item (Regular property) $\mathcal{F}\cap C_c(X)$ is uniformly dense in $C_c(X)$ and $(\mathcal{E}(\cdot)^{1/p}+\lVert \cdot\rVert_{L^p(X;m)})$-dense in $\mathcal{F}$.
\item (Strongly local property) For any $f,g\in\mathcal{F}$ with compact support and $g$ constant in an open neighborhood of $\mathrm{supp}(f)$, we have $\mathcal{E}(f+g)=\mathcal{E}(f)+\mathcal{E}(g)$.
\item ($p$-Clarkson's inequality) For any $f,g\in\mathcal{F}$, we have
\begin{equation*}\label{eq_Cla}\tag*{Cla}
\begin{cases}
\mathcal{E}(f+g)+\mathcal{E}(f-g)\ge2 \left(\mathcal{E}(f)^{\frac{1}{p-1}}+\mathcal{E}(g)^{\frac{1}{p-1}}\right)^{p-1}&\text{if }p\in(1,2],\\
\mathcal{E}(f+g)+\mathcal{E}(f-g)\le2 \left(\mathcal{E}(f)^{\frac{1}{p-1}}+\mathcal{E}(g)^{\frac{1}{p-1}}\right)^{p-1}&\text{if }p\in[2,+\infty).\\
\end{cases}
\end{equation*}
\end{enumerate}
For $p=2$, the above definition coincides with the definition of strongly local regular Dirichlet forms.

By \cite[Theorem 2.4]{Sas26}, a $p$-energy $(\mathcal{E},\mathcal{F})$ corresponds to a (canonical) $p$-energy measure $\Gamma:\mathcal{F}\times\mathcal{B}(X)\to[0,+\infty)$, $(f,A)\mapsto\Gamma(f)(A)$ satisfying the following conditions.
\begin{enumerate}[label=(\arabic*)]
\item For any $f\in\mathcal{F}$, $\Gamma(f)(\cdot)$ is a positive Radon measure on $X$ with $\Gamma(f)(X)=\mathcal{E}(f)$.
\item For any $A\in\mathcal{B}(X)$, $\Gamma(\cdot)(A)^{1/p}$ is a semi-norm on $\mathcal{F}$.
\item For any $f,g\in\mathcal{F}\cap C_c(X)$, $A\in\mathcal{B}(X)$, if $f-g$ is constant on $A$, then $\Gamma(f)(A)=\Gamma(g)(A)$.
\item ($p$-Clarkson's inequality) For any $f,g\in\mathcal{F}$ and any $A\in\mathcal{B}(X)$, we have
\begin{equation*}
\begin{cases}
\Gamma(f+g)(A)+\Gamma(f-g)(A)\ge2 \left(\Gamma(f)(A)^{\frac{1}{p-1}}+\Gamma(g)(A)^{\frac{1}{p-1}}\right)^{p-1}&\text{if }p\in(1,2],\\
\Gamma(f+g)(A)+\Gamma(f-g)(A)\le2 \left(\Gamma(f)(A)^{\frac{1}{p-1}}+\Gamma(g)(A)^{\frac{1}{p-1}}\right)^{p-1}&\text{if }p\in[2,+\infty).\\
\end{cases}
\end{equation*}
\item (Chain rule) For any $f\in\mathcal{F}\cap C_c(X)$, for any piecewise $C^1$ function $\varphi:\mathbb{R}\to \mathbb{R}$, we have $\mathrm{d}\Gamma(\varphi(f))=|\varphi'(f)|^p\mathrm{d}\Gamma(f)$.
\end{enumerate}
For $p=2$, the above definition coincides with the definition of energy measures with respect to strongly local regular Dirichlet forms. Using the chain rule, we have the following conditions.
\begin{enumerate}[label=(\arabic*)]
\item (Strong sub-additivity) For any $f,g\in\mathcal{F}$, we have $f\vee g, f\wedge g\in\mathcal{F}$ and
\begin{equation*}\label{eq_SubAdd}\tag*{SubAdd}
\mathcal{E}(f\vee g)+\mathcal{E}(f\wedge g)\le\mathcal{E}(f)+\mathcal{E}(g).
\end{equation*}
\item ($\mathcal{F}\cap L^\infty(X;m)$ is an algebra) For any $f,g\in\mathcal{F}\cap L^\infty(X;m)$, we have $fg\in\mathcal{F}$ and
\begin{equation*}\label{eq_Alg}\tag*{Alg}
\mathcal{E}(fg)^{{1}/{p}}\le C_p\max\left\{\lVert f\rVert_{L^\infty(X;m)},\lVert g\rVert_{L^\infty(X;m)}\right\}\left(\mathcal{E}(g)^{1/p}+\mathcal{E}(f)^{1/p}\right),
\end{equation*}
where $C_p>0$ is some constant depending only on $p$.
\end{enumerate}

We say that the Poincar\'e inequality \ref{eq_PI} holds if there exist $C_{PI}>0$, $A_{PI}\ge1$ such that for any ball $B$ with radius $r\in(0,\mathrm{diam}(X)/A_{PI})$, for any $f\in\mathcal{F}$, we have
\begin{equation*}\label{eq_PI}\tag*{PI($\Psi$)}
\int_B\lvert f-f_B\rvert^p\mathrm{d} m\le C_{PI}\Psi(r)\int_{A_{PI}B}\mathrm{d}\Gamma(f).
\end{equation*}
For $\beta_p>0$, we say that the Poincar\'e inequality \hypertarget{eq_PIbeta}{PI($\beta_p$)} holds if \ref{eq_PI} holds with $\Psi:r\mapsto r^{\beta_p}$.

Let $U,V$ be two open subsets of $X$ satisfying $U\subseteq \overline{U}\subseteq V$. We say that $\phi\in\mathcal{F}$ is a cutoff function for $U\subseteq V$ if $0\le\phi\le1$ in $X$, $\phi=1$ in an open neighborhood of $\overline{U}$ and $\mathrm{supp}(\phi)\subseteq V$, where $\mathrm{supp}(f)$ refers to the support of the measure of $|f|\mathrm{d} m$ for any given function $f$.

We say that the cutoff Sobolev inequality \ref{eq_CS} holds if there exist $C_{1},C_{2}>0$, $A_{S}>1$ such that for any ball $B(x,r)$, there exists a cutoff function $\phi\in\mathcal{F}$ for $B(x,r)\subseteq B(x,A_Sr)$ such that for any $f\in\mathcal{F}$, we have
\begin{equation*}\label{eq_CS}\tag*{CS($\Psi$)}
\int_{B(x,A_{S}r)}|\widetilde{f}|^p\mathrm{d}\Gamma(\phi)\le C_{1}\int_{B(x,A_{S}r)}\mathrm{d}\Gamma(f)+\frac{C_{2}}{\Psi(r)}\int_{B(x,A_{S}r)}|f|^p\mathrm{d} m,
\end{equation*}
where $\widetilde{f}$ is a quasi-continuous modification of $f$, such that $\widetilde{f}$ is uniquely determined $\Gamma(\phi)$-a.e. in $X$. In Section \ref{sec_quasi}, we will provide some necessary results from potential theory in the $p$-energy setting, in parallel with the Dirichlet form framework. For $\beta_p>0$, we say that the cutoff Sobolev inequality \hypertarget{eq_CSbeta}{CS($\beta_p$)} holds if \ref{eq_CS} holds with $\Psi:r\mapsto r^{\beta_p}$.

Let $A_1,A_2\in\mathcal{B}(X)$. We define the capacity between $A_1$ and $A_2$ as
\begin{align*}
&\mathrm{cap}(A_1,A_2)=\inf\left\{\mathcal{E}(\varphi):\varphi\in\mathcal{F},
\begin{array}{l}
\varphi=1\text{ in an open neighborhood of }A_1,\\
\varphi=0\text{ in an open neighborhood of }A_2
\end{array}
\right\},
\end{align*}

We say that the two-sided capacity bounds \hypertarget{eq_cap}{$\text{cap}(\Psi)$} hold if both the capacity upper bound \ref{eq_ucap} and the capacity lower bound \ref{eq_lcap} hold as follows. There exist $C_{cap}>0$, $A_{cap}>1$ such that for any ball $B(x,r)$, we have
\begin{align*}
\mathrm{cap}\left(B(x,r),X\backslash B(x,A_{cap}r)\right)&\le C_{cap} \frac{V(x,r)}{\Psi(r)},\label{eq_ucap}\tag*{$\text{cap}(\Psi)_{\le}$}\\
\mathrm{cap}\left(B(x,r),X\backslash B(x,A_{cap}r)\right)&\ge \frac{1}{C_{cap}} \frac{V(x,r)}{\Psi(r)}.\label{eq_lcap}\tag*{$\text{cap}(\Psi)_{\ge}$}
\end{align*}
For $\beta_p>0$, we say that \hypertarget{eq_capbeta}{$\text{cap}(\beta_p)$} (resp. \hypertarget{eq_ucapbeta}{$\text{cap}(\beta_p)_{\le}$}, \hypertarget{eq_lcapbeta}{$\text{cap}(\beta_p)_{\ge}$}) holds if \hyperlink{eq_cap}{$\text{cap}(\Psi)$} (resp. \ref{eq_ucap}, \ref{eq_lcap}) holds with $\Psi:r\mapsto r^{\beta_p}$. Under \ref{eq_VD}, by taking $f\equiv1$ in $B(x,A_Sr)$, it is easy to see that \ref{eq_CS} (resp. \hyperlink{eq_CSbeta}{$\text{CS}(\beta_p)$}) implies \ref{eq_ucap} (resp. \hyperlink{eq_ucapbeta}{$\text{cap}(\beta_p)_{\le}$}).

Our first result is that certain geometric and functional conditions impose constraints on doubling functions as follows.

\begin{proposition}\label{prop_PhiPsi_bound}
Let $(X,d,m)$ be a metric measure space and $(\mathcal{E},\mathcal{F})$ a $p$-energy with a $p$-energy measure $\Gamma$. Assume that \ref{eq_CC}, \ref{eq_VPhi}, \ref{eq_PI} and \ref{eq_ucap} hold. Then there exists $C>0$ such that for any $R,r\in(0,\mathrm{diam}(X))$ with $r\le R$, we have
\begin{equation*}
\frac{1}{C}\left(\frac{R}{r}\right)^p\le\frac{\Psi(R)}{\Psi(r)}\le C\left(\frac{R}{r}\right)^{p-1}\frac{\Phi(R)}{\Phi(r)}.
\end{equation*}
In particular, assume that \ref{eq_CC}, \hyperlink{eq_Valpha}{V($d_h$)}, \hyperlink{eq_PIbeta}{PI($\beta_p$)} and \hyperlink{eq_ucapbeta}{$\text{cap}(\beta_p)_{\le}$} hold, then
$$p\le\beta_p\le d_h+(p-1),$$
see Figure \ref{fig_adm}.
\begin{figure}[ht]
\begin{center}
\begin{tikzpicture}
\draw[->] (0,0)--(5,0);
\draw[->] (0,0)--(0,7);

\draw[dashed] (1,0)--(1,3);
\draw[dashed] (0,3)--(1,3);
\draw[dashed] (0,2)--(1,3);

\fill[blue] (1,3)--(5,7)--(5,3)--(1,3);

\node[left] at (0,7) {$\beta_p$};
\node[left] at (0,3) {$p$};
\node[left] at (0,2) {$p-1$};
\node[left,below] at (0,0) {$0$};
\node[below] at (1,0) {$1$};
\node[below] at (5,0) {$d_h$};
\end{tikzpicture}
\end{center}
\caption{The admissible region of $(d_h,\beta_p)$}\label{fig_adm}
\end{figure}
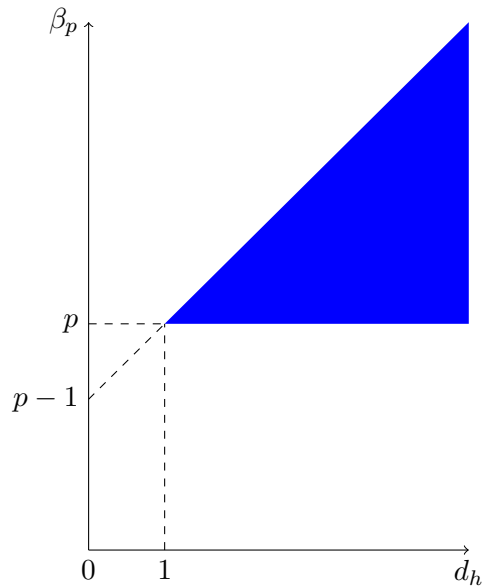
\end{proposition}

\begin{remark}
Our proof will follow arguments similar to those in \cite{Mur20,Mur24}, where the case $p=2$ was considered. As will be clear from the proof, the lower bound require only the assumption of \ref{eq_VD}, rather than \ref{eq_VPhi}. In fact, we will provide two proofs showing that
\begin{equation}\label{eq_lbd1}
\text{\ref{eq_CC} + \ref{eq_VD} + \ref{eq_PI} + \ref{eq_ucap} $\Rightarrow$ the lower bound,}
\end{equation}
and
\begin{equation}\label{eq_lbd2}
\text{\ref{eq_CC} + \ref{eq_VD} + \hyperlink{eq_cap}{$\text{cap}(\Psi)$} $\Rightarrow$ the lower bound.}
\end{equation}
Both implications are of independent interest, even though \ref{eq_PI} already implies \ref{eq_lcap} under \ref{eq_CC} and \ref{eq_VD} (see, for example, \cite[Lemma 5.1 (a)]{BB04}). In particular, the lower bound means that the parameter $\beta^{(1)}_\Psi$ in Equation (\ref{eq_beta12}) satisfies $\beta^{(1)}_\Psi\ge p$.
\end{remark}

Our main result is that the ``converse" of Proposition \ref{prop_PhiPsi_bound} also holds as follows.

\begin{theorem}\label{thm_main}
Assume that there exists $C>0$ such that for any $R,r>0$ with $r\le R$, we have
$$\frac{1}{C}\left(\frac{R}{r}\right)^p\le \frac{\Psi(R)}{\Psi(r)}\le C \left(\frac{R}{r}\right)^{p-1}\frac{\Phi(R)}{\Phi(r)}.$$
Then there exists an unbounded geodesic metric measure space $(X,d,m)$ and a $p$-energy $(\mathcal{E},\mathcal{F})$ with a $p$-energy measure $\Gamma$ such that \ref{eq_VPhi}, \ref{eq_PI} and \ref{eq_CS} hold. In particular, for any $d_h, \beta_p>0$ satisfying
$$p\le \beta_p\le d_h+(p-1),$$
there exists an unbounded geodesic metric measure space $(X,d,m)$ and a $p$-energy $(\mathcal{E},\mathcal{F})$ with a $p$-energy measure $\Gamma$ such that \hyperlink{eq_Valpha}{V($d_h$)}, \hyperlink{eq_PIbeta}{PI($\beta_p$)} and \hyperlink{eq_CSbeta}{$\text{CS}(\beta_p)$} hold.
\end{theorem}

As an application, we consider the critical exponent of Besov spaces as follows. For any $\alpha>0$, we have the following definition of Besov spaces.
\begin{align*}
&B^{p,\alpha}(X,d,m)\\
&=\left\{f\in L^p(X;m):\sup_{r\in(0,\mathrm{diam}(X))}\frac{1}{r^{p\alpha}}\int_X\dashint_{B(x,r)}|f(x)-f(y)|^pm(\mathrm{d} y)m(\mathrm{d} x)<+\infty\right\}.
\end{align*}
Obviously, $B^{p,\alpha}(X,d,m)$ is decreasing in $\alpha$ and may become trivial if $\alpha$ is too large. We define the following critical exponent
$$\alpha_p(X,d,m)=\sup \left\{\alpha>0:B^{p,\alpha}(X,d,m)\text{ contains non-constant functions}\right\}.$$
Let
$$\beta_p(X,d,m)=p\alpha_p(X,d,m).$$

We have a relatively cheap bound for $\beta_p(X,d,m)$ as follows.

\begin{proposition}\label{prop_cri_bound}
Let $(X,d,m)$ be a metric measure space satisfying \ref{eq_CC} and \hyperlink{eq_Valpha}{V($d_h$)}. Then
\begin{enumerate}[label=(\arabic*)]
\item (\cite[Theorem 4.1]{Bau24})
$$p\le\beta_p(X,d,m).$$
\item (\cite[Theorem 4.8]{Shi24a})
$$\beta_p(X,d,m)\le d_h+(p-1).$$
\end{enumerate}
\end{proposition}

\begin{remark}
The bound $\beta_p(X,d,m)\le d_h+p$ was proved in \cite[Theorem 4.2]{Bau24}, and the bound $\beta_p(X,d,m)\le d_h+(p-1)$ was left as a question in \cite[Page 524, Line 3]{Bau24}. The proof of (2) in \cite{Shi24a} proceeds similarly to the proof of \cite[Theorem 4.8 (ii)]{GHL03}, where the case $p=2$ was considered.
\end{remark}

A direct but important consequence of Theorem \ref{thm_main} is that the converse of Proposition \ref{prop_cri_bound} holds as follows.

\begin{corollary}\label{cor_cri_main}
For any $d_h, \beta_p>0$ satisfying
$$p\le \beta_p\le d_h+(p-1),$$
there exists an unbounded geodesic metric measure space $(X,d,m)$ satisfying \hyperlink{eq_Valpha}{V($d_h$)} such that $\beta_p(X,d,m)=\beta_p$.
\end{corollary}

\begin{proof}
Let $(X,d,m)$, $(\mathcal{E},\mathcal{F})$, $\Gamma$ be given by Theorem \ref{thm_main}, on which \hyperlink{eq_Valpha}{V($d_h$)}, \hyperlink{eq_PIbeta}{PI($\beta_p$)} and \hyperlink{eq_ucapbeta}{$\text{cap}(\beta_p)_{\le}$} hold. By \cite[Theorem 4.6]{Shi24a}, the conjunction of \ref{eq_VD}, \hyperlink{eq_PIbeta}{PI($\beta_p$)}, \hyperlink{eq_ucapbeta}{$\text{cap}(\beta_p)_{\le}$} implies $\beta_p(X,d,m)=\beta_p$.
\end{proof}

This paper is organized as follows. In Section \ref{sec_PhiPsi_bound}, we give the proof of Proposition \ref{prop_PhiPsi_bound}. In Section \ref{sec_tree}, we construct an $\mathbb{R}$-tree, introduce a $p$-energy with a $p$-energy measure, and prove the Poincar\'e inequality, the capacity upper bound and the cutoff Sobolev inequality. In Section \ref{sec_Laa_metric}, we introduce a Laakso-type space along with a geodesic metric. In Section \ref{sec_Laa_energy}, we introduce a $p$-energy with a $p$-energy measure on Laakso-type spaces, and prove the Poincar\'e inequality, the capacity upper bound and the cutoff Sobolev inequality. In Section \ref{sec_proof_main}, we give the proof of Theorem \ref{thm_main}. In Section \ref{sec_quasi}, we provide some necessary results from potential theory.

\section{Proof of Proposition \ref{prop_PhiPsi_bound}}\label{sec_PhiPsi_bound}

We start with the proof of the upper bound. We need the following result to pick up approximately evenly spaced points.

\begin{lemma}[{\cite[Lemma 4.10]{GHL03}}]\label{lem_subseq}
Let $\{x_k:0\le k\le n\}$ be a sequence of points in $X$. Let $\rho>0$. Assume that $d(x_0,x_n)>2\rho$ and
$$d(x_k,x_{k-1})<\rho\text{ for any }k=1,\ldots,n,$$
Then there exists a subsequence $\{x_{k_i}:0\le i\le l\}$ such that
\begin{enumerate}[label=(\alph*)]
\item $0=k_0<k_1<\ldots<k_l=n$.
\item $d(x_{k_{i-1}},x_{k_{i}})<5\rho$ for any $i=1,\ldots,l$.
\item $d(x_{k_i},x_{k_j})\ge 2\rho$ for any distinct $i,j=0,1,\ldots,l$.
\end{enumerate}
\end{lemma}

\begin{proof}[Proof of the upper bound]
For notational convenience, we may assume that $\mathrm{diam}(X)=+\infty$. For any $R,r>0$ with $r\le R$. If $R/r\lesssim1$, then since $\Psi$ is doubling and $\Phi$ is increasing, we have
$$\frac{\Psi(R)}{\Psi(r)}\lesssim1\le\left(\frac{R}{r}\right)^{p-1}\frac{\Phi(R)}{\Phi(r)}.$$
Assume that $R/r\gtrsim1$. Let $N$ be the integer satisfying $Nr\le R<(N+1)r$, then $N\gtrsim1$. Take $x,y\in X$ satisfying $Nr<d(x,y)<(N+1)r$. By \ref{eq_ucap}, there exists $u\in\mathcal{F}$ with $u=1$ in $B(x,R/(2A_{cap}))$ and $u=0$ on $X\backslash B(x,R/2)$ such that
$$\mathcal{E}(u)\le 2C_{cap}\frac{V(x,\frac{R}{2A_{cap}})}{\Psi(\frac{R}{2A_{cap}})}\asymp\frac{\Phi(R)}{\Psi(R)}.$$
By \ref{eq_CC}, there exists a sequence $\{y_k:0\le k\le 10N\}$ with $y_0=x$ and $y_{10N}=y$ such that
$$d(y_k,y_{k-1})\le C_{cc}\frac{d(x,y)}{10N}\le C_{cc}\frac{(N+1)r}{10N}<C_{cc}\frac{2Nr}{10N}=\frac{1}{5}C_{cc}r\text{ for any }k=1,\ldots,10N.$$
Since $d(y_0,y_{10N})=d(x,y)>Nr$, by Lemma \ref{lem_subseq}, there exists a subsequence of $\{y_k:0\le k\le 10N\}$, denoted as $\{x_k:0\le k\le M\}$ such that $x_0=x$, $x_M=y$, $d(x_k,x_{k-1})<C_{cc}r$ for any $k=1,\ldots,M$ and $d(x_k,x_l)\ge \frac{2}{5}C_{cc}r$ for any distinct $k,l=0,1,\ldots,M$. Obviously, $M\le 10N$. On the other hand
$$Nr< d(x,y)=d(x_0,x_M)\le\sum_{k=1}^Md(x_k,x_{k-1})\le M C_{cc}r,$$
hence $M\ge N/C_{cc}$, which gives $M\asymp N\asymp R/r$.

For any $k=0,\ldots,M$, let $B_k=B(x_k,r)$, then $u\equiv1$ in $B_0$ and $u\equiv0$ in $B_M$, which gives $u_{B_0}=1$ and $u_{B_M}=0$. For any $k=1,\ldots,M$, by \ref{eq_VPhi}, \ref{eq_PI}, and H\"older's inequality, we have
\begin{align*}
&|u_{B_k}-u_{B_{k-1}}|\le\left(\dashint_{B_k}\dashint_{B_{k-1}}|u(x)-u(y)|^pm(\mathrm{d} y)m(\mathrm{d} x)\right)^{1/p}\\
&\le \frac{1}{(m(B_k)m(B_{k-1}))^{1/p}}\left(\int_{B(x_k,2C_{cc}r)}\int_{{B(x_k,2C_{cc}r)}}|u(x)-u(y)|^pm(\mathrm{d} y)m(\mathrm{d} x)\right)^{1/p}\\
&\lesssim \frac{1}{\Phi(r)^{1/p}}\left(\int_{B(x_k,2C_{cc}r)}|u-u_{B(x_k,2C_{cc}r)}|^p\mathrm{d} m\right)^{1/p}\\
&\lesssim\left(\frac{\Psi(r)}{\Phi(r)}\right)^{1/p}\Gamma(u)(B(x_k,2A_{PI}C_{cc}r))^{1/p},
\end{align*}
hence
\begin{align*}
&1=|u_{B_0}-u_{B_M}|\le\sum_{k=1}^M|u_{B_k}-u_{B_{k-1}}|\\
&\lesssim\left(\frac{\Psi(r)}{\Phi(r)}\right)^{1/p}\sum_{k=1}^M\Gamma(u)(B(x_k,2A_{PI}C_{cc}r))^{1/p}\\
&\le\left(\frac{\Psi(r)}{\Phi(r)}\right)^{1/p}M^{1-{1}/{p}}\left(\sum_{k=1}^M\Gamma(u)(B(x_k,2A_{PI}C_{cc}r))\right)^{1/p}\\
&\asymp\left(\frac{\Psi(r)}{\Phi(r)}\right)^{1/p}\left(\frac{R}{r}\right)^{1-{1}/{p}}\left(\int_X\sum_{k=1}^M1_{B(x_k,2A_{PI}C_{cc}r)}\mathrm{d}\Gamma(u)\right)^{1/p}.
\end{align*}
Since $d(x_k,x_l)\ge \frac{2}{5}C_{cc}r$ for any distinct $k,l=0,1,\ldots,M$, by \ref{eq_VPhi}, there exists some positive integer $K$ depending only on $C_\Phi, A_{PI}, C_{cc}$ such that
$$\sum_{k=1}^M1_{B(x_k,2A_{PI}C_{cc}r)}\le K1_{\cup_{k=1}^MB(x_k,2A_{PI}C_{cc}r)}.$$
Hence
$$1\lesssim\left(\frac{\Psi(r)}{\Phi(r)}\right)^{1/p}\left(\frac{R}{r}\right)^{1-{1}/{p}}\mathcal{E}(u)^{1/p}\lesssim\left(\frac{\Psi(r)}{\Phi(r)}\right)^{1/p}\left(\frac{R}{r}\right)^{1-{1}/{p}}\left(\frac{\Phi(R)}{\Psi(R)}\right)^{1/p},$$
which gives
$$\frac{\Psi(R)}{\Psi(r)}\lesssim \left(\frac{R}{r}\right)^{p-1}\frac{\Phi(R)}{\Phi(r)}.$$
\end{proof}

We provide two proofs of the lower bound via the implications (\ref{eq_lbd1}) and (\ref{eq_lbd2}). To this end, we make several preparations. Let $\varepsilon>0$. We say that $V\subseteq X$ is an $\varepsilon$-net if for any distinct $x,y\in V$, we have $d(x,y)\ge\varepsilon$, and for any $z\in X$, there exists $x\in V$ such that $d(x,z)<\varepsilon$. Since $(X,d)$ is separable, all $\varepsilon$-nets are countable.

We have the following partition of unity with controlled energy.

\begin{lemma}\label{lem_cutoff}
Assume that \ref{eq_VD} and \ref{eq_ucap} hold. Then we have the following controlled cutoff condition. There exists $C_{cut}>0$ depending only on $p, C_{\Psi}, C_{VD}, C_{cap}, A_{cap}$ such that for any $\varepsilon\in(0,\mathrm{diam}(X))$, for any $\varepsilon$-net $V$, there exists a family of functions $\{\psi_z\in\mathcal{F}:z\in V\}$ satisfying the following conditions.
\begin{enumerate}[label=(CO\arabic*)]
\item\label{item_COspt} For any $z\in V$, $0\le\psi_z\le 1$ in $X$, $\psi_z=1$ in $B(z,\varepsilon/4)$, and $\psi_z=0$ on $X\backslash B(z,5\varepsilon/4)$.
\item\label{item_COunit} $\sum_{z\in V}\psi_z=1$.
\item\label{item_COenergy} For any $z\in V$, $\mathcal{E}(\psi_z)\le C_{cut}\frac{V(z,\varepsilon)}{\Psi(\varepsilon)}$.
\end{enumerate}
\end{lemma}

\begin{remark}
By \ref{eq_VD}, in any bounded open set, the summation in \ref{item_COunit} is indeed a finite summation.
\end{remark}

The proof is essentially the same as the proof of \cite[Lemma 2.5]{Mur20} for $p=2$. We give the proof here for completeness.

\begin{proof}
For any $z\in V$, let
$$R_z=\{x\in X:d(x,z)=d(x,V)\},$$
then $B(z,\varepsilon/2)\subseteq R_z\subseteq B(z,\varepsilon)$ and $\{R_z:z\in V\}$ is a cover of $X$. Let $N_z$ be an $(\varepsilon/(4A_{cap}))$-net of $R_z$. By \ref{eq_VD}, there exists some positive integer $M_1$ depending only on $C_{VD}, A_{cap}$ such that $\# N_z\le M_1$. For any $w\in N_z$, by \ref{eq_VD} and \ref{eq_ucap}, there exists $\rho_w\in\mathcal{F}$ with $0\le\rho_w\le1$ in $X$, $\rho_w=1$ in $B(w,\varepsilon/(4A_{cap}))$, $\rho_w=0$ on $X\backslash B(w,\varepsilon/4)$ such that
$$\mathcal{E}(\rho_w)\le 2C_{cap}\frac{V(w,\varepsilon/(4A_{cap}))}{\Psi(\varepsilon/(4A_{cap}))}\le C_1 \frac{V(z,\varepsilon)}{\Psi(\varepsilon)},$$
where $C_1>0$ depends only on $C_{\Psi}, C_{VD}, C_{cap}, A_{cap}$.

For any $z\in V$, let $\varphi_z=\max_{w\in N_z}\rho_w$, then $0\le\varphi_z\le1$ in $X$, $\varphi_z=1$ on $R_z$, $\varphi_z=0$ on $\{x\in X:d(x,R_z)\ge\varepsilon/4\}$. By \ref{eq_SubAdd}, we have $\varphi_z\in\mathcal{F}$ and
$$\mathcal{E}(\varphi_z)\le C\sum_{w\in N_z}\mathcal{E}(\rho_w)\le CM_1C_1 \frac{V(z,\varepsilon)}{\Psi(\varepsilon)}=C_2\frac{V(z,\varepsilon)}{\Psi(\varepsilon)},$$
where $C>0$ depends only on $p, M_1$, and $C_2=CM_1C_1$. Since $\{R_z:z\in V\}$ is a cover of $X$, we have $\sum_{z\in V}\varphi_z\ge1$ in $X$. By \ref{eq_VD}, there exists some positive integer $M_2$ depending only on $C_{VD}$ such that $\#\{w\in V:\varphi_w\ne0\text{ in }B(z,\frac{5}{4}\varepsilon)\}\le M_2$ for any $z\in V$, which implies that $\sum_{z\in V}\varphi_z\le M_2$ in $X$.

For any $z\in V$, let $\psi_z=\frac{\varphi_z}{\sum_{z\in V}\varphi_z}$. Since $\varphi_z=0$ on $X\backslash B(z,5\varepsilon/4)$, we have
$$\psi_z=\frac{\varphi_z}{\sum_{w\in V:\varphi_w\ne0\text{ in }B(z,\frac{5}{4}\varepsilon)}\varphi_w}.$$
By \cite[Proposition 2.3 (c)]{Shi24a}, we have $\psi_z\in\mathcal{F}$, and there exists $C_3>0$ depending only on $p, M_2$ such that
$$\mathcal{E}(\psi_z)\le C_3\sum_{w\in V:\varphi_w\ne0\text{ in }B(z,\frac{5}{4}\varepsilon)}\mathcal{E}(\varphi_w)\le C_3M_2C_2 \frac{V(z,\varepsilon)}{\Psi(\varepsilon)}=C_4\frac{V(z,\varepsilon)}{\Psi(\varepsilon)},$$
where $C_4=C_3M_2C_2$, hence we have \ref{item_COenergy}.

For any distinct $z,w\in V$, we claim that
$$B(w,\frac{1}{4}\varepsilon)\subseteq\{x\in X:d(x,R_z)\ge \frac{1}{4}\varepsilon\}.$$
Otherwise there exists $x\in X$ such that $d(x,w)<\varepsilon/4$, $d(x,R_z)<\varepsilon/4$, there exists $y\in R_z$ such that $d(x,y)<\varepsilon/4$, then $d(y,w)\le d(x,y)+d(x,w)<\varepsilon/2$. But $y\in R_z$ implies $d(y,z)=d(y,V)\le d(y,w)<\varepsilon/2$, hence $d(z,w)\le d(y,z)+d(y,w)<\varepsilon$, contradicting to that $V$ is an $\varepsilon$-net. Hence for any $z\in V$, $\sum_{z\in V}\varphi_z=\varphi_z$ in $B(z,\varepsilon/4)$, $\psi_z=\varphi_z/\varphi_z=1$ in $B(z,\varepsilon/4)$, $\psi_z=0$ on $\{x\in X:d(x,R_z)\ge\varepsilon/4\}\supseteq X\backslash B(z,5\varepsilon/4)$, hence we have \ref{item_COspt}.

Finally, \ref{item_COunit} is obvious.
\end{proof}

Let
\begin{align*}
&\mathcal{F}_{loc}=\left\{u:
\begin{array}{l}
\text{for any relatively compact open set }U,\\
\text{there exists }u^\#\in\mathcal{F}\text{ such that }u=u^\#\text{ }m\text{-a.e. in }U
\end{array}
\right\}.
\end{align*}
For any $u\in\mathcal{F}_{loc}$, let $\Gamma(u)|_U=\Gamma(u^\#)|_U$, where $u^\#$, $U$ are given as above, then $\Gamma(u)$ is a well-defined positive Radon measure on $X$.

We have a characterization of \ref{eq_PI} in terms of maximal functions as follows.

\begin{lemma}\label{lem_PI_max}
Assume that \ref{eq_VD} and \ref{eq_PI} hold. Then there exist $C_{max}>0$, $A_{max}>1$ such that for any $u\in\mathcal{F}_{loc}$, for any ball $B$ with radius $R$, for $m$-a.e. $x,y\in A_{max}^{-1}B$, we have
$$|u(x)-u(y)|^p\le C_{max}\Psi(d(x,y))\left(M_R\Gamma(u)(x)+M_R\Gamma(u)(y)\right),$$
where
$$M_R\Gamma(u)(x)=\sup_{r\in(0,R)}\frac{\Gamma(u)(B(x,r))}{V(x,r)}.$$
\end{lemma}

\begin{proof}
The proof is standard using the telescopic technique. Since $u\in\mathcal{F}_{loc}$, we have $u\in L^p_{loc}(X;m)$, $m$-almost all points are Lebesgue points of $u$. Let $x,y\in(4A_{PI})^{-1}B$ be two distinct Lebesgue points of $u$. Denote $D=d(x,y)\in(0,R/(2A_{PI}))$. For any $n\ge0$, let $B_n=B(x,2^{-n}D)$, then by \ref{eq_VD}, \ref{eq_PI}, and H\"older's inequality, we have
\begin{align*}
&|u_{B_n}-u_{B_{n+1}}|\lesssim\left(\dashint_{B_{n}}|u-u_{B_n}|^p\mathrm{d} m\right)^{1/p}\lesssim\left(\frac{\Psi(2^{-n}D)}{m(B_n)}\int_{A_{PI}B_n}\mathrm{d}\Gamma(u)\right)^{1/p}\\
&\lesssim \left(\Psi(2^{-n}D)\frac{\Gamma(u)(A_{PI}B_n)}{m(A_{PI}B_n)}\right)^{1/p}\le\left(\Psi(2^{-n}D)M_R\Gamma(u)(x)\right)^{1/p},
\end{align*}
hence
\begin{align*}
&|u(x)-u_{B(x,D)}|=\lim_{n\to+\infty}|u_{B_n}-u_{B_0}|\le\sum_{n=0}^\infty|u_{B_n}-u_{B_{n+1}}|\\
&\lesssim\sum_{n=0}^\infty\left(\Psi(2^{-n}D)M_R\Gamma(u)(x)\right)^{1/p}\lesssim\left(\Psi(D)M_R\Gamma(u)(x)\right)^{1/p},
\end{align*}
where the last inequality follows from Equation (\ref{eq_beta12}). Similarly, we have
$$|u(y)-u_{B(y,D)}|\lesssim\left(\Psi(D)M_R\Gamma(u)(y)\right)^{1/p}.$$
Moreover
\begin{align*}
&|u_{B(x,D)}-u_{B(y,D)}|\le\left(\dashint_{B(x,D)}\dashint_{B(y,D)}|u(z_1)-u(z_2)|^pm(\mathrm{d} z_1)m(\mathrm{d} z_2)\right)^{1/p}\\
&\lesssim\left(\dashint_{B(x,2D)}\dashint_{B(x,2D)}|u(z_1)-u(z_2)|^pm(\mathrm{d} z_1)m(\mathrm{d} z_2)\right)^{1/p}\\
&\lesssim\left(\frac{1}{V(x,2D)}\int_{B(x,2D)}|u-u_{B(x,2D)}|^p\mathrm{d} m\right)^{1/p}\lesssim\left(\frac{\Psi(2D)}{V(x,2D)}\int_{B(x,2A_{PI}D)}\mathrm{d}\Gamma(u)\right)^{1/p}\\
&\lesssim\left(\Psi(D)\frac{\Gamma(u)(B(x,2A_{PI}D))}{V(x,2A_{PI}D)}\right)^{1/p}\le\left(\Psi(D)M_{R}\Gamma(u)(x)\right)^{1/p}.
\end{align*}
In summary, we have
\begin{align*}
&|u(x)-u(y)|\le|u(x)-u_{B(x,D)}|+|u(y)-u_{B(y,D)}|+|u_{B(x,D)}-u_{B(y,D)}|\\
&\lesssim\Psi(D)^{1/p}\left(M_{R}\Gamma(u)(x)^{1/p}+M_{R}\Gamma(u)(y)^{1/p}\right),
\end{align*}
which implies
$$|u(x)-u(y)|^p\lesssim\Psi(D)\left(M_{R}\Gamma(u)(x)+M_{R}\Gamma(u)(y)\right).$$
\end{proof}

Let $\varepsilon>0$ and $x,y\in X$. We say that $\{x_k:0\le k\le N\}$ is an $\varepsilon$-chain between $x$ and $y$ if $x_0=x$, $x_N=y$, and $d(x_k,x_{k-1})<\varepsilon$ for any $k=1,\ldots,N$. Let
$$N_\varepsilon(x,y)=\inf\left\{N:\{x_k:0\le k\le N\}\text{ is an }\varepsilon\text{-chain between }x\text{ and }y\right\},$$
where $\inf\emptyset=+\infty$, see \cite[Section 6]{GT12} for more details. By the triangle inequality, we have $N_\varepsilon(x,y)\ge\frac{d(x,y)}{\varepsilon}$. Assume \ref{eq_CC}, then $N_\varepsilon(x,y)\le C_{cc}\frac{d(x,y)}{\varepsilon}+1$, hence $N_\varepsilon(x,y)\asymp \frac{d(x,y)}{\varepsilon}$.

\begin{proof}[Proof of the lower bound through the implication (\ref{eq_lbd1})]
For notational convenience, we may assume that $\mathrm{diam}(X)=+\infty$. For any $R,\varepsilon>0$ with $\varepsilon\le R$. If $R/\varepsilon\lesssim1$, then since $\Psi$ is increasing, we have
$$\left(\frac{R}{\varepsilon}\right)^p\lesssim1\le \frac{\Psi(R)}{\Psi(\varepsilon)}.$$
We show that if $R/\varepsilon\gtrsim1$, then
$$\frac{\Psi(R)}{\Psi(\varepsilon)}\gtrsim \left(\frac{R}{\varepsilon}\right)^p.$$

Take $\overline{x},\overline{y}\in X$ with $R<d(\overline{x},\overline{y})< R+\varepsilon$, take an $\varepsilon$-net $V$ with $\{\overline{x},\overline{y}\}\in V$, let $\{\psi_z:z\in V\}$ be the family of functions given by Lemma \ref{lem_cutoff}, let
$$u=\sum_{z\in V}N_\varepsilon(\overline{x},z)\psi_z.$$
By \ref{item_COspt}, we have $u\equiv u(\overline{x})=N_\varepsilon(\overline{x},\overline{x})=0$ in $B(\overline{x},\varepsilon/4)$ and $u\equiv u(\overline{y})=N_\varepsilon(\overline{x},\overline{y})\asymp \frac{d(\overline{x},\overline{y})}{\varepsilon}\asymp \frac{R}{\varepsilon}$ in $B(\overline{y},\varepsilon/4)$. By \ref{eq_VD}, the above summation is locally a finite summation of functions in $\mathcal{F}$, hence $u\in\mathcal{F}_{loc}$.

Firstly, we show that for any $\overline{z}\in V$
$$\Gamma(u)(B(\overline{z},\frac{5}{4}\varepsilon))\lesssim \frac{V(\overline{z},\varepsilon)}{\Psi(\varepsilon)}.$$
Let $N_{\overline{z}}=\{z\in V:B(z,\frac{5}{4}\varepsilon)\cap B(\overline{z},\frac{5}{4}\varepsilon)\ne\emptyset\}$. By \ref{eq_VD}, there exists some positive integer $N$ depending only on $C_{VD}$ such that $\# N_{\overline{z}}\le N$. By \ref{item_COunit}, we have
$$u=\sum_{z\in V}N_\varepsilon(\overline{x},z)\psi_z=\sum_{z\in V}\left( N_\varepsilon(\overline{x},z)-N_\varepsilon(\overline{x},\overline{z})\right)\psi_z+N_\varepsilon(\overline{x},\overline{z}).$$
In particular, in $B(\overline{z},\frac{5}{4}\varepsilon)$, by \ref{item_COspt}, we have
$$u=\sum_{z\in N_{\overline{z}}}\left( N_\varepsilon(\overline{x},z)-N_\varepsilon(\overline{x},\overline{z})\right)\psi_z+N_\varepsilon(\overline{x},\overline{z}),$$
hence
\begin{align*}
&\Gamma(u)(B(\overline{z},\frac{5}{4}\varepsilon))=\Gamma\left(\sum_{z\in N_{\overline{z}}}\left( N_\varepsilon(\overline{x},z)-N_\varepsilon(\overline{x},\overline{z})\right)\psi_z+N_\varepsilon(\overline{x},\overline{z})\right)(B(\overline{z},\frac{5}{4}\varepsilon))\\
&=\Gamma\left(\sum_{z\in N_{\overline{z}}}\left( N_\varepsilon(\overline{x},z)-N_\varepsilon(\overline{x},\overline{z})\right)\psi_z\right)(B(\overline{z},\frac{5}{4}\varepsilon))\\
&\le C_1\sum_{z\in N_{\overline{z}}}\Gamma\left(\left( N_\varepsilon(\overline{x},z)-N_\varepsilon(\overline{x},\overline{z})\right)\psi_z\right)(B(\overline{z},\frac{5}{4}\varepsilon))\\
&\le C_1\sum_{z\in N_{\overline{z}}}|N_\varepsilon(\overline{x},z)-N_\varepsilon(\overline{x},\overline{z})|^p\mathcal{E}(\psi_z),
\end{align*}
where $C_1>0$ depends only on $p,N$. For any $z\in N_{\overline{z}}$, we have $d(z,\overline{z})<\frac{5}{2}\varepsilon$, then by \ref{eq_CC}, we have $|N_\varepsilon(\overline{x},z)-N_\varepsilon(\overline{x},\overline{z})|\le \lceil {\frac{5}{2}C_{cc}}\rceil$. Combining this with \ref{item_COenergy} and \ref{eq_VD}, we have
\begin{align*}
&\Gamma(u)(B(\overline{z},\frac{5}{4}\varepsilon))\le C_1\sum_{z\in N_{\overline{z}}}\lceil{\frac{5}{2}C_{cc}} \rceil^pC_{cut}\frac{V(z,\varepsilon)}{\Psi(\varepsilon)}\\
&\le C_1\sum_{z\in N_{\overline{z}}}\lceil {\frac{5}{2}C_{cc}}\rceil^pC_{cut}C_2\frac{V(\overline{z},\varepsilon)}{\Psi(\varepsilon)}\le C_3\frac{V(\overline{z},\varepsilon)}{\Psi(\varepsilon)},
\end{align*}
where $C_2>0$ depends only on $C_{VD}$, and $C_3=\lceil {\frac{5}{2}C_{cc}}\rceil^pNC_1C_2C_{cut}$.

Secondly, we show that for any $\overline{z}\in V$ and any $r>0$
$$\Gamma(u)(B(\overline{z},r))\lesssim \frac{V(\overline{z},r)}{\Psi(\varepsilon)}.$$
If $r\le\varepsilon/4$, then in $B(\overline{z},r)$, we have $u\equiv N_\varepsilon(\overline{x},\overline{z})$, hence $\Gamma(u)(B(\overline{z},r))=0$, the result is trivial. We may assume that $r>\varepsilon/4$, then
\begin{align*}
&\Gamma(u)(B(\overline{z},r))\le\sum_{z\in V:B(z,\frac{5}{4}\varepsilon)\cap B(\overline{z},r)\ne\emptyset}\Gamma(u)(B(z,\frac{5}{4}\varepsilon))\\
&\le \sum_{z\in V:B(z,\frac{5}{4}\varepsilon)\cap B(\overline{z},r)\ne\emptyset}C_3 \frac{V(z,\varepsilon)}{\Psi(\varepsilon)}=\frac{C_3}{\Psi(\varepsilon)}\int_X\sum_{z\in V:B(z,\frac{5}{4}\varepsilon)\cap B(\overline{z},r)\ne\emptyset}1_{B(z,\varepsilon)}\mathrm{d} m.
\end{align*}
By \ref{eq_VD}, there exists some positive integer $M$ depending only on $C_{VD}$ such that
$$\sum_{z\in V:B(z,\frac{5}{4}\varepsilon)\cap B(\overline{z},r)\ne\emptyset}1_{B(z,\varepsilon)}\le M1_{\cup_{z\in V:B(z,\frac{5}{4}\varepsilon)\cap B(\overline{z},r)\ne\emptyset}B(z,\varepsilon)},$$
hence
\begin{align*}
&\Gamma(u)(B(\overline{z},r))\le\frac{C_3}{\Psi(\varepsilon)}Mm\left({\bigcup_{z\in V:B(z,\frac{5}{4}\varepsilon)\cap B(\overline{z},r)\ne\emptyset}B(z,\varepsilon)}\right)\\
&\le \frac{C_3M}{\Psi(\varepsilon)}m(B(\overline{z},r+\frac{5}{2}\varepsilon))\le \frac{C_3M}{\Psi(\varepsilon)}m(B(\overline{z},11r))\le\frac{C_4}{\Psi(\varepsilon)}V(\overline{z},r),
\end{align*}
where $C_4=C_3C_{VD}^4M$.

Finally, in Lemma \ref{lem_PI_max}, let $B=B(\overline{x},4A_{max}R)$, then $\overline{x},\overline{y}\in A_{max}^{-1}B=B(\overline{x},4R)$, since $B(\overline{x},\frac{\varepsilon}{4})\cup B(\overline{y},\frac{\varepsilon}{4})\subseteq B(\overline{x},4R)$, and $u\equiv u(\overline{x})=N_\varepsilon(\overline{x},\overline{x})=0$ in $B(\overline{x},\frac{\varepsilon}{4})$, $u\equiv u(\overline{y})=N_\varepsilon(\overline{x},\overline{y})\asymp \frac{R}{\varepsilon}$ in $B(\overline{y},\frac{\varepsilon}{4})$, we may assume that $\overline{x}, \overline{y}$ are Lebesgue points of $u$, then
\begin{align*}
|u(\overline{x})-u(\overline{y})|^p\le C_{max}\Psi(d(\overline{x},\overline{y}))\left(M_{4A_{max}R}\Gamma(u)(\overline{x})+M_{4A_{max}R}\Gamma(u)(\overline{y})\right),
\end{align*}
where
\begin{align*}
&M_{4A_{max}R}\Gamma(u)(\overline{x})=\sup_{r\in(0,4A_{max}R)}\frac{\Gamma(u)(B(\overline{x},r))}{V(\overline{x},r)}\le \frac{C_4}{\Psi(\varepsilon)},\\
&M_{4A_{max}R}\Gamma(u)(\overline{y})=\sup_{r\in(0,4A_{max}R)}\frac{\Gamma(u)(B(\overline{y},r))}{V(\overline{y},r)}\le \frac{C_4}{\Psi(\varepsilon)}.
\end{align*}
Hence
$$\left(\frac{R}{\varepsilon}\right)^p\asymp N_\varepsilon(\overline{x},\overline{y})^p=|u(\overline{x})-u(\overline{y})|^p\le C_{max}\Psi(d(\overline{x},\overline{y}))\frac{2C_4}{\Psi(\varepsilon)}\asymp\frac{\Psi(R)}{\Psi(\varepsilon)},$$
which gives
$$\left(\frac{R}{\varepsilon}\right)^p\lesssim \frac{\Psi(R)}{\Psi(\varepsilon)}.$$
\end{proof}

We need the following capacity upper bound.

\begin{lemma}\label{lem_ucapball}
Assume that \ref{eq_VD} and \ref{eq_ucap} hold. Then there exists $C>0$ such that for any $x\in X$ and any $R,r>0$, we have
$$\mathrm{cap}(B(x,R),X\backslash B(x,R+r))\le C \frac{m \left(B(x,R+r)\backslash \overline{B(x,R)}\right)}{\Psi(r)}.$$
\end{lemma}

\begin{proof}
Let $C_1=C_{cap}$, $A=A_{cap}$ be the constants in \ref{eq_ucap}. Let $L=4A+4$ and let $V$ be an $\frac{r}{L}$-net. For any $v\in V$, by \ref{eq_ucap}, there exists a cutoff function $\phi_v\in \mathcal{F}$ for $B(v,\frac{r}{L})\subseteq B(v,\frac{Ar}{L})$ such that
$$\mathcal{E}(\phi_v)=\int_{B(v,\frac{Ar}{L})}\mathrm{d}\Gamma(\phi_v)\le 2C_1\frac{V(v,\frac{r}{L})}{\Psi(\frac{r}{L})}.$$
Let ${V}_1=V\cap B(x,R+\frac{r}{2})$, then by \ref{eq_VD}, we have $\#V_1<+\infty$. Let
$$\phi=\max_{v\in{V}_1}\phi_v,$$
then by \ref{eq_SubAdd}, we have $\phi\in \mathcal{F}$. Moreover, $0\le\phi\le1$ in $X$, $\phi=1$ in $B(x,R+\frac{r}{2}-\frac{r}{L})\supseteq\overline{B(x,R)}$, $\mathrm{supp}(\phi)\subseteq B(x,R+\frac{r}{2}+\frac{Ar}{L})\subseteq B(x,R+r)$, hence $\phi\in \mathcal{F}$ is a cutoff function for $B(x,R)\subseteq B(x,R+r)$ and
$$\mathrm{supp}(\Gamma(\phi))\subseteq B(x,R+\frac{r}{2}+\frac{Ar}{L})\backslash{B(x,R+\frac{r}{2}-\frac{r}{L})}.$$
Let $W={V}_1\backslash B(x,R+\frac{r}{2}-\frac{r}{L}-\frac{Ar}{L})$, then
$$\phi=\max_{v\in W}\phi_v\text{ in }B(x,R+\frac{r}{2}+\frac{Ar}{L})\backslash{B(x,R+\frac{r}{2}-\frac{r}{L})},$$
hence
\begin{align*}
&\mathcal{E}(\phi)=\int_{B(x,R+\frac{r}{2}+\frac{Ar}{L})\backslash{B(x,R+\frac{r}{2}-\frac{r}{L})}}\mathrm{d}\Gamma(\phi)\overset{(\star)}{\scalebox{2}[1]{$\le$}}\sum_{v\in W}\int_{{B(x,R+\frac{r}{2}+\frac{Ar}{L})\backslash{B(x,R+\frac{r}{2}-\frac{r}{L})}}}\mathrm{d}\Gamma(\phi_v)\\
&\le\sum_{v\in W}\int_{{B(v,\frac{Ar}{L})}}\mathrm{d}\Gamma(\phi_v)\le 2C_1\sum_{v\in W}\frac{V(v,\frac{r}{L})}{\Psi(\frac{r}{L})}=\frac{2C_1}{\Psi(\frac{r}{L})}\int_X\left(\sum_{v\in W}1_{B(v,\frac{r}{L})}\right)\mathrm{d} m,
\end{align*}
where $(\star)$ follows from \ref{eq_SubAdd}. By \ref{eq_VD}, there exists some positive integer $N$ depending only on $C_{VD}$ such that
$$\sum_{v\in W}1_{B(v,\frac{r}{L})}\le N1_{\cup_{v\in W}B(v,\frac{r}{L})},$$
where
$${\bigcup_{v\in W}B(v,\frac{r}{L})}\subseteq B(x,R+\frac{r}{2}+\frac{r}{L})\backslash \overline{B(x,R+\frac{r}{2}-\frac{r}{L}-\frac{Ar}{L}-\frac{r}{L})}\subseteq B(x,R+r)\backslash\overline{B(x,R)}.$$
Therefore, we have
\begin{align*}
&\mathrm{cap}\left(B(x,R),X\backslash B(x,R+r)\right)\le \mathcal{E}(\phi)\le \frac{2C_1N}{\Psi(\frac{r}{L})}m\left(B(x,R+r)\backslash\overline{B(x,R)}\right)\\
&\le \frac{2C_1NC_2}{\Psi({r})}m \left(B(x,R+r)\backslash\overline{B(x,R)}\right)=C \frac{m \left(B(x,R+r)\backslash\overline{B(x,R)}\right)}{\Psi(r)},
\end{align*}
where $C_2>0$ depends only on $C_\Psi,L$, and $C=2C_1NC_2$.
\end{proof}

\begin{proof}[Proof of the lower bound through the implication (\ref{eq_lbd2})]
For notational convenience, we may assume that $\mathrm{diam}(X)=+\infty$. Let $A=A_{cap}$ be the constant in \hyperlink{eq_cap}{$\text{cap}(\Psi)$}. For any $R,r>0$ with $r\le R$, without loss of generality, we may assume that $R/r\gtrsim1$, let $N$ be the integer satisfying $Nr\le (A-1)R<(N+1)r$, then $N\asymp R/r\gtrsim1$. By \ref{eq_VD} and \ref{eq_ucap}, applying Lemma \ref{lem_ucapball}, there exists $C_1>0$ such that for any $n=0,\ldots,N-1$, there exists a cutoff function $\phi_n\in \mathcal{F}$ for $B(x,R+nr)\subseteq B(x,R+(n+1)r)$ such that
$$\mathcal{E}(\phi_n)\le C_1 \frac{m \left(B(x,R+(n+1)r)\backslash \overline{B(x,R+nr)}\right)}{\Psi(r)}.$$
Let $\phi=\frac{1}{N}\sum_{n=0}^{N-1}\phi_n$, then $\phi\in \mathcal{F}$ is a cutoff function for $B(x,R)\subseteq B(x,R+Nr)$, hence also a cutoff function for $B(x,R)\subseteq B(x,AR)$. By \ref{eq_lcap}, we have
\begin{align*}
&\frac{1}{C_{cap}}\frac{V(x,R)}{\Psi(R)}\le \mathrm{cap}\left(B(x,R),X\backslash B(x,AR)\right)\le \mathcal{E}(\phi)\overset{(\star)}{\scalebox{2}[1]{$\le$}}\frac{1}{N^p}\sum_{n=0}^{N-1}\mathcal{E}\left(\phi_n\right)\\
&\le \frac{C_1}{N^p} \sum_{n=0}^{N-1}\frac{m \left(B(x,R+(n+1)r)\backslash \overline{B(x,R+nr)}\right)}{\Psi(r)}\le \frac{C_1}{N^p} \frac{V (x,AR)}{\Psi(r)}\overset{(\diamond)}{\scalebox{2}[1]{$\le$}} \frac{C_1C_2}{N^p} \frac{V(x,R)}{\Psi(r)},
\end{align*}
where $(\star)$ follows from \ref{eq_SubAdd}, $(\diamond)$ follows from \ref{eq_VD}, and $C_2>0$ depends only on $C_{VD}, A$. Therefore, we have
$$\frac{\Psi(R)}{\Psi(r)}\gtrsim N^p\asymp \left(\frac{R}{r}\right)^p.$$
\end{proof}

\section{Construction of \texorpdfstring{$\mathbb{R}$}{R}-trees}\label{sec_tree}

In this section, we construct an $\mathbb{R}$-tree, introduce a $p$-energy with a $p$-energy measure, and prove the Poincar\'e inequality, the capacity upper bound and the cutoff Sobolev inequality.

Fix a branching function $\mathbf{b}:\mathbb{Z}\to\mathbb{N}$, which satisfies
$$2\le\inf_{k\in\mathbb{Z}}\mathbf{b}(k)\le\sup_{k\in\mathbb{Z}}\mathbf{b}(k)<+\infty,$$
and a gluing function $\mathbf{g}:\mathbb{Z}\to\mathbb{N}$, which satisfies
$$1\le\inf_{k\in\mathbb{Z}}\mathbf{g}(k)\le\sup_{k\in\mathbb{Z}}\mathbf{g}(k)<+\infty.$$
Let
$$\mathcal{U}(\mathbf{g})=\left\{\mathbf{s}:\mathbb{Z}\to\mathbb{Z}|\mathbf{s}(k)\in\{0,1,\ldots,\mathbf{g}(k)-1\}\text{ for any }k\in\mathbb{Z}\text{ and }\lim_{k\to+\infty}\mathbf{s}(k)=0\right\}.$$
For any distinct $\mathbf{s},\mathbf{t}\in\mathcal{U}(\mathbf{g})$, let $\mathbf{d}_{\mathcal{U}(\mathbf{g})}(\mathbf{t},\mathbf{t})=0$ and
$$\mathbf{d}_{\mathcal{U}(\mathbf{g})}(\mathbf{s},\mathbf{t})=2^{\max\{k\in\mathbb{Z}:\mathbf{s}(k)\ne\mathbf{t}(k)\}}.$$
Then $(\mathcal{U}(\mathbf{g}),\mathbf{d}_{\mathcal{U}(\mathbf{g})})$ is an ultrametric space, that is, it is a metric space and the metric satisfies the following inequality.
$$\mathbf{d}_{\mathcal{U}(\mathbf{g})}(\mathbf{s},\mathbf{t})\le\max \left\{\mathbf{d}_{\mathcal{U}(\mathbf{g})}(\mathbf{s},\mathbf{r}),\mathbf{d}_{\mathcal{U}(\mathbf{g})}(\mathbf{r},\mathbf{t})\right\}\text{ for any }\mathbf{r},\mathbf{s},\mathbf{t}\in\mathcal{U}(\mathbf{g}).$$

Let us list some ``unusual" properties of ultrametric spaces as follows, which are easy consequences of the above inequality.
\begin{itemize}
\item Any two balls with the same radius are either disjoint or coincide. More generally, any two balls are either disjoint or one contains the other.
\item Any point inside a ball is its center. Any ball is simultaneously open and closed.
\end{itemize}

For any $\mathbf{s}\in\mathcal{U}(\mathbf{g})$ and any $n\in\mathbb{Z}$, the ball $\{\mathbf{t}\in\mathcal{U}(\mathbf{g}):\mathbf{d}_{\mathcal{U}(\mathbf{g})}(\mathbf{t},\mathbf{s})\le 2^n\}$ coincides with the cylindrical set $\{\mathbf{t}\in\mathcal{U}(\mathbf{g}):\mathbf{t}(k)=\mathbf{s}(k)\text{ for any }k\ge n+1\}$. For any $r>0$, let
$$
V_{\mathbf{g}}(r)=
\begin{cases}
\left(\prod_{k=n}^0\mathbf{g}(k)\right)^{-1}&\text{if }2^{n-1}\le r<2^n,n\le0,\\
1&\text{if }1\le r<2,\\
\prod_{k=1}^{n-1}\mathbf{g}(k)&\text{if }2^{n-1}\le r<2^n,n\ge2.
\end{cases}
$$
Then there exists a unique measure $\mathbf{m}_{\mathcal{U}(\mathbf{g})}$ on $(\mathcal{U}(\mathbf{g}),\mathbf{d}_{\mathcal{U}(\mathbf{g})})$ such that any closed ball with radius $2^n$ has measure $V_{\mathbf{g}}(2^n)$. It is obvious that there exists some positive constant $C$ depending only on $\sup_\mathbb{Z}\mathbf{g}$ such that
$$\frac{1}{C}V_{\mathbf{g}}(r)\le\mathbf{m}_{\mathcal{U}(\mathbf{g})}(B_{\mathcal{U}(\mathbf{g})}(\mathbf{t},r))\le CV_{\mathbf{g}}(r)\text{ for any }\mathbf{t}\in\mathcal{U}(\mathbf{g}),r>0.$$

Similarly, we also have an ultrametric space $(\mathcal{U}(\mathbf{b}),\mathbf{d}_{\mathcal{U}(\mathbf{b})})$ endowed with a measure $\mathbf{m}_{\mathcal{U}(\mathbf{b})}$, under which any closed ball with radius $2^n$ has measure $V_{\mathbf{b}}(2^n)$, and there exists some positive constant $C$ depending only on $\sup_\mathbb{Z}\mathbf{b}$ such that
$$\frac{1}{C}V_{\mathbf{b}}(r)\le\mathbf{m}_{\mathcal{U}(\mathbf{b})}(B_{\mathcal{U}(\mathbf{b})}(\mathbf{t},r))\le CV_{\mathbf{b}}(r)\text{ for any }\mathbf{t}\in\mathcal{U}(\mathbf{b}),r>0.$$

We construct an $\mathbb{R}$-tree $(\mathcal{T}(\mathbf{b}),\mathbf{d}_{\mathcal{T}(\mathbf{b})})$ endowed with a measure $\mathbf{m}_{\mathcal{T}(\mathbf{b})}$, using approximation of finite trees. A detailed and rigorous construction was given in \cite[Subsections 3.1, 3.2, 3.3]{Mur24}. Here we only outline the main idea as follows.

For any $m,n\in\mathbb{Z}$ with $m\le n$, we construct a finite tree $T_{m,n}$ endowed with a metric $\mathbf{d}_{m,n}$ inductively on $n-m$, starting from $n-m=0$. For the sets $T_{m,n}$. Firstly, $T_{n,n}$ is a star graph $K_{1,\mathbf{b}(n)}$. Secondly, $T_{m,n+1}$ is obtained
\begin{enumerate}[label=(TREE\arabic*)]%, ref=(C.\arabic*)]
\item\label{item_tree1} either by gluing $\mathbf{b}(n+1)$ copies of $T_{m,n}$ at one vertex,
\item\label{item_tree2} or by replacing each edge of $T_{m+1,n}$ by a star graph $K_{1,\mathbf{b}(m)}$.
\end{enumerate}

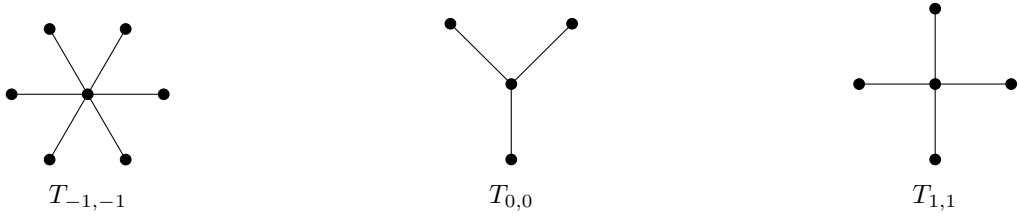
\begin{figure}[htp]
\centering
\begin{subfigure}{0.3\textwidth}
\centering
\begin{tikzpicture}
\foreach \angle in {0,60,120,180,240,300} {
\draw[fill] (\angle:1) circle (2pt);
\draw (0,0) -- (\angle:1);
}
\draw[fill] (0,0) circle (2pt);
\end{tikzpicture}
\caption{$T_{-1,-1}$}
\end{subfigure}
\hfill
\begin{subfigure}{0.3\textwidth}
\centering
\begin{tikzpicture}
\draw[fill] (0,0) circle (2pt);
\draw[fill] (-0.8,0.8) circle (2pt);
\draw[fill] (0.8,0.8) circle (2pt);
\draw[fill] (0,-1) circle (2pt);
\draw (0,-1) -- (0,0);
\draw (0,0) -- (-0.8,0.8);
\draw (0,0) -- (0.8,0.8);
\end{tikzpicture}
\caption{$T_{0,0}$}
\end{subfigure}
\hfill
\begin{subfigure}{0.3\textwidth}
\centering
\begin{tikzpicture}
\draw[fill] (0,0) circle (2pt);
\draw[fill] (0,1) circle (2pt);
\draw[fill] (1,0) circle (2pt);
\draw[fill] (-1,0) circle (2pt);
\draw[fill] (0,-1) circle (2pt);
\draw (0,-1) -- (0,1);
\draw (-1,0) -- (1,0);
\end{tikzpicture}        
\caption{$T_{1,1}$}
\end{subfigure}
\caption{The figures of $T_{-1,-1}$, $T_{0,0}$ and $T_{1,1}$}\label{fig_Tn0}
\end{figure}

\begin{figure}[htp]
\centering
\begin{tikzpicture}
% Define star command
\coordinate (AA) at (0,1.5);
\coordinate (BB) at (-1.5,0);
\coordinate (CC) at (0,-1.5);
\coordinate (DD) at (1.5,0);

\begin{scope}[shift={(AA)}, rotate=60]
\newcommand{\starshape}[2]{
    \begin{scope}[shift={(#1)}, rotate=#2]
        \foreach \angle in {0,60,120,180,240,300} {
            \draw (#1) -- ++(\angle:0.4);
            \draw[fill] ($(#1)+(\angle:0.4)$) circle (1.5pt);
        }
        \draw[fill] (#1) circle (1.5pt);
    \end{scope}
}

% Central point of Y-shape
\coordinate (O) at (0,0);

% Endpoints of Y-shape
\coordinate (A) at (90:0.6);
\coordinate (B) at (210:0.6);
\coordinate (C) at (330:0.6);

% Draw stars at endpoints
\starshape{A}{30}
\starshape{B}{30}
\starshape{C}{30}

\fill[white] (0,0) circle (0.3);

\draw (A)--(0,0);
\draw (B)--(0,0);
\draw (C)--(0,0);

% Mark central intersection point
\draw[fill] (O) circle (1.5pt);
\end{scope}

\begin{scope}[shift={(BB)}, rotate=30]
\newcommand{\starshape}[2]{
    \begin{scope}[shift={(#1)}, rotate=#2]
        \foreach \angle in {0,60,120,180,240,300} {
            \draw (#1) -- ++(\angle:0.4);
            \draw[fill] ($(#1)+(\angle:0.4)$) circle (1.5pt);
        }
        \draw[fill] (#1) circle (1.5pt);
    \end{scope}
}

% Central point of Y-shape
\coordinate (O) at (0,0);

% Endpoints of Y-shape
\coordinate (A) at (90:0.6);
\coordinate (B) at (210:0.6);
\coordinate (C) at (330:0.6);

% Draw stars at endpoints
\starshape{A}{30}
\starshape{B}{30}
\starshape{C}{30}

\fill[white] (0,0) circle (0.3);

\draw (A)--(0,0);
\draw (B)--(0,0);
\draw (C)--(0,0);

% Mark central intersection point
\draw[fill] (O) circle (1.5pt);
\end{scope}

\begin{scope}[shift={(CC)}, rotate=0]
\newcommand{\starshape}[2]{
    \begin{scope}[shift={(#1)}, rotate=#2]
        \foreach \angle in {0,60,120,180,240,300} {
            \draw (#1) -- ++(\angle:0.4);
            \draw[fill] ($(#1)+(\angle:0.4)$) circle (1.5pt);
        }
        \draw[fill] (#1) circle (1.5pt);
    \end{scope}
}

% Central point of Y-shape
\coordinate (O) at (0,0);

% Endpoints of Y-shape
\coordinate (A) at (90:0.6);
\coordinate (B) at (210:0.6);
\coordinate (C) at (330:0.6);

% Draw stars at endpoints
\starshape{A}{30}
\starshape{B}{30}
\starshape{C}{30}

\fill[white] (0,0) circle (0.3);

\draw (A)--(0,0);
\draw (B)--(0,0);
\draw (C)--(0,0);

% Mark central intersection point
\draw[fill] (O) circle (1.5pt);
\end{scope}

\begin{scope}[shift={(DD)}, rotate=90]
\newcommand{\starshape}[2]{
    \begin{scope}[shift={(#1)}, rotate=#2]
        \foreach \angle in {0,60,120,180,240,300} {
            \draw (#1) -- ++(\angle:0.4);
            \draw[fill] ($(#1)+(\angle:0.4)$) circle (1.5pt);
        }
        \draw[fill] (#1) circle (1.5pt);
    \end{scope}
}

% Central point of Y-shape
\coordinate (O) at (0,0);

% Endpoints of Y-shape
\coordinate (A) at (90:0.6);
\coordinate (B) at (210:0.6);
\coordinate (C) at (330:0.6);

% Draw stars at endpoints
\starshape{A}{30}
\starshape{B}{30}
\starshape{C}{30}

\fill[white] (0,0) circle (0.3);

\draw (A)--(0,0);
\draw (B)--(0,0);
\draw (C)--(0,0);

% Mark central intersection point
\draw[fill] (O) circle (1.5pt);
\end{scope}

\fill[white] (0,0) circle (0.6);

\draw (AA)--(0,0);
\draw (BB)--(0,0);
\draw (CC)--(0,0);
\draw (DD)--(0,0);

\draw[fill] (0,0) circle (1.5pt);

\end{tikzpicture}
\caption{The figure of $T_{-1,1}$}\label{fig_Tn2}
\end{figure}
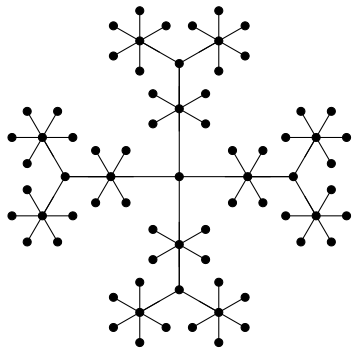

We give a simple example of $\mathbf{b}$ with $\mathbf{b}(-1)=6$, $\mathbf{b}(0)=3$, $\mathbf{b}(1)=4$, to illustrate how the construction works, see Figure \ref{fig_Tn0}, Figure \ref{fig_Tn1}, Figure \ref{fig_Tn2} for the figures of $T_{-1,-1}$, \ldots, $T_{-1,1}$.

\captionsetup[subfigure]{labelformat=empty}

\begin{figure}[htp]
\centering
\begin{subfigure}{0.4\textwidth}
\centering
\begin{tikzpicture}
% Define star command
\newcommand{\starshape}[2]{
    \begin{scope}[shift={(#1)}, rotate=#2]
        \foreach \angle in {0,60,120,180,240,300} {
            \draw (#1) -- ++(\angle:0.4);
            \draw[fill] ($(#1)+(\angle:0.4)$) circle (1.5pt);
        }
        \draw[fill] (#1) circle (1.5pt);
    \end{scope}
}

% Central point of Y-shape
\coordinate (E) at (0,0);

% Endpoints of Y-shape
\coordinate (A) at (90:0.6);
\coordinate (B) at (210:0.6);
\coordinate (C) at (330:0.6);

% Draw stars at endpoints
\starshape{A}{30}
\starshape{B}{30}
\starshape{C}{30}

\fill[white] (0,0) circle (0.3);

\draw (0,0.6)--(0,0);
\draw (B)--(0,0);
\draw (C)--(0,0);

% Mark central intersection point
\draw[fill] (E) circle (1.5pt);

\end{tikzpicture}
\caption{$T_{-1,0}$}
\end{subfigure}
\hfill
\begin{subfigure}{0.4\textwidth}
\centering
\begin{tikzpicture}
\newcommand{\Yshape}[2]{
    \begin{scope}[shift={(#1)}, rotate=#2]
        \foreach \angle in {90,210,330} {
            \draw (#1) -- ++(\angle:0.6);
            \draw[fill] ($(#1)+(\angle:0.6)$) circle (1.5pt);
        }
        \draw[fill] (#1) circle (1.5pt);
    \end{scope}
}

\coordinate (A) at (0:0.6);
\coordinate (B) at (90:0.6);
\coordinate (C) at (180:0.6);
\coordinate (D) at (270:0.6);

\Yshape{A}{90}
\Yshape{B}{60}
\Yshape{C}{-90}
\Yshape{D}{-120}

\end{tikzpicture}
\caption{$T_{0,1}$}
\end{subfigure}
\caption{The figures of $T_{-1,0}$ and $T_{0,1}$}\label{fig_Tn1}
\end{figure}
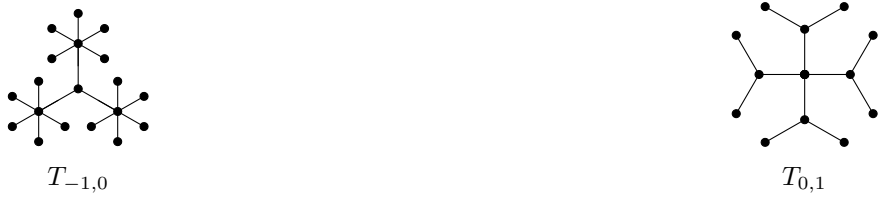

\FloatBarrier

For the metrics $\mathbf{d}_{m,n}$. On $T_{n,n}$, let $\mathbf{d}_{n,n}$ be the standard graph metric multiplied by $2^{n-1}$, under which $T_{n,n}$ has diameter $2^{n}$. On $T_{m,n+1}$, by \ref{item_tree1}, $T_{m,n+1}$ is obtained by gluing $\mathbf{b}(n+1)$ copies of $T_{m,n}$, we give $\mathbf{d}_{m,n+1}$ using $\mathbf{d}_{m,n}$ as follows. Let $\mathbf{r}\in T_{m,n+1}$ be the common vertex that lies in all copies of $T_{m,n}$. For any pair of points $\mathbf{t},\mathbf{s}\in T_{m,n+1}$.
\begin{itemize}
\item If $\mathbf{t}, \mathbf{s}$ lie in the same copy of $T_{m,n}$, then let
$$\mathbf{d}_{m,n+1}(\mathbf{t},\mathbf{s})=\mathbf{d}_{m,n}(\mathbf{t},\mathbf{s}).$$
\item If $\mathbf{t}, \mathbf{s}$ lie in two distinct copies of $T_{m,n}$,  then let
$$\mathbf{d}_{m,n+1}(\mathbf{t},\mathbf{s})=\mathbf{d}_{m,n}(\mathbf{t},\mathbf{r})+\mathbf{d}_{m,n}(\mathbf{r},\mathbf{s}).$$
\end{itemize}
It is obvious that $\mathbf{d}_{m,n+1}$ is a well-defined metric on $T_{m,n+1}$, under which $T_{m,n+1}$ has diameter $2^{n+1}$, and certain compatible conditions hold among $\{(T_{m,n},\mathbf{d}_{m,n}):-\infty<m\le n<+\infty\}$. Letting $m\to-\infty$ or $n\to+\infty$, we obtain $\{(T_{m,n},\mathbf{d}_{m,n}):-\infty\le m\le n\le+\infty\}$. Doing completion to $(T_{-\infty,+\infty},\mathbf{d}_{-\infty,+\infty})$ and $(T_{-\infty,n},\mathbf{d}_{-\infty,n})$ for $n\in\mathbb{Z}$, we obtain an $\mathbb{R}$-tree $(\mathcal{T}(\mathbf{b}),\mathbf{d}_{\mathcal{T}(\mathbf{b})})$ and also sub-$\mathbb{R}$-trees $(\mathcal{T}_n(\mathbf{b}),\mathbf{d}_{\mathcal{T}(\mathbf{b})})$ for $n\in\mathbb{Z}$.

By \ref{item_tree1}, for any $n\in\mathbb{Z}$, $\mathcal{T}_{n}(\mathbf{b})$ consists of $\mathbf{b}(n)$ branches, where each branch is a copy of $\mathcal{T}_{n-1}(\mathbf{b})$, indexed from $0$ to $\mathbf{b}(n)-1$, while each copy of $\mathcal{T}_{n-1}(\mathbf{b})$ consists of $\mathbf{b}(n-1)$ branches, where each branch is a copy of $\mathcal{T}_{n-2}(\mathbf{b})$, indexed from $0$ to $\mathbf{b}(n-1)-1$, hence $\mathcal{T}_{n}(\mathbf{b})$ consists of $\mathbf{b}(n-1)\mathbf{b}(n)$ branches, where each branch is a copy of $\mathcal{T}_{n-2}(\mathbf{b})$, indexed by $\{(\mathbf{s}({n-1}),\mathbf{s}({n})):\mathbf{s}(k)\in\{0,\ldots,\mathbf{b}(k)-1\},k=n-1,n\}$. The index can be chosen such that the $(\mathbf{s}(n-1),\mathbf{s}(n))$-th copy of $\mathcal{T}_{n-2}(\mathbf{b})$ is contained in the $\mathbf{s}(n)$-th copy of $\mathcal{T}_{n-1}(\mathbf{b})$. Similarly, for any $m,n\in\mathbb{Z}$ with $m<n$, $\mathcal{T}_{n}(\mathbf{b})$ consists of $\prod_{k=m+1}^n\mathbf{b}(k)$ branches, where each branch is a copy of $\mathcal{T}_{m}(\mathbf{b})$, which can be index by $\{(\mathbf{s}({m+1}),\ldots,\mathbf{s}({n})):\mathbf{s}(k)\in\{0,\ldots,\mathbf{b}(k)-1\},k=m+1,\ldots,n\}$ accordingly, such that the $(\mathbf{s}({m+1}),\mathbf{s}(m+2),\ldots,\mathbf{s}({n}))$-th copy of $\mathcal{T}_m(\mathbf{b})$ is contained in the $(\mathbf{s}({m+2}),\ldots,\mathbf{s}({n}))$-th copy of $\mathcal{T}_{m+1}(\mathbf{b})$. We say that any such copy of $\mathcal{T}_{m}(\mathbf{b})$ is an $m$-cell. Then any $m$-cell is compact and has diameter $2^m$ in $(\mathcal{T}(\mathbf{b}),\mathbf{d}_{\mathcal{T}(\mathbf{b})})$, it consists of $\mathbf{b}(m)$ $(m-1)$-cells by gluing at one point, which is called the center of the $m$-cell. Roughly speaking, an $m$-cell is comparable to a ball with radius $2^{m-1}$.

A natural projection $\chi:\mathcal{U}(\mathbf{b})\to\mathcal{T}(\mathbf{b})$ can be given as follows. Let $\mathbf{0}\in\mathcal{U}(\mathbf{b})$ be given by $\mathbf{0}(k)=0$ for any $k\in\mathbb{Z}$. Since for any $n\in\mathbb{Z}$, $\mathcal{T}_n(\mathbf{b})$ is compact, has diameter $2^n$ in $(\mathcal{T}(\mathbf{b}),\mathbf{d}_{\mathcal{T}(\mathbf{b})})$, and $\mathcal{T}_n(\mathbf{b})\supseteq\mathcal{T}_{n-1}(\mathbf{b})$, we have $\cap_{n\in\mathbb{Z}}\mathcal{T}_n(\mathbf{b})$ is a one-point set, define $\chi(\mathbf{0})$ as the point in the set. For any $\mathbf{s}\in\mathcal{U}(\mathbf{b})\backslash\{\mathbf{0}\}$, let $n=\max\{k\in\mathbb{Z}:\mathbf{s}(k)\ne0\}$. For any $m<n$, let $K_m$ be the $(\mathbf{s}(m+1),\ldots,\mathbf{s}(n))$-th $m$-cell in $\mathcal{T}_n(\mathbf{b})$, then $K_m$ is compact and has diameter $2^m$ in $(\mathcal{T}(\mathbf{b}),\mathbf{d}_{\mathcal{T}(\mathbf{b})})$. Since $K_m\supseteq K_{m-1}$ for any $m<n$, we have $\cap_{m=n-1}^{-\infty}K_m$ is a one-point set, define $\chi(\mathbf{s})$ as the point in the set. It is obvious that the pre-image of an $n$-cell in $\mathcal{T}(\mathbf{b})$ under $\chi$ is a closed ball with radius $2^{n-1}$ in $\mathcal{U}(\mathbf{b})$, which is also a cylindrical set in $\mathcal{U}(\mathbf{b})$, hence $\chi$ is measurable.

Let $\mathbf{m}_{\mathcal{T}(\mathbf{b})}$ be the pushforward of the measure $\mathbf{m}_{\mathcal{U}(\mathbf{b})}$ on $(\mathcal{U}(\mathbf{b}),\mathbf{d}_{\mathcal{U}(\mathbf{b})})$ under $\chi$, then any $n$-cell has measure $V_{\mathbf{b}}(2^{n-1})$ under $\mathbf{m}_{\mathcal{T}(\mathbf{b})}$, while each $n$-cell has diameter $2^n$ under $\mathbf{d}_{\mathcal{T}(\mathbf{b})}$, hence there exists some positive constant $C$ depending only on $\sup_\mathbb{Z}\mathbf{b}$ such that
$$\frac{1}{C}V_{\mathbf{b}}(r)\le\mathbf{m}_{\mathcal{T}(\mathbf{b})}(B_{\mathcal{T}(\mathbf{b})}(\mathbf{t},r))\le CV_{\mathbf{b}}(r)\text{ for any }\mathbf{t}\in\mathcal{T}(\mathbf{b}),r>0.$$

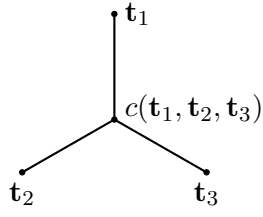
\begin{figure}[ht]
\centering
\begin{tikzpicture}[scale=0.7]

\draw[thick] (0,0)--(0,2);
\draw[thick] (0,0)--(-1.732,-1);
\draw[thick] (0,0)--(1.732,-1);

\filldraw[fill=black] (0,2) circle (0.05);
\filldraw[fill=black] (0,0) circle (0.05);
\filldraw[fill=black] (-1.732,-1) circle (0.05);
\filldraw[fill=black] (1.732,-1) circle (0.05);

\node[right] at (0,0.2) {$c(\mathbf{t}_1,\mathbf{t}_2,\mathbf{t}_3)$};
\node[right] at (0,2) {$\mathbf{t}_1$};
\node[below] at (-1.732,-1) {$\mathbf{t}_2$};
\node[below] at (1.732,-1) {$\mathbf{t}_3$};

\end{tikzpicture}
\caption{$c(\mathbf{t}_1,\mathbf{t}_2,\mathbf{t}_3)$}\label{fig_ct123}
\end{figure}

Now we obtain an $\mathbb{R}$-tree $(\mathcal{T}(\mathbf{b}),\mathbf{d}_{\mathcal{T}(\mathbf{b})})$ endowed with a measure $\mathbf{m}_{\mathcal{T}(\mathbf{b})}$. For any $\mathbf{t}_1,\mathbf{t}_2\in\mathcal{T}(\mathbf{b})$, there exists a unique geodesic $[\mathbf{t}_1,\mathbf{t}_2]$ connecting $\mathbf{t}_1,\mathbf{t}_2$ with length $\mathbf{d}_{\mathcal{T}(\mathbf{d})}(\mathbf{t}_1,\mathbf{t}_2)$. For any $\mathbf{t}_1,\mathbf{t}_2,\mathbf{t}_3\in\mathcal{T}(\mathbf{b})$, we have $[\mathbf{t}_2,\mathbf{t}_3]\subseteq[\mathbf{t}_1,\mathbf{t}_2]\cup[\mathbf{t}_1,\mathbf{t}_3]$, and there exists a unique $c(\mathbf{t}_1,\mathbf{t}_2,\mathbf{t}_3)\in\mathcal{T}(\mathbf{b})$ such that $[\mathbf{t}_1,\mathbf{t}_2]\cap[\mathbf{t}_1,\mathbf{t}_3]=[\mathbf{t}_1,c(\mathbf{t}_1,\mathbf{t}_2,\mathbf{t}_3)]$, see \cite[Lemma 3.20]{Eva08} and Figure \ref{fig_ct123}.

There exists a unique $\sigma$-finite Borel measure $\lambda_{\mathcal{T}(\mathbf{b})}$ on $(\mathcal{T}(\mathbf{b}),\mathbf{d}_{\mathcal{T}(\mathbf{b})})$, called length measure, satisfying that
$$\lambda_{\mathcal{T}(\mathbf{b})}(]\mathbf{t}_1,\mathbf{t}_2[)=\mathbf{d}_{\mathcal{T}(\mathbf{b})}(\mathbf{t}_1,\mathbf{t}_2)\text{ for any distinct }\mathbf{t}_1,\mathbf{t}_2\in\mathcal{T}(\mathbf{b}),$$
where $]\mathbf{t}_1,\mathbf{t}_2[=[\mathbf{t}_1,\mathbf{t}_2]\backslash\{\mathbf{t}_1,\mathbf{t}_2\}$. Indeed, let $\mathcal{T}^o(\mathbf{b})=\cup_{\mathbf{t}_1,\mathbf{t}_2\in\mathcal{T}(\mathbf{b})}]\mathbf{t}_1,\mathbf{t}_2[$ be the skeleton of $\mathcal{T}(\mathbf{b})$, then $\lambda_{\mathcal{T}(\mathbf{b})}(\mathcal{T}(\mathbf{b})\backslash\mathcal{T}^o(\mathbf{b}))=0$, in particular, $\lambda_{\mathcal{T}(\mathbf{b})}$ is the trace onto $\mathcal{T}^o(\mathbf{b})$ of $1$-dimensional Hausdorff measure on $\mathcal{T}(\mathbf{b})$, see \cite[Subsection 2.4]{EPW06}.

We say that $f\in C(\mathcal{T}(\mathbf{b}))$ is locally absolutely continuous if for any $\varepsilon>0$, for any Borel set $\mathcal{S}\subseteq\mathcal{T}(\mathbf{b})$ with $\lambda_{\mathcal{T}(\mathbf{b})}(\mathcal{S})\in(0,+\infty)$, there exists $\delta=\delta(\varepsilon,\mathcal{S})>0$ such that for any integer $n\ge1$, for any disjoint geodesics $[\mathbf{t}_1,\mathbf{s}_1]$, \ldots, $[\mathbf{t}_n,\mathbf{s}_n]$ contained in $\mathcal{S}$ with $\sum_{k=1}^n\mathbf{d}_{\mathcal{T}(\mathbf{b})}(\mathbf{t}_k,\mathbf{s}_k)<\delta$, we have $\sum_{k=1}^n|f(\mathbf{t}_k)-f(\mathbf{s}_k)|<\varepsilon$. Let $\mathcal{A}(\mathcal{T}(\mathbf{b}))$ be the family of all locally absolutely continuous functions.

Fix an arbitrary point $\mathbf{\rho}\in\mathcal{T}(\mathbf{b})$, we can define
\begin{itemize}
\item A wedge product $\wedge_\rho$ by setting $x\wedge_\rho y=c(\rho,x,y)$.
\item An orientation-sensitive integration by
$$\int_x^y g(z)\lambda_{\mathcal{T}(\mathbf{b})}(\mathrm{d} z)=\int_{[x\wedge_\rho y,y]}g(z)\lambda_{\mathcal{T}(\mathbf{b})}(\mathrm{d} z)-\int_{[x\wedge_\rho y,x]}g(z)\lambda_{\mathcal{T}(\mathbf{b})}(\mathrm{d} z).$$
\end{itemize}

Then for any $f\in\mathcal{A}(\mathcal{T}(\mathbf{b}))$, there exists a unique $g\in L^1_{loc}(\mathcal{T}(\mathbf{b});\lambda_{\mathcal{T}(\mathbf{b})})$, up to $\lambda_{\mathcal{T}(\mathbf{b})}$-measure zero sets, such that
$$f(y)-f(x)=\int_x^yg(z)\lambda_{\mathcal{T}(\mathbf{b})}(\mathrm{d} z)$$
for any $x,y\in\mathcal{T}(\mathbf{b})$, define the gradient $\nabla_\mathcal{T} f=\nabla_\mathcal{T}^\rho f$ as the function $g$. Note that the choice of the point $\rho\in\mathcal{T}(\mathbf{b})$ may affect the sign of $\nabla_\mathcal{T} f$, but it does not affect the value of $|\nabla_\mathcal{T} f|$, see \cite[Proposition 1.1, Definition 1.2, Remark 1.3]{AEW13} for further details.

Let
\begin{align*}
&\mathcal{E}^{\mathcal{T}}(f)=\int_{\mathcal{T}(\mathbf{b})}|\nabla_\mathcal{T} f|^p\mathrm{d}\lambda_{\mathcal{T}(\mathbf{b})},\\
&\mathcal{F}^{\mathcal{T}}=\left\{f\in L^p(\mathcal{T}(\mathbf{b});\mathbf{m}_{\mathcal{T}(\mathbf{b})})\cap \mathcal{A}(\mathcal{T}(\mathbf{b})):\int_{\mathcal{T}(\mathbf{b})}|\nabla_\mathcal{T} f|^p\mathrm{d}\lambda_{\mathcal{T}(\mathbf{b})}<+\infty\right\}.
\end{align*}

By virtue of the tree property, $\mathcal{F}^\mathcal{T}$ automatically contains a large family of Lipshcitz functions arising from distance functions, as demonstrated by the following result.

\begin{lemma}\label{lem_tree_cutoff}
For any $n\in\mathbb{Z}$ and any $n$-cell ${K}$, let $\mathcal{N}\left({{K}}\right)=\bigcup_{\widetilde{K}:n\text{-cell,}\widetilde{K}\cap {K}\ne\emptyset}\widetilde{K}$ be the $n$-cell neighborhood of $K$ in $\mathcal{T}(\mathbf{b})$, then there exists $\phi_{{K}}\in\mathcal{F}^{\mathcal{T}}\cap C_c(\mathcal{T}(\mathbf{b}))$ with $0\le\phi_{{K}}\le1$ in $\mathcal{T}(\mathbf{b})$, $\phi_{{K}}=1$ on ${{K}}$, $\phi_{K}=0$ on $\mathcal{T}(\mathbf{b})\backslash\mathcal{N}({{K}})$ such that
$$\mathcal{E}^\mathcal{T}(\phi_{{{K}}})\le \frac{2\left(\sup_{\mathbb{Z}}\mathbf{b}-1\right)}{2^{(p-1)n}}.$$
Moreover, for any $f\in\mathcal{F}^\mathcal{T}$, we have
\begin{align*}
&\int_{\mathcal{N}(K)}|f|^p|\nabla_\mathcal{T}\phi_K|^p\mathrm{d}\lambda_{\mathcal{T}(\mathbf{b})}\\
&\le2^p3^{p-1}(\sup_{\mathbb{Z}}\mathbf{b}-1)\int_{\mathcal{N}(K)}|\nabla_\mathcal{T} f|^p\mathrm{d}\lambda_{\mathcal{T}(\mathbf{b})}+\frac{2^p(\sup_{\mathbb{Z}}\mathbf{b}-1)}{2^{(p-1)n}}\dashint_{\mathcal{N}(K)}|f|^p\mathrm{d}\mathbf{m}_{\mathcal{T}(\mathbf{b})}.
\end{align*}
\end{lemma}

\begin{proof}
By \ref{item_tree2}, any two distinct $n$-cells intersect at most one point and any $n$-cell contains at most two points shared with other $n$-cells, hence there exist non-negative integers $M,N\le\sup_{\mathbb{Z}}\mathbf{b}-1$, $\mathbf{t},\mathbf{s}\in K$, and distinct $n$-cells $\widetilde{K}_1$, \ldots $\widetilde{K}_M$ and $\widehat{K}_1$, \ldots $\widehat{K}_N$ such that $\widetilde{K}_k\cap K=\{\mathbf{t}\}$ for any $k=1,\ldots,M$, and $\widehat{K}_l\cap K=\{\mathbf{s}\}$ for any $l=1,\ldots,N$, here it is possible that $M$ or $N$ is 0, which means that $K$ has only one point shared with other $n$-cells. Let $\mathbf{t}^{(k)}\in\widetilde{K}_k$ (resp. $\mathbf{s}^{(l)}\in\widehat{K}_l$) denote the other point in $\widetilde{K}_k$ (resp. $\widehat{K}_l$) which is shared with other $n$-cells, if such a point exists; see Figure \ref{fig_KNK}.

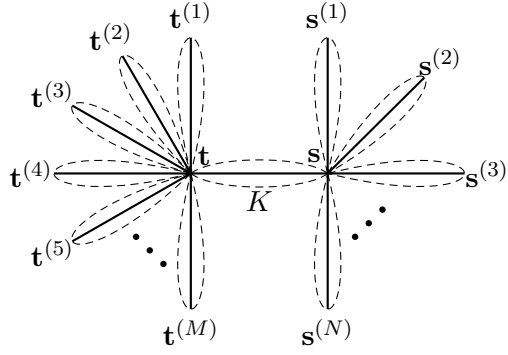
\begin{figure}[ht]
\centering
\begin{tikzpicture}[scale=1.2]

\draw[densely dashed] (-0.75,0) .. controls (-0.5,0.2) and (0.5,0.2) .. (0.75,0);
\draw[densely dashed] (-0.75,0) .. controls (-0.5,-0.2) and (0.5,-0.2) .. (0.75,0);
\draw[thick] (-0.75,0)--(0.75,0);

\node at (-0.6,0.2) {$\mathbf{t}$};
\node at (0.6,0.2) {$\mathbf{s}$};

% Draw outward petals on the right
\begin{scope}[shift={(0.75,0)}]
\foreach \angle [count=\i from 1] in {90,45,0,-90}
{
\begin{scope}[rotate=\angle]
\draw[densely dashed] (0,0) .. controls (2,0.5) and (2,-0.5) .. (0,0);
\draw[thick] (0,0)--(1.5,0);
\ifnum\i<4
\node at (1.75,0) {$\mathbf{s}^{(\i)}$};
\else
\node at (1.75,0) {$\mathbf{s}^{(N)}$};
\fi
\end{scope}
}
\begin{scope}[shift={(0.3,-0.7)}]
\fill (0,0) circle (1pt);
\fill (0.15,0.15) circle (1pt);
\fill (0.3,0.3) circle (1pt);
\end{scope}
\end{scope}

% Draw outward petals on the left
\begin{scope}[shift={(-0.75,0)}]
\foreach \angle [count=\i from 1] in {90, 120, 150, 180, 210, 270}
{
\begin{scope}[rotate=\angle]
\draw[densely dashed] (0,0) .. controls (2,0.5) and (2,-0.5) .. (0,0);
\draw[thick] (0,0)--(1.5,0);
\ifnum\i<6
\node at (1.75,0) {$\mathbf{t}^{(\i)}$};
\else
\node at (1.75,0) {$\mathbf{t}^{(M)}$};
\fi
\end{scope}
}
\begin{scope}[shift={(-0.3,-1)},rotate=90]
\fill (0,0) circle (1pt);
\fill (0.15,0.15) circle (1pt);
\fill (0.3,0.3) circle (1pt);
\end{scope}
\end{scope}

\node at (0,-0.3) {$K$};
\end{tikzpicture}
\caption{$\mathcal{N}(K)$}\label{fig_KNK}
\end{figure}

It is obvious that $\mathcal{N}(K)=K\cup\cup_{k=1}^M\widetilde{K}_k\cup\cup_{l=1}^N\widehat{K}_l$. Let $\phi_K$ be given by $\phi_K=1$ on $K$, $\phi_K=0$ on $\mathcal{T}(\mathbf{b})\backslash\mathcal{N}(K)$, on each $\widetilde{K}_k$,
$$\phi_K=\frac{\mathbf{d}_{\mathcal{T}(\mathbf{b})}(c(\cdot,\mathbf{t},\mathbf{t}^{(k)}),\mathbf{t}^{(k)})}{\mathbf{d}_{\mathcal{T}(\mathbf{b})}(\mathbf{t},\mathbf{t}^{(k)})}=\frac{1}{2^n}{\mathbf{d}_{\mathcal{T}(\mathbf{b})}(c(\cdot,\mathbf{t},\mathbf{t}^{(k)}),\mathbf{t}^{(k)})},$$
and on each $\widehat{K}_l$,
$$\phi_K=\frac{\mathbf{d}_{\mathcal{T}(\mathbf{b})}(c(\cdot,\mathbf{s},\mathbf{s}^{(l)}),\mathbf{s}^{(l)})}{\mathbf{d}_{\mathcal{T}(\mathbf{b})}(\mathbf{s},\mathbf{s}^{(l)})}=\frac{1}{2^n}{\mathbf{d}_{\mathcal{T}(\mathbf{b})}(c(\cdot,\mathbf{s},\mathbf{s}^{(l)}),\mathbf{s}^{(l)})}.$$
It is obvious that $\phi_K\in C_c(\mathcal{T}(\mathbf{b}))$ is well-defined. For $\widetilde{K}_k$, on $[\mathbf{t},\mathbf{t}^{(k)}]$, we have $\phi_K=\frac{\mathbf{d}_{\mathcal{T}(\mathbf{b})}(\cdot,\mathbf{t}^{(k)})}{\mathbf{d}_{\mathcal{T}(\mathbf{b})}(\mathbf{t},\mathbf{t}^{(k)})}=\frac{1}{2^n}{\mathbf{d}_{\mathcal{T}(\mathbf{b})}(\cdot,\mathbf{t}^{(k)})}$, and on each component of $\widetilde{K}_k\backslash[\mathbf{t},\mathbf{t}^{(k)}]$, since $c(\cdot,\mathbf{t},\mathbf{t}^{(k)})$ is constant, we have $\phi_K$ is constant, which gives
$$\int_{\widetilde{K}_k}|\nabla_\mathcal{T}\phi_K|^p\mathrm{d}\lambda_{\mathcal{T}(\mathbf{b})}=\int_{[\mathbf{t},\mathbf{t}^{(k)}]}\frac{1}{2^{pn}}\mathrm{d}\lambda_{\mathcal{T}(\mathbf{b})}=\frac{1}{2^{(p-1)n}}.$$
Similarly, we have
$$\int_{\widehat{K}_l}|\nabla_\mathcal{T}\phi_K|^p\mathrm{d}\lambda_{\mathcal{T}(\mathbf{b})}=\frac{1}{2^{(p-1)n}}.$$
Hence
\begin{align*}
&\int_{\mathcal{T}(\mathbf{b})}|\nabla_\mathcal{T}\phi_K|^p\mathrm{d}\lambda_{\mathcal{T}(\mathbf{b})}\\
&=\sum_{k=1}^M\int_{\widetilde{K}_k}|\nabla_\mathcal{T}\phi_K|^p\mathrm{d}\lambda_{\mathcal{T}(\mathbf{b})}+\sum_{l=1}^N\int_{\widehat{K}_l}|\nabla_\mathcal{T}\phi_K|^p\mathrm{d}\lambda_{\mathcal{T}(\mathbf{b})}\\
&=(M+N)\frac{1}{2^{(p-1)n}}\le \frac{2\left(\sup_{\mathbb{Z}}\mathbf{b}-1\right)}{2^{(p-1)n}},
\end{align*}
and $\phi_K\in\mathcal{F}^\mathcal{T}$. Moreover, for any $f\in\mathcal{F}^\mathcal{T}\subseteq C(\mathcal{T}(\mathbf{b}))$, there exists $x_0\in\mathcal{N}(K)$ such that $|f(x_0)|^p=\dashint_{\mathcal{N}(K)}|f|^p\mathrm{d} \mathbf{m}_{\mathcal{T}(\mathbf{b})}$, and
$$\int_{\mathcal{N}(K)}|f|^p|\nabla_\mathcal{T}\phi_K|^p\mathrm{d}\lambda_{\mathcal{T}(\mathbf{b})}=\frac{1}{2^{pn}}\left(\sum_{k=1}^M\int_{[\mathbf{t},\mathbf{t}^{(k)}]}|f|^p\mathrm{d}\lambda_{\mathcal{T}(\mathbf{b})}+\sum_{l=1}^N\int_{[\mathbf{s},\mathbf{s}^{(l)}]}|f|^p\mathrm{d}\lambda_{\mathcal{T}(\mathbf{b})}\right).$$
For any $x\in\mathcal{N}(K)$, we have $[x,x_0]\subseteq\mathcal{N}(K)$ and $\mathbf{d}_{\mathcal{T}(\mathbf{b})}(x,x_0)\le3\cdot2^n$, hence by H\"older's inequality, we have
\begin{align*}
&|f(x)-f(x_0)|\le\int_{[x,x_0]}|\nabla_\mathcal{T} f|\mathrm{d}\lambda_{\mathcal{T}(\mathbf{b})}\\
&\le \mathbf{d}_{\mathcal{T}(\mathbf{b})}(x,x_0)^{1-\frac{1}{p}}\left(\int_{[x,x_0]}|\nabla_\mathcal{T} f|^p\mathrm{d}\lambda_{\mathcal{T}(\mathbf{b})}\right)^{\frac{1}{p}}\le\left(3\cdot2^n\right)^{1-\frac{1}{p}}\left(\int_{\mathcal{N}(K)}|\nabla_\mathcal{T} f|^p\mathrm{d}\lambda_{\mathcal{T}(\mathbf{b})}\right)^{\frac{1}{p}},
\end{align*}
which gives
\begin{align*}
&|f(x)|^p\le2^{p-1}\left(|f(x)-f(x_0)|^p+|f(x_0)|^p\right)\\
&\le2^{p-1}\left(\left(3\cdot2^n\right)^{p-1}\int_{\mathcal{N}(K)}|\nabla_\mathcal{T} f|^p\mathrm{d}\lambda_{\mathcal{T}(\mathbf{b})}+\dashint_{\mathcal{N}(K)}|f|^p\mathrm{d} \mathbf{m}_{\mathcal{T}(\mathbf{b})}\right).
\end{align*}
Hence
\begin{align*}
&\int_{\mathcal{N}(K)}|f|^p|\nabla_\mathcal{T}\phi_K|^p\mathrm{d}\lambda_{\mathcal{T}(\mathbf{b})}\\
&\le \frac{1}{2^{pn}}\cdot2^{p-1}\left(\left(3\cdot2^n\right)^{p-1}\int_{\mathcal{N}(K)}|\nabla_\mathcal{T} f|^p\mathrm{d}\lambda_{\mathcal{T}(\mathbf{b})}+\dashint_{\mathcal{N}(K)}|f|^p\mathrm{d} \mathbf{m}_{\mathcal{T}(\mathbf{b})}\right)\cdot(M+N)\cdot2^n\\
&\le2^p3^{p-1}(\sup_{\mathbb{Z}}\mathbf{b}-1)\int_{\mathcal{N}(K)}|\nabla_\mathcal{T} f|^p\mathrm{d}\lambda_{\mathcal{T}(\mathbf{b})}+\frac{2^p(\sup_{\mathbb{Z}}\mathbf{b}-1)}{2^{(p-1)n}}\dashint_{\mathcal{N}(K)}|f|^p\mathrm{d} \mathbf{m}_{\mathcal{T}(\mathbf{b})}.
\end{align*}
\end{proof}

It is easy to have the following result.

\begin{lemma}\label{lem_tree_energy}
$(\mathcal{E}^{\mathcal{T}},\mathcal{F}^{\mathcal{T}})$ is a $p$-energy on $(\mathcal{T}(\mathbf{b}),\mathbf{d}_{\mathcal{T}(\mathbf{b})},\mathbf{m}_{\mathcal{T}(\mathbf{b})})$ with a $p$-energy measure $\Gamma^\mathcal{T}$ given by
$$\Gamma^\mathcal{T}(f)(A)=\int_A|\nabla_\mathcal{T} f|^p\mathrm{d}\lambda_{\mathcal{T}(\mathbf{b})}$$
for any $f\in\mathcal{F}^\mathcal{T}$, $A\in\mathcal{B}(\mathcal{T}(\mathbf{b}))$.
\end{lemma}

\begin{proof}
We only give the proof of the uniform denseness of $\mathcal{F}^\mathcal{T}\cap C_c(\mathcal{T}(\mathbf{b}))$ in $C_c(\mathcal{T}(\mathbf{b}))$. The other parts are easy and similar to the proof of $(\int_{\mathbb{R}}|f'(x)|^p\mathrm{d} x,W^{1,p}(\mathbb{R}))$ is a $p$-energy on $\mathbb{R}$ with a $p$-energy measure given by $(f,A)\mapsto\int_A|f'(x)|^p\mathrm{d} x$.

Indeed, since $\mathcal{F}^\mathcal{T}\cap C_c(\mathcal{T}(\mathbf{b}))$ is a sub-algebra of $C_c(\mathcal{T}(\mathbf{b}))$, by the Stone-Weierstra{\ss} theorem, we only need to show that $\mathcal{F}^\mathcal{T}\cap C_c(\mathcal{T}(\mathbf{b}))$ separates points and vanishes nowhere. For any distinct $\mathbf{t},\mathbf{s}\in\mathcal{T}(\mathbf{b})$, there exist $n\in\mathbb{Z}$ sufficiently small and an $n$-cell $K$ such that $\mathbf{t}\in K$ and $\mathbf{s}\not\in\mathcal{N}(K)$. Let $\phi_K\in\mathcal{F}^\mathcal{T}\cap C_c(\mathcal{T}(\mathbf{b}))$ be given by Lemma \ref{lem_tree_cutoff}, then $\phi_K(\mathbf{t})=1\ne0=\phi_K(\mathbf{s})$.
\end{proof}

We have the Poincar\'e inequality on $\mathcal{T}(\mathbf{b})$ as follows.

\begin{proposition}\label{prop_tree_PI}
There exists $C>0$ such that for any ball ${B_{\mathcal{T}(\mathbf{b})}}$ with radius $r$, for any $f\in\mathcal{F}^{\mathcal{T}}$, we have
$$\int_{B_{\mathcal{T}(\mathbf{b})}}\left| f-\dashint_{B_{\mathcal{T}(\mathbf{b})}}f\mathrm{d}\mathbf{m}_{\mathcal{T}(\mathbf{b})}\right|^p\mathrm{d}\mathbf{m}_{\mathcal{T}(\mathbf{b})}\le Cr^{p-1}V_{\mathbf{b}}(r)\int_{{B_{\mathcal{T}(\mathbf{b})}}}|\nabla_\mathcal{T} f|^p\mathrm{d}\lambda_{\mathcal{T}(\mathbf{b})}.$$
\end{proposition}

\begin{proof}
For any $x,y\in {B_{\mathcal{T}(\mathbf{b})}}$, by the tree property, we have $[x,y]\subseteq{B_{\mathcal{T}(\mathbf{b})}}$, by H\"older's inequality, we have
\begin{align*}
&|f(x)-f(y)|\le\int_{[x,y]}|\nabla_\mathcal{T} f|\mathrm{d}\lambda_{\mathcal{T}(\mathbf{b})}\\
&\le \mathbf{d}_{\mathcal{T}(\mathbf{b})}(x,y)^{1-\frac{1}{p}}\left(\int_{[x,y]}|\nabla_\mathcal{T} f|^p\mathrm{d}\lambda_{\mathcal{T}(\mathbf{b})}\right)^{\frac{1}{p}}\le(2r)^{1-\frac{1}{p}}\left(\int_{{B_{\mathcal{T}(\mathbf{b})}}}|\nabla_\mathcal{T} f|^p\mathrm{d}\lambda_{\mathcal{T}(\mathbf{b})}\right)^{\frac{1}{p}}.
\end{align*}
Hence
\begin{align*}
&\int_{B_{\mathcal{T}(\mathbf{b})}}\left| f-\dashint_{B_{\mathcal{T}(\mathbf{b})}}f\mathrm{d}\mathbf{m}_{\mathcal{T}(\mathbf{b})}\right|^p\mathrm{d}\mathbf{m}_{\mathcal{T}(\mathbf{b})}\\
&\le \frac{1}{\mathbf{m}_{\mathcal{T}(\mathbf{b})}({B_{\mathcal{T}(\mathbf{b})}})}\int_{B_{\mathcal{T}(\mathbf{b})}}\int_{B_{\mathcal{T}(\mathbf{b})}}|f(x)-f(y)|^p\mathbf{m}_{\mathcal{T}(\mathbf{b})}(\mathrm{d} x)\mathbf{m}_{\mathcal{T}(\mathbf{b})}(\mathrm{d} y)\\
&\le {\mathbf{m}_{\mathcal{T}(\mathbf{b})}({B_{\mathcal{T}(\mathbf{b})}})}(2r)^{p-1}\int_{{B_{\mathcal{T}(\mathbf{b})}}}|\nabla_\mathcal{T} f|^p\mathrm{d}\lambda_{\mathcal{T}(\mathbf{b})}\asymp r^{p-1}V_{\mathbf{b}}(r)\int_{{B_{\mathcal{T}(\mathbf{b})}}}|\nabla_\mathcal{T} f|^p\mathrm{d}\lambda_{\mathcal{T}(\mathbf{b})}.
\end{align*}
\end{proof}

We have the capacity upper bound and the cutoff Sobolev inequality on $\mathcal{T}(\mathbf{b})$ as follows.

\begin{proposition}\label{prop_tree_CS}
There exist $C_1, C_2, C_3, C_4>0$ such that for any ball $B_{\mathcal{T}(\mathbf{b})}$ with radius $r$, there exists $\phi_\mathcal{T}\in\mathcal{F}^\mathcal{T}\cap C_c(\mathcal{T}(\mathbf{b}))$ with $0\le \phi_\mathcal{T}\le 1$ in $\mathcal{T}(\mathbf{b})$, $\phi_\mathcal{T}=1$ in $B_{\mathcal{T}(\mathbf{b})}$, $\phi_\mathcal{T}=0$ on $\mathcal{T}(\mathbf{b})\backslash(8B_{\mathcal{T}(\mathbf{b})})$ such that
$$\mathcal{E}^\mathcal{T}(\phi_\mathcal{T})\le \frac{C_1}{r^{p-1}}.$$
Hence
$$\mathrm{cap}_{\mathcal{T}(\mathbf{b})}(B_{\mathcal{T}(\mathbf{b})},\mathcal{T}(\mathbf{b})\backslash(8B_{\mathcal{T}(\mathbf{b})}))\le \frac{C_1}{r^{p-1}}\le C_2\frac{\mathbf{m}_{\mathcal{T}(\mathbf{b})}(B_{\mathcal{T}(\mathbf{b})})}{r^{p-1}V_\mathbf{b}(r)}.$$
Moreover, for any $f\in\mathcal{F}^\mathcal{T}$, we have
$$\int_{8B_{\mathcal{T}(\mathbf{b})}}|f|^p|\nabla_\mathcal{T}\phi_\mathcal{T}|^p\mathrm{d}\lambda_{\mathcal{T}(\mathbf{b})}\le C_3\int_{8B_{\mathcal{T}(\mathbf{b})}}|\nabla_\mathcal{T} f|^p\mathrm{d}\lambda_{\mathcal{T}(\mathbf{b})}+\frac{C_4}{r^{p-1}V_\mathbf{b}(r)}\int_{8B_{\mathcal{T}(\mathbf{b})}}|f|^p\mathrm{d}\mathbf{m}_{\mathcal{T}(\mathbf{b})}.$$
\end{proposition}

\begin{proof}
Write $B_{\mathcal{T}(\mathbf{b})}=B_{\mathcal{T}(\mathbf{b})}(x,r)$. Let $n$ be the integer satisfying $2^{n-1}\le r<2^{n}$, there exists an $n$-cell $K$ such that $x\in K$. Let $\mathcal{N}\left({{K}}\right)=\bigcup_{\widetilde{K}:n\text{-cell,}\widetilde{K}\cap {K}\ne\emptyset}\widetilde{K}$ be the $n$-cell neighborhood of $K$ in $\mathcal{T}(\mathbf{b})$, then $B_{\mathcal{T}(\mathbf{b})}(x,r)\subseteq\mathcal{N}(K)$, see Figure \ref{fig_KNK}. For any $n$-cell $\widetilde{K}$ with $\widetilde{K}\cap K\ne\emptyset$, by Lemma \ref{lem_tree_cutoff}, there exists $\phi_{\widetilde{K}}\in\mathcal{F}^\mathcal{T}\cap C_c(\mathcal{T}(\mathbf{b}))$ with $0\le\phi_{\widetilde{K}}\le1$ in $\mathcal{T}(\mathbf{b})$, $\phi_{{\widetilde{K}}}=1$ on ${\widetilde{K}}$, $\phi_{{\widetilde{K}}}=0$ on $\mathcal{T}(\mathbf{b})\backslash\mathcal{N}({\widetilde{K}})$ such that $\mathcal{E}^\mathcal{T}(\phi_{{\widetilde{K}}})\le \frac{2(\sup_\mathbb{Z}\mathbf{b}-1)}{2^{(p-1)n}}$.

Let
\begin{equation}\label{eq_tree_cutoff}
\phi_\mathcal{T}=\max_{\widetilde{K}:n\text{-cell, }\widetilde{K}\cap {K}\ne\emptyset}\phi_{\widetilde{K}}.
\end{equation}
Since
\begin{equation}\label{eq_tree_num}
\#\{\widetilde{K}:n\text{-cell, }\widetilde{K}\cap {K}\ne\emptyset\}\le 2\sup_{\mathbb{Z}}\mathbf{b}-1,
\end{equation}
we have $\phi_\mathcal{T}\in\mathcal{F}^\mathcal{T}\cap C_c(\mathcal{T}(\mathbf{b}))$ is well-defined, $0\le\phi_\mathcal{T}\le1$ in $\mathcal{T}(\mathbf{b})$, $\phi_\mathcal{T}=1$ on $\mathcal{N}(K)\supseteq B_{\mathcal{T}(\mathbf{b})}(x,r)$, $\phi_\mathcal{T}=0$ on $\mathcal{T}(\mathbf{b})\backslash\cup_{\widetilde{K}:n\text{-cell, }\widetilde{K}\cap {K}\ne\emptyset}\mathcal{N}(\widetilde{K})$. Since $\mathcal{N}(\widetilde{K})\subseteq B_{\mathcal{T}(\mathbf{b})}(x,2^{n+2})$ for any $n$-cell $\widetilde{K}$ with $\widetilde{K}\cap K\ne\emptyset$, we have $\phi_\mathcal{T}=0$ on $\mathcal{T}(\mathbf{b})\backslash B_{\mathcal{T}(\mathbf{b})}(x,2^{n+2})\supseteq \mathcal{T}(\mathbf{b})\backslash B_{\mathcal{T}(\mathbf{b})}(x,8r)$. Moreover, there exists $C\ge1$ depending only on $p, \sup_{\mathbb{Z}}\mathbf{b}$ such that
\begin{align*}
&\mathcal{E}^\mathcal{T}(\phi_\mathcal{T})\le C\sum_{\widetilde{K}:n\text{-cell, }\widetilde{K}\cap {K}\ne\emptyset}\mathcal{E}^\mathcal{T}(\phi_{\widetilde{K}})\\
&\le C\left(2\sup_{\mathbb{Z}}\mathbf{b}-1\right)\frac{2(\sup_\mathbb{Z}\mathbf{b}-1)}{2^{(p-1)n}}\le \frac{4C\left(\sup_{\mathbb{Z}}\mathbf{b}\right)^2}{r^{p-1}}.
\end{align*}
Hence
\begin{align*}
&\mathrm{cap}_{\mathcal{T}(\mathbf{b})}(B_{\mathcal{T}(\mathbf{b})}(x,r),\mathcal{T}(\mathbf{b})\backslash(B_{\mathcal{T}(\mathbf{b})}(x,8r)))\le\frac{4C\left(\sup_{\mathbb{Z}}\mathbf{b}\right)^2}{r^{p-1}}\\
&={4C\left(\sup_{\mathbb{Z}}\mathbf{b}\right)^2}\frac{V_\mathbf{b}(r)}{r^{p-1}V_\mathbf{b}(r)}\asymp\frac{\mathbf{m}_{\mathcal{T}(\mathbf{b})}(B_{\mathcal{T}(\mathbf{b})}(x,r))}{r^{p-1}V_\mathbf{b}(r)}.
\end{align*}
Moreover, for any $f\in\mathcal{F}^\mathcal{T}$, by Lemma \ref{lem_tree_cutoff}, we have
$$\int_{\mathcal{N}(\widetilde{K})}|f|^p|\nabla_\mathcal{T}\phi_{\widetilde{K}}|^p\mathrm{d}\lambda_{\mathcal{T}(\mathbf{b})}\lesssim\int_{\mathcal{N}(\widetilde{K})}|\nabla_\mathcal{T} f|^p\mathrm{d}\lambda_{\mathcal{T}(\mathbf{b})}+\frac{1}{2^{(p-1)n}}\dashint_{\mathcal{N}(\widetilde{K})}|f|^p\mathrm{d}\mathbf{m}_{\mathcal{T}(\mathbf{b})}.$$
Therefore, we have
\begin{align*}
&\int_{B_{\mathcal{T}(\mathbf{b})}(x,8r)}|f|^p|\nabla_\mathcal{T}\phi_\mathcal{T}|^p\mathrm{d}\lambda_{\mathcal{T}(\mathbf{b})}\\
&\overset{(\star)}{\scalebox{2}[1]{$\lesssim$}}\sum_{\widetilde{K}:n\text{-cell, }\widetilde{K}\cap {K}\ne\emptyset}\int_{B_{\mathcal{T}(\mathbf{b})}(x,8r)}|f|^p|\nabla_\mathcal{T}\phi_{\widetilde{K}}|^p\mathrm{d}\lambda_{\mathcal{T}(\mathbf{b})}\\
&\overset{(\dagger)}{\scalebox{2}[1]{$=$}}\sum_{\widetilde{K}:n\text{-cell, }\widetilde{K}\cap {K}\ne\emptyset}\int_{\mathcal{N}(\widetilde{K})}|f|^p|\nabla_\mathcal{T}\phi_{\widetilde{K}}|^p\mathrm{d}\lambda_{\mathcal{T}(\mathbf{b})}\\
&\lesssim\sum_{\widetilde{K}:n\text{-cell, }\widetilde{K}\cap {K}\ne\emptyset}\left(\int_{\mathcal{N}(\widetilde{K})}|\nabla_\mathcal{T} f|^p\mathrm{d}\lambda_{\mathcal{T}(\mathbf{b})}+\frac{1}{2^{(p-1)n}}\dashint_{\mathcal{N}(\widetilde{K})}|f|^p\mathrm{d}\mathbf{m}_{\mathcal{T}(\mathbf{b})}\right)\\
&\overset{(\diamond)}{\scalebox{2}[1]{$\lesssim$}}\int_{B_{\mathcal{T}(\mathbf{b})}(x,8r)}|\nabla_\mathcal{T} f|^p\mathrm{d}\lambda_{\mathcal{T}(\mathbf{b})}+\frac{1}{r^{p-1}V_\mathbf{b}(r)}\int_{B_{\mathcal{T}(\mathbf{b})}(x,8r)}|f|^p\mathrm{d}\mathbf{m}_{\mathcal{T}(\mathbf{b})},
\end{align*}
where $(\star)$ follows from Equations (\ref{eq_tree_cutoff}) and (\ref{eq_tree_num}), $(\dagger)$ follows from the fact that $\mathcal{N}(\widetilde{K})\subseteq B_{\mathcal{T}(\mathbf{b})}(x,8r)$, and $(\diamond)$ follows from Equation (\ref{eq_tree_num}) and the fact that $\mathbf{m}_{\mathcal{T}(\mathbf{b})}(\mathcal{N}(\widetilde{K}))\asymp V_\mathbf{b}(r)$.
\end{proof}

\section{Construction of Laakso-type spaces}\label{sec_Laa_metric}

In this section, we introduce a Laakso-type space along with a geodesic metric. We will give a canonical construction of geodesics, which will play a crucial role in the subsequent analysis.

Let $\mathcal{P}(\mathbf{g},\mathbf{b})=\mathcal{U}(\mathbf{g})\times\mathcal{T}(\mathbf{b})$ be endowed with a metric $\mathbf{d}_{\mathcal{P}(\mathbf{g},\mathbf{b})}$ and a measure $\mathbf{m}_{\mathcal{P}(\mathbf{g},\mathbf{b})}$ naturally inherited from $(\mathcal{U}(\mathbf{g}),\mathbf{d}_{\mathcal{U}(\mathbf{g})},\mathbf{m}_{\mathcal{U}(\mathbf{g})})$ and $(\mathcal{T}(\mathbf{b}),\mathbf{d}_{\mathcal{T}(\mathbf{b})},\mathbf{m}_{\mathcal{T}(\mathbf{b})})$, that is, $\mathbf{m}_{\mathcal{P}(\mathbf{g},\mathbf{b})}=\mathbf{m}_{\mathcal{U}(\mathbf{g})}\times\mathbf{m}_{\mathcal{T}(\mathbf{b})}$ and we use the convention
$$\mathbf{d}_{\mathcal{P}(\mathbf{g},\mathbf{b})}((\mathbf{u}^{(1)},\mathbf{t}^{(1)}),(\mathbf{u}^{(2)},\mathbf{t}^{(2)}))=\max\left\{\mathbf{d}_{\mathcal{U}(\mathbf{g})}(\mathbf{u}^{(1)},\mathbf{u}^{(2)}),\mathbf{d}_{\mathcal{T}(\mathbf{b})}(\mathbf{t}^{(1)},\mathbf{t}^{(2)})\right\}.$$
Then there exists some positive constant $C$ depending only on $\sup_\mathbb{Z}\mathbf{g}$ and $\sup_\mathbb{Z}\mathbf{b}$ such that
\begin{align*}
\frac{1}{C}V_\mathbf{g}(r)V_{\mathbf{b}}(r)\le&\mathbf{m}_{\mathcal{P}(\mathbf{g},\mathbf{b})}(B_{{\mathcal{P}(\mathbf{g},\mathbf{b})}}((\mathbf{u},\mathbf{t}),r)\le CV_\mathbf{g}(r)V_{\mathbf{b}}(r)\\
&\hspace{40pt}\text{ for any }(\mathbf{u},\mathbf{t})\in\mathcal{P}(\mathbf{g},\mathbf{b}),r>0.
\end{align*}

For any $n\in\mathbb{Z}$, we define the set of level-$n$ wormholes as
\begin{center}
$\mathcal{W}_n=$ the set of all centers of $n$-cells in $\mathcal{T}(\mathbf{b})$.
\end{center}
Any point in $\mathcal{W}_n$ is called a level-$n$ wormhole. It is obvious that $\mathcal{W}_n\cap\mathcal{W}_m=\emptyset$ for any $n\ne m$. We define a relation $R_\mathcal{L}$ on $\mathcal{P}(\mathbf{g},\mathbf{b})$ as follows. We say $(\mathbf{u}^{(1)},\mathbf{t}^{(1)}), (\mathbf{u}^{(2)},\mathbf{t}^{(2)})\in\mathcal{P}(\mathbf{g},\mathbf{b})$ are $R_\mathcal{L}$-related if, either $(\mathbf{u}^{(1)},\mathbf{t}^{(1)})=(\mathbf{u}^{(2)},\mathbf{t}^{(2)})$, or
\begin{center}
$\mathbf{u}^{(1)}|_{\mathbb{Z}\backslash\{n\}}=\mathbf{u}^{(2)}|_{\mathbb{Z}\backslash\{n\}}$ and $\mathbf{t}^{(1)}=\mathbf{t}^{(2)}\in\mathcal{W}_n$.
\end{center}
Then $R_\mathcal{L}$ is obviously an equivalence relation. Let $\mathcal{L}(\mathbf{g},\mathbf{b})=\mathcal{P}(\mathbf{g},\mathbf{b})/R_\mathcal{L}$ be the quotient space, called a Laakso-type space. Let $\mathcal{Q}:\mathcal{P}(\mathbf{g},\mathbf{b})\to\mathcal{L}(\mathbf{g},\mathbf{b})$, $\mathbf{p}\mapsto[\mathbf{p}]$ be the quotient map, where $[\mathbf{p}]$ is the equivalence class containing $\mathbf{p}$. Since $[(\mathbf{u}^{(1)},\mathbf{t}^{(1)})]=[(\mathbf{u}^{(2)},\mathbf{t}^{(2)})]$ always implies $\mathbf{t}^{(1)}=\mathbf{t}^{(2)}$, we have the map $\pi^\mathcal{T}:\mathcal{L}(\mathbf{g},\mathbf{b})\to\mathcal{T}(\mathbf{b})$, $[(\mathbf{u},\mathbf{t})]\mapsto\mathbf{t}$ is well-defined.

We endow the quotient space $\mathcal{L}(\mathbf{g},\mathbf{b})$ with the quotient topology from $\mathcal{P}(\mathbf{g},\mathbf{b})$. We say that $\gamma:[0,1]\to\mathcal{L}(\mathbf{g},\mathbf{b})$ is a path connecting $x,y\in\mathcal{L}(\mathbf{g},\mathbf{b})$ if $\gamma(0)=x$, $\gamma(1)=y$, and $\gamma$ is continuous. For convenience, we also say the image $\gamma([0,1])$ is a path.

We introduce $\mathbf{d}_{\mathcal{L}(\mathbf{g},\mathbf{b})}$ as follows. For any $x,y\in\mathcal{L}(\mathbf{g},\mathbf{b})$, let
$$\mathbf{d}_{\mathcal{L}(\mathbf{g},\mathbf{b})}\left(x,y\right)=\inf\left\{\mathcal{H}^{1}(\Gamma):\Gamma\subseteq\mathcal{P}(\mathbf{g},\mathbf{b})\text{ such that }\mathcal{Q}(\Gamma)\text{ is a path connecting }x,y\right\},$$
where $\mathcal{H}^1$ is the $1$-dimensional Hausdorff measure on $(\mathcal{P}(\mathbf{g},\mathbf{b}),\mathbf{d}_{\mathcal{P}(\mathbf{g},\mathbf{b})})$. It is easy to see that for any $[(\mathbf{u},\mathbf{t})], [(\mathbf{v},\mathbf{s})]\in\mathcal{L}(\mathbf{g},\mathbf{b})$
$$\mathbf{d}_{\mathcal{T}(\mathbf{b})}(\mathbf{t},\mathbf{s})\le\mathbf{d}_{\mathcal{L}(\mathbf{g},\mathbf{b})}([(\mathbf{u},\mathbf{t})],[(\mathbf{v},\mathbf{s})])\le+\infty.$$

In the following two results, we prove that $(\mathcal{L}(\mathbf{g},\mathbf{b}),\mathbf{d}_{\mathcal{L}(\mathbf{g},\mathbf{b})})$ is indeed a geodesic metric space and give a characterization of geodesics. These two results are similar to \cite[PROPOSITION 1.1, PROPOSITION 1.2]{Laa00}, where the tree is the unit interval and the ultrametric space is a Cantor set. The results in \cite{Laa00} were given without a formal proof and a detailed proof was recently given in \cite[Section 2]{Cap24}. We use a similar argument to that in \cite{Cap24}.

Firstly, we show that any two points in $\mathcal{L}(\mathbf{g},\mathbf{b})$ can be connected by a path $\mathcal{Q}(\Gamma)$, where $\Gamma\subseteq\mathcal{P}(\mathbf{g},\mathbf{b})$ is the union of \emph{countably}\footnote{Here countable refers to finite or countably infinite.} many geodesics in $\mathcal{T}(\mathbf{b})$. Here we do not require the path $\mathcal{Q}(\Gamma)$ to have minimal length.

\begin{lemma}\label{lem_Laa_path}
For any distinct $x=[(\mathbf{u}^{(1)},\mathbf{t}^{(1)})],y=[(\mathbf{u}^{(2)},\mathbf{t}^{(2)})]\in\mathcal{L}(\mathbf{g},\mathbf{b})$, there exists $\Gamma\subseteq\mathcal{P}(\mathbf{g},\mathbf{b})$ with
$$\Gamma=\bigcup_{\mathbf{u}\in U}\left(\{\mathbf{u}\}\times\gamma^{(\mathbf{u})}\right),$$
where $U\subseteq\mathcal{U}(\mathbf{g})$ is a countable subset and $\gamma^{(\mathbf{u})}$ is a geodesic in $(\mathcal{T}(\mathbf{b}),\mathbf{d}_{\mathcal{T}(\mathbf{b})})$ for any $\mathbf{u}\in U$, such that $\mathcal{Q}(\Gamma)$ is a path connecting $x,y$ and $\mathcal{H}^1(\Gamma)=\sum_{\mathbf{u}\in U}\lambda_{\mathcal{T}(\mathbf{b})}(\gamma^{(\mathbf{u})})<+\infty$.
\end{lemma}

\begin{proof}
We fix two points $(\mathbf{u}^{(1)},\mathbf{t}^{(1)})$ and $(\mathbf{u}^{(2)},\mathbf{t}^{(2)})$ from the equivalence classes $[(\mathbf{u}^{(1)},\mathbf{t}^{(1)})]$ and $[(\mathbf{u}^{(2)},\mathbf{t}^{(2)})]$. If $\mathbf{u}^{(1)}=\mathbf{u}^{(2)}$, let $\Gamma=\{\mathbf{u}^{(1)}\}\times[\mathbf{t}^{(1)},\mathbf{t}^{(2)}]$, then $\mathcal{Q}(\Gamma)$ is a path connecting $x,y$, and $\mathcal{H}^1(\Gamma)=\mathbf{d}_{\mathcal{T}(\mathbf{b})}(\mathbf{t}^{(1)},\mathbf{t}^{(2)})<+\infty$. Hence we may assume that $\mathbf{u}^{(1)}\ne\mathbf{u}^{(2)}$.

Let $I=\left\{k\in\mathbb{Z}:\mathbf{u}^{(1)}(k)\ne\mathbf{u}^{(2)}(k)\right\}$. Since $\mathbf{u}^{(1)}(k)\to0$, $\mathbf{u}^{(2)}(k)\to0$ as $k\to+\infty$, we have $I\subseteq\mathbb{Z}$ is bounded from above. Write $I=\{n_1,n_2,\ldots:n_1>n_2>\ldots\}$. Roughly speaking, $I$ is the set of all the levels of wormholes needed to jump to go from $x$ to $y$, our proof is to jump through wormholes sorted by level: $n_1,n_2,\ldots$.

Let $\mathbf{s}^{(0)}=\mathbf{t}^{(1)}$, $\mathbf{v}^{(0)}=\mathbf{u}^{(1)}$. Let $\mathbf{s}^{(1)}\in\mathcal{T}(\mathbf{b})$ be a nearest level-$n_1$ wormhole to $\mathbf{s}^{(0)}$\footnote{Recall that each $n$-cell has diameter $2^n$ in $(\mathcal{T}(\mathbf{b}),\mathbf{d}_{\mathcal{T}(\mathbf{b})})$ and each level-$n$ wormhole is the center of an $n$-cell.}, then $\mathbf{d}_{\mathcal{T}(\mathbf{b})}(\mathbf{s}^{(0)},\mathbf{s}^{(1)})\le2^{n_1-1}$. Let $\mathbf{v}^{(1)}\in\mathcal{U}(\mathbf{g})$ be given by $\mathbf{v}^{(1)}|_{\mathbb{Z}\backslash\{n_1\}}=\mathbf{v}^{(0)}|_{\mathbb{Z}\backslash\{n_1\}}$ and $\mathbf{v}^{(1)}(n_1)=\mathbf{u}^{(2)}(n_1)$, then $[(\mathbf{v}^{(0)},\mathbf{s}^{(1)})]=[(\mathbf{v}^{(1)},\mathbf{s}^{(1)})]$.

Assume we have constructed $\mathbf{s}^{(1)},\ldots,\mathbf{s}^{(l)}$, $\mathbf{v}^{(1)},\ldots,\mathbf{v}^{(l)}$. Let $\mathbf{s}^{(l+1)}\in\mathcal{T}(\mathbf{b})$ be a nearest level-$n_{l+1}$ wormhole to $\mathbf{s}^{(l)}$, then $\mathbf{d}_{\mathcal{T}(\mathbf{b})}(\mathbf{s}^{(l)},\mathbf{s}^{(l+1)})\le2^{n_{l+1}-1}$. Let $\mathbf{v}^{(l+1)}\in\mathcal{U}(\mathbf{g})$ be given by $\mathbf{v}^{(l+1)}|_{\mathbb{Z}\backslash\{n_{l+1}\}}=\mathbf{v}^{(l)}|_{\mathbb{Z}\backslash\{n_{l+1}\}}$ and $\mathbf{v}^{(l+1)}(n_{l+1})=\mathbf{u}^{(2)}(n_{l+1})$, then $[(\mathbf{v}^{(l)},\mathbf{s}^{(l+1)})]=[(\mathbf{v}^{(l+1)},\mathbf{s}^{(l+1)})]$. Since $\mathbf{v}^{(0)}=\mathbf{u}^{(1)}$ and $\mathbf{v}^{(l+1)}$ is obtained by replacing the $n_{l+1}$-th term of $\mathbf{v}^{(l)}$ by $\mathbf{u}^{(2)}(n_{l+1})$, we have $\mathbf{d}_{\mathcal{U}(\mathbf{g})}(\mathbf{v}^{(l)},\mathbf{u}^{(2)})\le2^{n_{l+1}}$.

If $I$ is a finite set, then $\mathbf{v}^{(\# I)}=\mathbf{u}^{(2)}$. Let
$$\Gamma=\left(\bigcup_{k=1}^{\# I}\{\mathbf{v}^{(k-1)}\}\times[\mathbf{s}^{(k-1)},\mathbf{s}^{(k)}]\right)\cup\left(\{\mathbf{u}^{(2)}\}\times[\mathbf{s}^{(\# I)},\mathbf{t}^{(2)}]\right),$$
then $\mathcal{Q}(\Gamma)$ is a path connecting $x,y$, and
$$\mathcal{H}^1(\Gamma)\le\sum_{k=1}^{\# I}2^{n_k-1}+\mathbf{d}_{\mathcal{T}(\mathbf{b})}(\mathbf{s}^{(\#I)},\mathbf{t}^{(2)})<+\infty.$$

If $I$ is an infinite set, then we obtain infinite sequences $\{\mathbf{s}^{(l)}\}$, $\{\mathbf{v}^{(l)}\}$. Since\\
\noindent $\mathbf{d}_{\mathcal{T}(\mathbf{b})}(\mathbf{s}^{(l)},\mathbf{s}^{(l+1)})\le2^{n_{l+1}-1}$, we have $\{\mathbf{s}^{(l)}\}$ converges to some $\mathbf{s}^{(\infty)}$ in $(\mathcal{T}(\mathbf{b}),\mathbf{d}_{\mathcal{T}(\mathbf{b})})$. Since $\mathbf{d}_{\mathcal{U}(\mathbf{g})}(\mathbf{v}^{(l)},\mathbf{u}^{(2)})\le2^{n_{l+1}}$, we have $\{\mathbf{v}^{(l)}\}$ converges to $\mathbf{u}^{(2)}$ in $(\mathcal{U}(\mathbf{g}),\mathbf{d}_{\mathcal{U}(\mathbf{g})})$. Let
$$\Gamma=\left(\bigcup_{k=1}^{+\infty}\{\mathbf{v}^{(k-1)}\}\times[\mathbf{s}^{(k-1)},\mathbf{s}^{(k)}]\right)\cup\left(\{\mathbf{u}^{(2)}\}\times[\mathbf{s}^{(\infty)},\mathbf{t}^{(2)}]\right),$$
then $\mathcal{Q}(\Gamma)$ is a path connecting $x,y$, and
$$\mathcal{H}^1(\Gamma)\le\sum_{k=1}^{+\infty}2^{n_k-1}+\mathbf{d}_{\mathcal{T}(\mathbf{b})}(\mathbf{s}^{(\infty)},\mathbf{t}^{(2)})<+\infty.$$
\end{proof}

Secondly, we show that any two points in $\mathcal{L}(\mathbf{g},\mathbf{b})$ can be connected by a geodesic which can constructed in a canonical way. We will see that such geodesic is either monotone, or has one or two inversions.

\begin{proposition}\label{prop_Laa_geo}
$(\mathcal{L}(\mathbf{g},\mathbf{b}),\mathbf{d}_{\mathcal{L}(\mathbf{g},\mathbf{b})})$ is a geodesic space. Moreover, for any distinct $[(\mathbf{u},\mathbf{t})]$, $[(\mathbf{v},\mathbf{s})]\in\mathcal{L}(\mathbf{g},\mathbf{b})$, there exists a geodesic connecting $[(\mathbf{u},\mathbf{t})]$, $[(\mathbf{v},\mathbf{s})]$, which falls into one of the following three cases.

\begin{enumerate}[label=(GEO\arabic*),ref=(GEO\arabic*)]
\item\label{item_geo_mono} There exists a geodesic $\gamma$ connecting $[(\mathbf{u},\mathbf{t})]$, $[(\mathbf{v},\mathbf{s})]$ with length $\mathbf{d}_{\mathcal{T}(\mathbf{b})}(\mathbf{t},\mathbf{s})$, we say such a geodesic $\gamma$ is \emph{monotone}, see Figure \ref{fig_mono}.
\item\label{item_geo_oneinv} There exist $[(\mathbf{u}^{(1)},\mathbf{t}^{(1)})]\in\mathcal{L}(\mathbf{g},\mathbf{b})$ with $\mathbf{t}^{(1)}\not\in[\mathbf{t},\mathbf{s}]$, and two monotone geodesics $\gamma^{(1)}$ and $\gamma^{(2)}$, with $\gamma^{(1)}$ connecting $[(\mathbf{u},\mathbf{t})]$, $[(\mathbf{u}^{(1)},\mathbf{t}^{(1)})]$, and $\gamma^{(2)}$ connecting $[(\mathbf{u}^{(1)},\mathbf{t}^{(1)})]$, $[(\mathbf{v},\mathbf{s})]$, as in \ref{item_geo_mono}, such that $\gamma^{(1)}\cup\gamma^{(2)}$ is a geodesic connecting $[(\mathbf{u},\mathbf{t})]$, $[(\mathbf{v},\mathbf{s})]$ with length
$$\mathbf{d}_{\mathcal{T}(\mathbf{b})}(\mathbf{t},\mathbf{s})+2\mathbf{d}_{\mathcal{T}(\mathbf{b})}(\mathbf{t}^{(1)},[\mathbf{t},\mathbf{s}]),$$
we say such a geodesic has \emph{one} inversion, see Figure \ref{fig_oneinv}.
\item\label{item_geo_twoinv} There exist $[(\mathbf{u}^{(1)},\mathbf{t}^{(1)})]$, $[(\mathbf{u}^{(2)},\mathbf{t}^{(2)})]\in\mathcal{L}(\mathbf{g},\mathbf{b})$ with $\mathbf{t}^{(1)},\mathbf{t}^{(2)}\not\in[\mathbf{t},\mathbf{s}]$, and three monotone geodesics $\gamma^{(1)}$, $\gamma^{(2)}$ and $\gamma^{(3)}$, with $\gamma^{(1)}$ connecting $[(\mathbf{u},\mathbf{t})]$, $[(\mathbf{u}^{(1)},\mathbf{t}^{(1)})]$, $\gamma^{(2)}$ connecting $[(\mathbf{u}^{(1)},\mathbf{t}^{(1)})]$, $[(\mathbf{u}^{(2)},\mathbf{t}^{(2)})]$, and $\gamma^{(3)}$ connecting $[(\mathbf{u}^{(2)},\mathbf{t}^{(2)})]$, $[(\mathbf{v},\mathbf{s})]$, as in\\
\noindent \ref{item_geo_mono}, such that $\gamma^{(1)}\cup\gamma^{(2)}\cup\gamma^{(3)}$ is a geodesic connecting $[(\mathbf{u},\mathbf{t})]$, $[(\mathbf{v},\mathbf{s})]$ with length
$$\mathbf{d}_{\mathcal{T}(\mathbf{b})}(\mathbf{t},\mathbf{s})+2\mathbf{d}_{\mathcal{T}(\mathbf{b})}(\mathbf{t}^{(1)},[\mathbf{t},\mathbf{s}])+2\mathbf{d}_{\mathcal{T}(\mathbf{b})}(\mathbf{t}^{(2)},[\mathbf{t},\mathbf{s}]),$$
we say such a geodesic has \emph{two} inversions, see Figure \ref{fig_twoinv}.
\end{enumerate}
\end{proposition}

\begin{remark}
We will provide explicit choices of $\mathbf{t}^{(1)}$ and $\mathbf{t}^{(2)}$ in \ref{item_geo_oneinv} and \ref{item_geo_twoinv}, as follows.
\begin{itemize}
\item In \ref{item_geo_oneinv}, we have
$$n^*=\max\{n\in\mathbb{Z}:\mathbf{u}(n)\ne\mathbf{v}(n)\text{ and }[\mathbf{t},\mathbf{s}]\cap\mathcal{W}_n=\emptyset\},$$
is well-defined, and $\mathbf{t}^{(1)}\in\mathcal{W}_{n^*}$ satisfies
$$\mathbf{d}_{\mathcal{T}(\mathbf{b})}(\mathbf{t}^{(1)},[\mathbf{t},\mathbf{s}])=\min_{\mathbf{r}\in\mathcal{W}_{n^*}}\mathbf{d}_{\mathcal{T}(\mathbf{b})}(\mathbf{r},[\mathbf{t},\mathbf{s}]).$$
\item In \ref{item_geo_twoinv}, we have
\begin{align*}
n^*&=\max\{n\in\mathbb{Z}:\mathbf{u}(n)\ne\mathbf{v}(n)\text{ and }[\mathbf{t},\mathbf{s}]\cap\mathcal{W}_n=\emptyset\},\\
n^{**}&=\max\{n\in\mathbb{Z}:n<n^*,\mathbf{u}(n)\ne\mathbf{v}(n)\},
\end{align*}
are well-defined, and $\mathbf{t}^{(1)}\in\mathcal{W}_{n^*}$, $\mathbf{t}^{(2)}\in\mathcal{W}_{n^{**}}$ satisfy
\begin{align*}
\mathbf{d}_{\mathcal{T}(\mathbf{b})}(\mathbf{t}^{(1)},[\mathbf{t},\mathbf{s}])&=\min_{\mathbf{r}\in\mathcal{W}_{n^*}}\mathbf{d}_{\mathcal{T}(\mathbf{b})}(\mathbf{r},[\mathbf{t},\mathbf{s}]),\\
\mathbf{d}_{\mathcal{T}(\mathbf{b})}(\mathbf{t}^{(2)},[\mathbf{t},\mathbf{s}])&=\min_{\mathbf{r}\in\mathcal{W}_{n^{**}}}\mathbf{d}_{\mathcal{T}(\mathbf{b})}(\mathbf{r},[\mathbf{t},\mathbf{s}]).
\end{align*}
\end{itemize}
\end{remark}

\begin{proof}
We fix two points $(\mathbf{u},\mathbf{t})$, $(\mathbf{v},\mathbf{s})$ from the equivalence classes $[(\mathbf{u},\mathbf{t})]$, $[(\mathbf{v},\mathbf{s})]$. If $\mathbf{u}=\mathbf{v}$, let $\Gamma=\{\mathbf{u}\}\times[\mathbf{t},\mathbf{s}]$, then $\mathcal{Q}(\Gamma)$ is obviously a geodesic connecting $[(\mathbf{u},\mathbf{t})]$, $[(\mathbf{v},\mathbf{s})]$, and $\mathcal{H}^1(\Gamma)=\mathbf{d}_{\mathcal{T}(\mathbf{b})}(\mathbf{t},\mathbf{s})=\mathbf{d}_{\mathcal{L}(\mathbf{g},\mathbf{b})}([(\mathbf{u},\mathbf{t})],[(\mathbf{v},\mathbf{s})])$, we say such a geodesic is \emph{monotone}. Hence we may assume that $\mathbf{u}\ne\mathbf{v}$.

Let $I=\left\{k\in\mathbb{Z}:\mathbf{u}(k)\ne\mathbf{v}(k)\right\}$. Since $\mathbf{u}(k)\to0$, $\mathbf{v}(k)\to0$ as $k\to+\infty$, we have $I\subseteq\mathbb{Z}$ is bounded from above. Write $I=\{n_1,n_2,\ldots:n_1>n_2>\ldots\}$. Roughly speaking, any path connecting $[(\mathbf{u},\mathbf{t})]$, $[(\mathbf{v},\mathbf{s})]$ needs to jump through wormholes with all the levels in $I$. We only have the following two cases.
\begin{enumerate}[label=(\arabic*)]
\item $[\mathbf{t},\mathbf{s}]\cap\mathcal{W}_{n}\ne\emptyset$ for \emph{any} $n\in I$, or $[\mathbf{t},\mathbf{s}]$ \emph{does} contain all the ``necessary" wormholes.
\item $[\mathbf{t},\mathbf{s}]\cap\mathcal{W}_{n}=\emptyset$ for \emph{some} $n\in I$, or $[\mathbf{t},\mathbf{s}]$ does \emph{not} contain all the ``necessary" wormholes.
\end{enumerate}

For (1), we construct a path connecting $[(\mathbf{u},\mathbf{t})]$, $[(\mathbf{v},\mathbf{s})]$ with length $\mathbf{d}_{\mathcal{T}(\mathbf{b})}(\mathbf{t},\mathbf{s})$, which would imply this path is indeed a geodesic connecting $[(\mathbf{u},\mathbf{t})]$, $[(\mathbf{v},\mathbf{s})]$ and\\
\noindent $\mathbf{d}_{\mathcal{L}(\mathbf{g},\mathbf{b})}([\mathbf{u},\mathbf{t}],[\mathbf{v},\mathbf{s}])=\mathbf{d}_{\mathcal{T}(\mathbf{b})}(\mathbf{t},\mathbf{s})$. If $\# I=1$, let $\mathbf{t}^{(1)}\in[\mathbf{t},\mathbf{s}]\cap\mathcal{W}_{n_1}$ and $\Gamma=(\{\mathbf{u}\}\times[\mathbf{t},\mathbf{t}^{(1)}])\cup(\{\mathbf{v}\}\times[\mathbf{t}^{(1)},\mathbf{s}])$, then $\mathcal{Q}(\Gamma)$ is a path connecting $[(\mathbf{u},\mathbf{t})]$, $[(\mathbf{v},\mathbf{s})]$ with length $\mathbf{d}_{\mathcal{T}(\mathbf{b})}(\mathbf{t},\mathbf{s})$. Assume that $\# I\ge2$. Take $\mathbf{t}^{(1)}\in[\mathbf{t},\mathbf{s}]\cap\mathcal{W}_{n_1}$ and $\mathbf{t}^{(2)}\in[\mathbf{t},\mathbf{s}]\cap\mathcal{W}_{n_2}$. By interchanging the roles of $\mathbf{t}$ and $\mathbf{s}$, we may assume that $[\mathbf{t},\mathbf{t}^{(1)}]\cap[\mathbf{t}^{(2)},\mathbf{s}]=\emptyset$. Since $n_1>n_2$, we have $[\mathbf{t}^{(1)},\mathbf{t}^{(2)}]\cap\mathcal{W}_n\ne\emptyset$ for any $n>n_2$.

If $\#I=+\infty$, for any $l\ge3$, there exists $\mathbf{t}^{(l)}\in[\mathbf{t}^{(1)},\mathbf{t}^{(2)}]\cap\mathcal{W}_{n_l}$ satisfying that $\mathbf{t}^{(l+1)}\in[\mathbf{t}^{(l)},\mathbf{t}^{(2)}]$, $\mathbf{d}_{\mathcal{T}(\mathbf{b})}(\mathbf{t}^{(l)},\mathbf{t}^{(l+1)})=2^{n_{l+1}-1}$, and $\mathbf{d}_{\mathcal{T}(\mathbf{b})}(\mathbf{t}^{(1)},\mathbf{t}^{(3)})=2^{n_3-1}$.\footnote{This follows from the fact that for any $n,m\in\mathbb{Z}$ with $n>m$, for any $n$-cell, there exist $\mathbf{b}(n)$ $m$-cells whose centers have the distance $2^{m-1}$ to the center of the $n$-cell.} Hence $\{\mathbf{t}^{(l)}\}$ converges to some $\mathbf{t}^{(\infty)}\in[\mathbf{t}^{(1)},\mathbf{t}^{(2)}]$ in $(\mathcal{T}(\mathbf{b}),\mathbf{d}_{\mathcal{T}(\mathbf{b})})$, see Figure \ref{fig_mono}.

\begin{figure}[ht]
\centering
\begin{tikzpicture}
\draw[thick] (0,0)--(5,0);

\draw[very thick] (0,-0.05)--(0,0.05);
\draw[very thick] (5,-0.05)--(5,0.05);
\draw[very thick] (1,-0.05)--(1,0.05);
\draw[very thick] (4,-0.05)--(4,0.05);
\draw[very thick] (1.8,-0.05)--(1.8,0.05);
\draw[very thick] (2.3,-0.05)--(2.3,0.05);
\draw[very thick] (3.7,-0.05)--(3.7,0.05);

\node[below] at (0,0) {$\mathbf{t}$};
\node[below] at (5,0) {$\mathbf{s}$};
\node[below] at (1,0) {$\mathbf{t}^{(1)}$};
\node[below] at (4,0) {$\mathbf{t}^{(2)}$};
\node[above] at (1.8,0) {\scriptsize{$\mathbf{t}^{(3)}$}};
\node[above] at (2.3,0) {\scriptsize{$\mathbf{t}^{(4)}$}};
\node[above] at (3,0) {\scriptsize{\ldots}};
\node[above] at (3.7,0) {\scriptsize{$\mathbf{t}^{(\infty)}$}};

\end{tikzpicture}
\caption{A monotone geodesic}\label{fig_mono}
\end{figure}
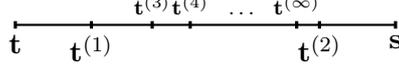

Let $\mathbf{u}^{(1)}\in\mathcal{U}(\mathbf{g})$ be given by $\mathbf{u}^{(1)}|_{\mathbb{Z}\backslash\{n_1\}}=\mathbf{u}|_{\mathbb{Z}\backslash\{n_1\}}$ and $\mathbf{u}^{(1)}(n_1)=\mathbf{v}(n_1)$. Let $\mathbf{u}^{(2)}\in\mathcal{U}(\mathbf{g})$ be given by $\mathbf{u}^{(2)}|_{\mathbb{Z}\backslash\{n_2\}}=\mathbf{v}|_{\mathbb{Z}\backslash\{n_2\}}$ and $\mathbf{u}^{(2)}(n_2)=\mathbf{u}(n_2)$. Let $\mathbf{u}^{(3)}\in\mathcal{U}(\mathbf{g})$ be given by $\mathbf{u}^{(3)}|_{\mathbb{Z}\backslash\{n_3\}}=\mathbf{u}^{(1)}|_{\mathbb{Z}\backslash\{n_3\}}$ and $\mathbf{u}^{(3)}(n_3)=\mathbf{v}(n_3)$. For any $l\ge3$, let $\mathbf{u}^{(l+1)}\in\mathcal{U}(\mathbf{g})$ be given by $\mathbf{u}^{(l+1)}|_{\mathbb{Z}\backslash\{n_{l+1}\}}=\mathbf{u}^{(l)}|_{\mathbb{Z}\backslash\{n_{l+1}\}}$ and $\mathbf{u}^{(l+1)}(n_{l+1})=\mathbf{v}(n_{l+1})$. Then $\{\mathbf{u}^{(l)}\}$ converges to $\mathbf{u}^{(2)}$ in $(\mathcal{U}(\mathbf{g}),\mathbf{d}_{\mathcal{U}(\mathbf{g})})$. Let
\begin{align*}
&\Gamma=\left(\{\mathbf{u}\}\times[\mathbf{t},\mathbf{t}^{(1)}]\right)\cup\left(\{\mathbf{u}^{(1)}\}\times[\mathbf{t}^{(1)},\mathbf{t}^{(3)}]\right)\\
&\cup\left(\bigcup_{l=3}^{+\infty}\{\mathbf{u}^{(l)}\}\times[\mathbf{t}^{(l)},\mathbf{t}^{(l+1)}]\right)\cup\left(\{\mathbf{u}^{(2)}\}\times[\mathbf{t}^{(\infty)},\mathbf{t}^{(2)}]\right)\cup\left(\{\mathbf{v}\}\times[\mathbf{t}^{(2)},\mathbf{s}]\right),
\end{align*}
then $\mathcal{Q}(\Gamma)$ is a path connecting $[(\mathbf{u},\mathbf{t})]$, $[(\mathbf{v},\mathbf{s})]$ with length $\mathbf{d}_{\mathcal{T}(\mathbf{b})}(\mathbf{t},\mathbf{s})$, we say $\mathcal{Q}(\Gamma)$ is \emph{monotone}. If $\#I<+\infty$, we construct $\Gamma$ using a similar finite union, which also gives the desired result, we also say $\mathcal{Q}(\Gamma)$ is \emph{monotone}.

For (2), $[\mathbf{t},\mathbf{s}]\cap\mathcal{W}_{n}=\emptyset$ for \emph{some} $n\in I$, let $n^*=n_l$ be the maximum of such $n$, take $\mathbf{t}^{(1)}\in\mathcal{W}_{n^*}$, which is closest to $[\mathbf{t},\mathbf{s}]$, while such $\mathbf{t}^{(1)}$ may not be unique, only finitely many exist, and we choose one arbitrarily. Let $\mathbf{c}=c(\mathbf{t}^{(1)},\mathbf{t},\mathbf{s})$, then $\mathbf{t}^{(1)}\ne\mathbf{c}$, see Figure \ref{fig_oneinv}.

\begin{figure}[ht]
\centering
\begin{tikzpicture}[scale=0.7]

\draw[thick] (0,0)--(-2,0);
\draw[thick] (0,0)--(1,1.732);
\draw[thick] (0,0)--(1,-1.732);

\filldraw[fill=black] (-2,0) circle (0.05);
\filldraw[fill=black] (0,0) circle (0.05);
\filldraw[fill=black] (1,1.732) circle (0.05);
\filldraw[fill=black] (1,-1.732) circle (0.05);

\node[below] at (-2,0) {$\mathbf{t}^{(1)}$};
\node[right] at (1,1.732) {$\mathbf{t}$};
\node[right] at (1,-1.732) {$\mathbf{s}$};
\node[right] at (0,0) {$\mathbf{c}$};

\end{tikzpicture}
\caption{A geodesic with one inversion}\label{fig_oneinv}
\end{figure}
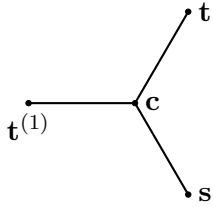

Firstly, if $([\mathbf{t}^{(1)},\mathbf{t}]\cup[\mathbf{t}^{(1)},\mathbf{s}])\cap\mathcal{W}_n\ne\emptyset$ for \emph{any} $n\in I$, let
\begin{align*}
I_1&=\left\{n\in I:[\mathbf{t},\mathbf{t}^{(1)}]\cap\mathcal{W}_n\ne\emptyset\right\},\\
I_2&=\left\{n\in I\backslash I_1:[\mathbf{t}^{(1)},\mathbf{s}]\cap\mathcal{W}_n\ne\emptyset\right\},
\end{align*}
then $I=I_1\sqcup I_2$. Let $\mathbf{u}^{(1)}\in\mathcal{U}(\mathbf{g})$ be given by $\mathbf{u}^{(1)}|_{\mathbb{Z}\backslash I_1}=\mathbf{u}|_{\mathbb{Z}\backslash I_1}$ and $\mathbf{u}^{(1)}|_{I_1}=\mathbf{v}|_{I_1}$, then by (1), there exist
\begin{itemize}
\item a monotone geodesic $\gamma^{(1)}$ connecting $[(\mathbf{u},\mathbf{t})]$, $[(\mathbf{u}^{(1)},\mathbf{t}^{(1)})]$ with length $\mathbf{d}_{\mathcal{T}(\mathbf{b})}(\mathbf{t},\mathbf{t}^{(1)})$,
\item a monotone geodesic $\gamma^{(2)}$ connecting $[(\mathbf{u}^{(1)},\mathbf{t}^{(1)})]$, $[(\mathbf{v},\mathbf{s})]$ with length $\mathbf{d}_{\mathcal{T}(\mathbf{b})}(\mathbf{t}^{(1)},\mathbf{s})$.
\end{itemize}
Hence $\gamma^{(1)}\cup\gamma^{(2)}$ is a path connecting $[(\mathbf{u},\mathbf{t})]$, $[(\mathbf{v},\mathbf{s})]$ with length
$$\mathbf{d}_{\mathcal{T}(\mathbf{b})}(\mathbf{t},\mathbf{t}^{(1)})+\mathbf{d}_{\mathcal{T}(\mathbf{b})}(\mathbf{t}^{(1)},\mathbf{s})=\mathbf{d}_{\mathcal{T}(\mathbf{b})}(\mathbf{t},\mathbf{s})+2\mathbf{d}_{\mathcal{T}(\mathbf{b})}(\mathbf{t}^{(1)},[\mathbf{t},\mathbf{s}]).$$
Since any path connecting $[(\mathbf{u},\mathbf{t})]$, $[(\mathbf{v},\mathbf{s})]$ needs to jump through at least one level-$n^*$ wormhole, we have
\begin{align*}
\mathbf{d}_{\mathcal{L}(\mathbf{g},\mathbf{b})}([(\mathbf{u},\mathbf{t})], [(\mathbf{v},\mathbf{s})])&\ge\mathbf{d}_{\mathcal{T}(\mathbf{b})}(\mathbf{t},\mathbf{s})+2\min_{\mathbf{r}\in\mathcal{W}_{n^*}}\mathbf{d}_{\mathcal{T}(\mathbf{b})}(\mathbf{r},[\mathbf{t},\mathbf{s}])\\
&=\mathbf{d}_{\mathcal{T}(\mathbf{b})}(\mathbf{t},\mathbf{s})+2\mathbf{d}_{\mathcal{T}(\mathbf{b})}(\mathbf{t}^{(1)},[\mathbf{t},\mathbf{s}]).
\end{align*}
Hence $\gamma^{(1)}\cup\gamma^{(2)}$ is a geodesic connecting $[(\mathbf{u},\mathbf{t})]$, $[(\mathbf{v},\mathbf{s})]$ and
$$\mathbf{d}_{\mathcal{L}(\mathbf{g},\mathbf{b})}([(\mathbf{u},\mathbf{t})], [(\mathbf{v},\mathbf{s})])=\mathbf{d}_{\mathcal{T}(\mathbf{b})}(\mathbf{t},\mathbf{s})+2\mathbf{d}_{\mathcal{T}(\mathbf{b})}(\mathbf{t}^{(1)},[\mathbf{t},\mathbf{s}]),$$
we say such a geodesic has \emph{one} inversion.

Secondly, if $([\mathbf{t}^{(1)},\mathbf{t}]\cup[\mathbf{t}^{(1)},\mathbf{s}])\cap\mathcal{W}_n=\emptyset$ for \emph{some} $n\in I$, let $n^{**}$ be the maximum of such $n$, then $n^{**}<n^*=n_l$. We claim that $n^{**}=n_{l+1}$. Otherwise, we have $([\mathbf{t}^{(1)},\mathbf{t}]\cup[\mathbf{t}^{(1)},\mathbf{s}])\cap\mathcal{W}_{n_{l+1}}\ne\emptyset$, take $\mathbf{t}^{(2)}\in([\mathbf{t}^{(1)},\mathbf{t}]\cup[\mathbf{t}^{(1)},\mathbf{s}])\cap\mathcal{W}_{n_{l+1}}$. Since $n_l>n_{l+1}$, we have $[\mathbf{t}^{(1)},\mathbf{t}^{(2)}]\cap\mathcal{W}_n\ne\emptyset$ for any $n>n_{l+1}$, in particular, $[\mathbf{t}^{(1)},\mathbf{t}]\cup[\mathbf{t}^{(1)},\mathbf{s}]\supseteq[\mathbf{t}^{(1)},\mathbf{t}^{(2)}]$ intersects $\mathcal{W}_n$ for any $n>n_{l+1}$, $([\mathbf{t}^{(1)},\mathbf{t}]\cup[\mathbf{t}^{(1)},\mathbf{s}])\cap\mathcal{W}_n\ne\emptyset$ for any $n\in I$, contradicting to our assumption.

Take $\mathbf{t}^{(2)}\in\mathcal{W}_{n^{**}}$, which is closest to $[\mathbf{t},\mathbf{s}]$. Let $\mathbf{c}^{(1)}=c(\mathbf{t},\mathbf{t}^{(1)},\mathbf{t}^{(2)})$, $\mathbf{c}^{(2)}=c(\mathbf{s},\mathbf{t}^{(1)},\mathbf{t}^{(2)})$. By assumption, we \emph{only} have the case $\mathbf{c}^{(1)},\mathbf{c}^{(2)}\in[\mathbf{t}^{(1)},\mathbf{t}^{(2)}]\backslash\{\mathbf{t}^{(1)},\mathbf{t}^{(2)}\}$, see Figure \ref{fig_twoinv}, however, we do \emph{not} exclude the cases $\mathbf{c}^{(1)}=\mathbf{c}^{(2)}$, $\mathbf{t}=\mathbf{c}^{(1)}$ or $\mathbf{s}=\mathbf{c}^{(2)}$. By interchanging the roles of $\mathbf{t}$ and $\mathbf{s}$, we may assume that $\mathbf{c}^{(1)}\subseteq[\mathbf{t}^{(1)},\mathbf{c}^{(2)}]$.

\begin{figure}[ht]
\centering
\begin{tikzpicture}[scale=1]

\draw[thick] (0,0)--(4,0);
\draw[thick] (1,0)--(1,2);
\draw[thick] (3,0)--(3,1);

\filldraw[fill=black] (0,0) circle (0.05);
\filldraw[fill=black] (4,0) circle (0.05);
\filldraw[fill=black] (1,0) circle (0.05);
\filldraw[fill=black] (1,2) circle (0.05);
\filldraw[fill=black] (3,0) circle (0.05);
\filldraw[fill=black] (3,1) circle (0.05);

\node[right] at (1,2) {$\mathbf{t}$};
\node[right] at (3,1) {$\mathbf{s}$};
\node[below] at (0,0) {$\mathbf{t}^{(1)}$};
\node[below] at (4,0) {$\mathbf{t}^{(2)}$};
\node[below] at (1,0) {$\mathbf{c}^{(1)}$};
\node[below] at (3,0) {$\mathbf{c}^{(2)}$};

\end{tikzpicture}
\caption{A geodesic with two inversions}\label{fig_twoinv}
\end{figure}
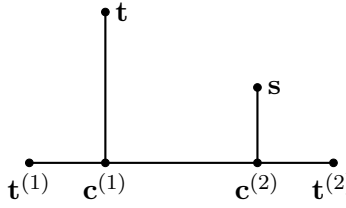

Since $\mathbf{t}^{(1)}\in\mathcal{W}_{n^*}=\mathcal{W}_{n_l}$, $\mathbf{t}^{(2)}\in\mathcal{W}_{n^{**}}=\mathcal{W}_{n_{l+1}}$ and $n_l>n_{l+1}$, we have $[\mathbf{t}^{(1)},\mathbf{t}^{(2)}]\cap\mathcal{W}_n\ne\emptyset$ for any $n>n_{l+1}$, which implies $([\mathbf{t},\mathbf{t}^{(1)}]\cup[\mathbf{t}^{(1)},\mathbf{t}^{(2)}]\cup[\mathbf{t}^{(2)},\mathbf{s}])\cap\mathcal{W}_n\ne\emptyset$ for any $n\in I$. Let
\begin{align*}
I_1&=\left\{n\in I:[\mathbf{t},\mathbf{t}^{(1)}]\cap\mathcal{W}_n\ne\emptyset\right\},\\
I_2&=\left\{n\in I\backslash I_1:[\mathbf{t}^{(1)},\mathbf{t}^{(2)}]\cap\mathcal{W}_n\ne\emptyset\right\},\\
I_3&=\left\{n\in I\backslash (I_1\sqcup I_2):[\mathbf{t}^{(2)},\mathbf{s}]\cap\mathcal{W}_n\ne\emptyset\right\},
\end{align*}
then $I=I_1\sqcup I_2\sqcup I_3$. Let $\mathbf{u}^{(1)}\in\mathcal{U}(\mathbf{g})$ be given by $\mathbf{u}^{(1)}|_{\mathbb{Z}\backslash I_1}=\mathbf{u}|_{\mathbb{Z}\backslash I_1}$ and $\mathbf{u}^{(1)}|_{I_1}=\mathbf{v}|_{I_1}$. Let $\mathbf{u}^{(2)}\in\mathcal{U}(\mathbf{g})$ be given by $\mathbf{u}^{(2)}|_{\mathbb{Z}\backslash I_2}=\mathbf{u}^{(1)}|_{\mathbb{Z}\backslash I_2}$ and $\mathbf{u}^{(2)}|_{I_2}=\mathbf{v}|_{I_2}$, then by (1), there exist
\begin{itemize}
\item a monotone geodesic $\gamma^{(1)}$ connecting $[(\mathbf{u},\mathbf{t})]$, $[(\mathbf{u}^{(1)},\mathbf{t}^{(1)})]$ with length $\mathbf{d}_{\mathcal{T}(\mathbf{b})}(\mathbf{t},\mathbf{t}^{(1)})$,
\item a monotone geodesic $\gamma^{(2)}$ connecting $[(\mathbf{u}^{(1)},\mathbf{t}^{(1)})]$, $[(\mathbf{u}^{(2)},\mathbf{t}^{(2)})]$\\
\noindent with length $\mathbf{d}_{\mathcal{T}(\mathbf{b})}(\mathbf{t}^{(1)},\mathbf{t}^{(2)})$,
\item a monotone geodesic $\gamma^{(3)}$ connecting $[(\mathbf{u}^{(2)},\mathbf{t}^{(2)})]$, $[(\mathbf{v},\mathbf{s})]$ with length $\mathbf{d}_{\mathcal{T}(\mathbf{b})}(\mathbf{t}^{(2)},\mathbf{s})$.
\end{itemize}
Hence $\gamma^{(1)}\cup\gamma^{(2)}\cup\gamma^{(3)}$ is a path connecting $[(\mathbf{u},\mathbf{t})]$, $[(\mathbf{v},\mathbf{s})]$ with length
\begin{align*}
&\mathbf{d}_{\mathcal{T}(\mathbf{b})}(\mathbf{t},\mathbf{t}^{(1)})+\mathbf{d}_{\mathcal{T}(\mathbf{b})}(\mathbf{t}^{(1)},\mathbf{t}^{(2)})+\mathbf{d}_{\mathcal{T}(\mathbf{b})}(\mathbf{t}^{(2)},\mathbf{s})\\
&=\mathbf{d}_{\mathcal{T}(\mathbf{b})}(\mathbf{t},\mathbf{s})+2\mathbf{d}_{\mathcal{T}(\mathbf{b})}(\mathbf{t}^{(1)},[\mathbf{t},
\mathbf{s}])+2\mathbf{d}_{\mathcal{T}(\mathbf{b})}(\mathbf{t}^{(2)},[\mathbf{t},\mathbf{s}]).
\end{align*}
Since any path connecting $[(\mathbf{u},\mathbf{t})]$, $[(\mathbf{v},\mathbf{s})]$ needs to jump through at least one level-$n^*$ wormhole and at least one level-$n^{**}$ wormhole, we have
\begin{align*}
\mathbf{d}_{\mathcal{L}(\mathbf{g},\mathbf{b})}([(\mathbf{u},\mathbf{t})], [(\mathbf{v},\mathbf{s})])&\ge\mathbf{d}_{\mathcal{T}(\mathbf{b})}(\mathbf{t},\mathbf{s})+2\min_{\mathbf{r}\in\mathcal{W}_{n^*}}\mathbf{d}_{\mathcal{T}(\mathbf{b})}(\mathbf{r},[\mathbf{t},\mathbf{s}])+2\min_{\mathbf{r}\in\mathcal{W}_{n^{**}}}\mathbf{d}_{\mathcal{T}(\mathbf{b})}(\mathbf{r},[\mathbf{t},\mathbf{s}])\\
&=\mathbf{d}_{\mathcal{T}(\mathbf{b})}(\mathbf{t},\mathbf{s})+2\mathbf{d}_{\mathcal{T}(\mathbf{b})}(\mathbf{t}^{(1)},[\mathbf{t},\mathbf{s}])+2\mathbf{d}_{\mathcal{T}(\mathbf{b})}(\mathbf{t}^{(2)},[\mathbf{t},\mathbf{s}]).
\end{align*}
Hence $\gamma^{(1)}\cup\gamma^{(2)}\cup\gamma^{(3)}$ is a geodesic connecting $[(\mathbf{u},\mathbf{t})]$, $[(\mathbf{v},\mathbf{s})]$ and
$$\mathbf{d}_{\mathcal{L}(\mathbf{g},\mathbf{b})}([(\mathbf{u},\mathbf{t})], [(\mathbf{v},\mathbf{s})])=\mathbf{d}_{\mathcal{T}(\mathbf{b})}(\mathbf{t},\mathbf{s})+2\mathbf{d}_{\mathcal{T}(\mathbf{b})}(\mathbf{t}^{(1)},[\mathbf{t},\mathbf{s}])+2\mathbf{d}_{\mathcal{T}(\mathbf{b})}(\mathbf{t}^{(2)},[\mathbf{t},\mathbf{s}]),$$
we say such a geodesic has \emph{two} inversions.
\end{proof}

From now on, all paths are assumed to be parametrized by arc length, and the length of a path $\gamma$ will be denoted by $l(\gamma)$. We list some easy results for later use.

\begin{lemma}\label{lem_onept}
\hspace{0em}
\begin{enumerate}[label=(\arabic*),ref=(\arabic*)]
\item \label{item_maxonept} Let $\mathbf{t}\in\mathcal{T}(\mathbf{b})$, $r>0$, and $n_0$ the integer satisfying $2^{n_0-1}\le r<2^{n_0}$, then $B_{\mathcal{T}(\mathbf{b})}(\mathbf{t},r)\cap\cup_{n\ge n_0+2}\mathcal{W}_n$ consists of at most one point.
\item\label{item_minonept} Let $\gamma$ be a geodesic in $\mathcal{L}(\mathbf{g},\mathbf{b})$ given by Proposition \ref{prop_Laa_geo}, and $n_0$ the integer satisfying $2^{n_0+2}\le l(\gamma)<2^{n_0+3}$, then for any $n\le n_0$, we have $\pi^\mathcal{T}(\gamma)\cap\mathcal{W}_n\ne\emptyset$.
\end{enumerate}

\end{lemma}

\begin{proof}
\ref{item_maxonept} Suppose there exist distinct $\mathbf{s}^{(1)},\mathbf{s}^{(2)}$ in this set, since $\mathbf{s}^{(1)},\mathbf{s}^{(2)}\in B_{\mathcal{T}(\mathbf{b})}(\mathbf{t},r)$, we have $\mathbf{d}_{\mathcal{T}(\mathbf{b})}(\mathbf{s}^{(1)},\mathbf{s}^{(2)})<2r<2^{n_0+1}$. However, since $\mathbf{s}^{(1)},\mathbf{s}^{(2)}\in\cup_{n\ge n_0+2}\mathcal{W}_n$ are distinct, we have $\mathbf{d}_{\mathcal{T}(\mathbf{b})}(\mathbf{s}^{(1)},\mathbf{s}^{(2)})\ge2^{n_0+1}$, which gives $2^{n_0+1}\le\mathbf{d}_{\mathcal{T}(\mathbf{b})}(\mathbf{s}^{(1)},\mathbf{s}^{(2)})<2^{n_0+1}$, contradiction.

\ref{item_minonept} If $\gamma$ is monotone as in \ref{item_geo_mono}, assume that $\gamma$ connects $[(\mathbf{u},\mathbf{t})]$, $[(\mathbf{v},\mathbf{s})]\in\mathcal{L}(\mathbf{g},\mathbf{b})$, then $\pi^\mathcal{T}(\gamma)=[\mathbf{t},\mathbf{s}]$ and $l(\gamma)=\mathbf{d}_{\mathcal{T}(\mathbf{b})}(\mathbf{t},\mathbf{s})$. Since for any $n\in\mathbb{Z}$, each $n$-cell in $(\mathcal{T}(\mathbf{b}),\mathbf{d}_{\mathcal{T}(\mathbf{b})})$ has diameter $2^n$ and any point in $(\mathcal{T}(\mathbf{b}),\mathbf{d}_{\mathcal{T}(\mathbf{b})})$ has distance at most $2^{n-1}$ from $\mathcal{W}_n$, we have for any $n\in\mathbb{Z}$ with $2^n\le \mathbf{d}_{\mathcal{T}(\mathbf{b})}(\mathbf{t},\mathbf{s})=l(\gamma)$, $\pi^\mathcal{T}(\gamma)\cap\mathcal{W}_n=[\mathbf{t},\mathbf{s}]\cap\mathcal{W}_n\ne\emptyset$. If $\gamma$ has one inversion as in \ref{item_geo_oneinv}, then $\gamma=\gamma^{(1)}\cup\gamma^{(2)}$, where $\gamma^{(1)}, \gamma^{(2)}$ are two monotone geodesics and $l(\gamma)=l(\gamma^{(1)})+l(\gamma^{(2)})$, then either $l(\gamma^{(1)})\ge l(\gamma)/2$ or $l(\gamma^{(2)})\ge l(\gamma)/2$, we may consider the first case, the second follows analogously. By the result for monotone geodesics, for any $n\in\mathbb{Z}$ with $2^n\le l(\gamma)/2\le l(\gamma^{(1)})$, we have $\pi^{\mathcal{T}}(\gamma)\cap\mathcal{W}_n\supseteq\pi^{\mathcal{T}}(\gamma^{(1)})\cap\mathcal{W}_n\ne\emptyset$. If $\gamma$ has two inversions as in \ref{item_geo_twoinv}, then $\gamma=\gamma^{(1)}\cup\gamma^{(2)}\cup\gamma^{(3)}$, where $\gamma^{(1)}, \gamma^{(2)}, \gamma^{(3)}$ are three monotone geodesics and $l(\gamma)=l(\gamma^{(1)})+l(\gamma^{(2)})+l(\gamma^{(3)})$, then there exists $i\in\{1,2,3\}$ such that $l(\gamma^{(i)})\ge l(\gamma)/3$. By the result for monotone geodesics, for any $n\in\mathbb{Z}$ with $2^n\le l(\gamma)/3\le l(\gamma^{(i)})$, we have $\pi^{\mathcal{T}}(\gamma)\cap\mathcal{W}_n\supseteq\pi^{\mathcal{T}}(\gamma^{(i)})\cap\mathcal{W}_n\ne\emptyset$. In summary, since $2^{n_0}\le l(\gamma)/4$, for any $n\le n_0$, we have $\pi^\mathcal{T}(\gamma)\cap\mathcal{W}_n\ne\emptyset$.
\end{proof}

To conclude this section, we prove that the quotient map $\mathcal{Q}:(\mathcal{P}(\mathbf{g},\mathbf{b}),\mathbf{d}_{\mathcal{P}(\mathbf{g},\mathbf{b})})\to(\mathcal{L}(\mathbf{g},\mathbf{b}),\mathbf{d}_{\mathcal{L}(\mathbf{g},\mathbf{b})})$ is David-Semmes regular.

\begin{definition}\label{def_DS}
Let $(X,\mathbf{d}_X)$, $(Y,\mathbf{d}_Y)$ be two metric spaces. Let $L,M,N$ be three positive integers. We say that a map $F:(X,\mathbf{d}_X)\to(Y,\mathbf{d}_Y)$ is David-Semmes regular with data $(L,M,N)$ if $F$ is $L$-Lipschitz and for any ball $B_Y(y,r)$ in $Y$, there exist $x_1,\ldots,x_M\in X$ such that
$$F^{-1}(B_Y(y,r))\subseteq\bigcup_{k=1}^MB_X(x_k,Nr).$$
We say that a map $F:(X,\mathbf{d}_X)\to(Y,\mathbf{d}_Y)$ is David-Semmes regular if it is David-Semmes regular with data $(L,M,N)$ for some positive integers $L,M,N$.
\end{definition}

\begin{lemma}\label{lem_DS}
The quotient map $\mathcal{Q}:(\mathcal{P}(\mathbf{g},\mathbf{b}),\mathbf{d}_{\mathcal{P}(\mathbf{g},\mathbf{b})})\to(\mathcal{L}(\mathbf{g},\mathbf{b}),\mathbf{d}_{\mathcal{L}(\mathbf{g},\mathbf{b})})$ is David-Semm-es regular. More precisely, $\mathcal{Q}$ is $3$-Lipschitz, and for any $[(\mathbf{u},\mathbf{t})]\in\mathcal{L}(\mathbf{g},\mathbf{b})$, $r>0$, let $n$ be the integer satisfying $2^{n-1}\le r<2^{n}$, by Lemma \ref{lem_onept} \ref{item_maxonept}, we have the following dichotomy.
\begin{itemize}
\item Either $B_{\mathcal{T}(\mathbf{b})}(\mathbf{t},r)\cap\cup_{k\ge n+2}\mathcal{W}_k=\emptyset$, then for any $(\mathbf{u},\mathbf{t})$ in $[(\mathbf{u},\mathbf{t})]$, we have
$$\mathcal{Q}^{-1}(B_{{\mathcal{L}(\mathbf{g},\mathbf{b})}}([(\mathbf{u},\mathbf{t})],r))\subseteq B_{\mathcal{U}(\mathbf{g})}(\mathbf{u},8r)\times B_{\mathcal{T}(\mathbf{b})}(\mathbf{t},r)\subseteq \mathcal{Q}^{-1}(B_{{\mathcal{L}(\mathbf{g},\mathbf{b})}}([(\mathbf{u},\mathbf{t})],32r)).$$
\item Or $B_{\mathcal{T}(\mathbf{b})}(\mathbf{t},r)\cap\cup_{k\ge n+2}\mathcal{W}_k=\{\mathbf{t}^{(1)}\}$, where $\mathbf{t}^{(1)}\in\mathcal{W}_m$ for some $m\ge n+2$, then for any $(\mathbf{u},\mathbf{t})$ in $[(\mathbf{u},\mathbf{t})]$, we have
\begin{align*}
&\mathcal{Q}^{-1}(B_{{\mathcal{L}(\mathbf{g},\mathbf{b})}}([(\mathbf{u},\mathbf{t})],r))\\
&\subseteq\left(\bigcup_{k=1}^{\mathbf{g}(m)}B_{\mathcal{U}(\mathbf{g})}(\mathbf{u}^{(k)},8r)\right)\times B_{\mathcal{T}(\mathbf{b})}(\mathbf{t},r)\\
&\subseteq\mathcal{Q}^{-1}\left(B_{\mathcal{L}(\mathbf{g},\mathbf{b})}([(\mathbf{u},\mathbf{t})],32r)\right),
\end{align*}
where $\mathbf{u}^{(1)}$, \ldots, $\mathbf{u}^{(\mathbf{g}(m))}\in\mathcal{U}(\mathbf{g})$ are \emph{all} the points satisfying $[(\mathbf{u}^{(1)},\mathbf{t}^{(1)})]=\ldots=[(\mathbf{u}^{(\mathbf{g}(m))},\mathbf{t}^{(1)})]=[(\mathbf{u},\mathbf{t}^{(1)})]$.
\end{itemize}
Hence $(\mathcal{L}(\mathbf{g},\mathbf{b}),\mathbf{d}_{\mathcal{L}(\mathbf{g},\mathbf{b})})$ is complete locally compact and separable.
\end{lemma}

\begin{proof}
Firstly, we show that $\mathcal{Q}$ is $3$-Lipschitz. For any $(\mathbf{u},\mathbf{t})$, $(\mathbf{v},\mathbf{s})\in\mathcal{P}(\mathbf{g},\mathbf{b})$, there exists a geodesic connecting $[(\mathbf{u},\mathbf{t})]$, $[(\mathbf{v},\mathbf{s})]\in\mathcal{L}(\mathbf{g},\mathbf{b})$ given by Proposition \ref{prop_Laa_geo}, along with three associated cases.

For \ref{item_geo_mono}, we have
$$\mathbf{d}_{\mathcal{L}(\mathbf{g},\mathbf{b})}([(\mathbf{u},\mathbf{t})], [(\mathbf{v},\mathbf{s})])=\mathbf{d}_{\mathcal{T}(\mathbf{b})}(\mathbf{t},\mathbf{s})\le\mathbf{d}_{\mathcal{P}(\mathbf{g},\mathbf{b})}((\mathbf{u},\mathbf{t}), (\mathbf{v},\mathbf{s})).$$

For \ref{item_geo_oneinv}, we have
$$\mathbf{d}_{\mathcal{L}(\mathbf{g},\mathbf{b})}([(\mathbf{u},\mathbf{t})], [(\mathbf{v},\mathbf{s})])=\mathbf{d}_{\mathcal{T}(\mathbf{b})}(\mathbf{t},\mathbf{s})+2\mathbf{d}_{\mathcal{T}(\mathbf{b})}(\mathbf{t}^{(1)},[\mathbf{t},\mathbf{s}]),$$
where $\mathbf{t}^{(1)}\in\mathcal{W}_{n^*}$ is given by the remark following Proposition \ref{prop_Laa_geo}. Since
$$\mathbf{d}_{\mathcal{T}(\mathbf{b})}(\mathbf{t}^{(1)},[\mathbf{t},\mathbf{s}])\le2^{n^*-1}\le 2^{\max\{n:\mathbf{u}(n)\ne\mathbf{v}(n)\}-1}=\frac{1}{2}\mathbf{d}_{\mathcal{U}(\mathbf{g})}(\mathbf{u},\mathbf{v}),$$
we have
$$\mathbf{d}_{\mathcal{L}(\mathbf{g},\mathbf{b})}([(\mathbf{u},\mathbf{t})], [(\mathbf{v},\mathbf{s})])\le\mathbf{d}_{\mathcal{T}(\mathbf{b})}(\mathbf{t},\mathbf{s})+\mathbf{d}_{\mathcal{U}(\mathbf{g})}(\mathbf{u},\mathbf{v})\le2\mathbf{d}_{\mathcal{P}(\mathbf{g},\mathbf{b})}((\mathbf{u},\mathbf{t}), (\mathbf{v},\mathbf{s})).$$

For \ref{item_geo_twoinv}, we have
$$\mathbf{d}_{\mathcal{L}(\mathbf{g},\mathbf{b})}([(\mathbf{u},\mathbf{t})], [(\mathbf{v},\mathbf{s})])=\mathbf{d}_{\mathcal{T}(\mathbf{b})}(\mathbf{t},\mathbf{s})+2\mathbf{d}_{\mathcal{T}(\mathbf{b})}(\mathbf{t}^{(1)},[\mathbf{t},\mathbf{s}])+2\mathbf{d}_{\mathcal{T}(\mathbf{b})}(\mathbf{t}^{(2)},[\mathbf{t},\mathbf{s}]),$$
where $\mathbf{t}^{(1)}\in\mathcal{W}_{n^*}$, $\mathbf{t}^{(2)}\in\mathcal{W}_{n^{**}}$ are given by the remark following Proposition \ref{prop_Laa_geo}. Since
\begin{align*}
\mathbf{d}_{\mathcal{T}(\mathbf{b})}(\mathbf{t}^{(1)},[\mathbf{t},\mathbf{s}])&\le2^{n^*-1}\le 2^{\max\{n:\mathbf{u}(n)\ne\mathbf{v}(n)\}-1}=\frac{1}{2}\mathbf{d}_{\mathcal{U}(\mathbf{g})}(\mathbf{u},\mathbf{v}),\\
\mathbf{d}_{\mathcal{T}(\mathbf{b})}(\mathbf{t}^{(2)},[\mathbf{t},\mathbf{s}])&\le2^{n^{**}-1}\le 2^{\max\{n:\mathbf{u}(n)\ne\mathbf{v}(n)\}-1}=\frac{1}{2}\mathbf{d}_{\mathcal{U}(\mathbf{g})}(\mathbf{u},\mathbf{v}),
\end{align*}
we have
$$\mathbf{d}_{\mathcal{L}(\mathbf{g},\mathbf{b})}([(\mathbf{u},\mathbf{t})], [(\mathbf{v},\mathbf{s})])\le\mathbf{d}_{\mathcal{T}(\mathbf{b})}(\mathbf{t},\mathbf{s})+\mathbf{d}_{\mathcal{U}(\mathbf{g})}(\mathbf{u},\mathbf{v})+\mathbf{d}_{\mathcal{U}(\mathbf{g})}(\mathbf{u},\mathbf{v})\le3\mathbf{d}_{\mathcal{P}(\mathbf{g},\mathbf{b})}((\mathbf{u},\mathbf{t}), (\mathbf{v},\mathbf{s})).$$

Secondly, for any $[(\mathbf{u},\mathbf{t})]\in\mathcal{L}(\mathbf{g},\mathbf{b})$, $r>0$, let $n$ be the integer satisfying $2^{n-1}\le r<2^n$. By Lemma \ref{lem_onept} \ref{item_maxonept}, $B_{\mathcal{T}(\mathbf{b})}(\mathbf{t},r)\cap\cup_{k\ge n+2}\mathcal{W}_k$ consists of at most one point.

If $B_{\mathcal{T}(\mathbf{b})}(\mathbf{t},r)\cap\cup_{k\ge n+2}\mathcal{W}_k=\emptyset$. Take $(\mathbf{u},\mathbf{t})$ in $[(\mathbf{u},\mathbf{t})]$. We claim that
$$\mathcal{Q}^{-1}(B_{\mathcal{L}(\mathbf{g},\mathbf{b})}([(\mathbf{u},\mathbf{t})],r))\subseteq\left\{\mathbf{v}\in\mathcal{U}(\mathbf{g}):\mathbf{v}|_{\mathbb{Z}\cap[n+2,+\infty)}=\mathbf{u}|_{\mathbb{Z}\cap[n+2,+\infty)}\right\}\times B_{\mathcal{T}(\mathbf{b})}(\mathbf{t},r).$$
Indeed, we prove that if $(\mathbf{v},\mathbf{s})\in\mathcal{P}(\mathbf{g},\mathbf{b})$ satisfies $\mathbf{d}_{\mathcal{T}(\mathbf{b})}(\mathbf{s},\mathbf{t})\ge r$ or $\mathbf{v}|_{\mathbb{Z}\cap[n+2,+\infty)}\ne\mathbf{u}|_{\mathbb{Z}\cap[n+2,+\infty)}$, then $\mathbf{d}_{\mathcal{L}(\mathbf{g},\mathbf{b})}([(\mathbf{v},\mathbf{s})],[(\mathbf{u},\mathbf{t})])\ge r$. For the first case, it is obvious that $\mathbf{d}_{\mathcal{L}(\mathbf{g},\mathbf{b})}([(\mathbf{v},\mathbf{s})],[(\mathbf{u},\mathbf{t})])\ge\mathbf{d}_{\mathcal{T}(\mathbf{b})}(\mathbf{s},\mathbf{t})\ge r$. For the second case, since $\mathbf{v}|_{\mathbb{Z}\cap[n+2,+\infty)}\ne\mathbf{u}|_{\mathbb{Z}\cap[n+2,+\infty)}$, any path connecting $[(\mathbf{v},\mathbf{s})],[(\mathbf{u},\mathbf{t})]$ needs to jump through at least one level-$k$ wormhole with $k\ge n+2$, however, $B_{\mathcal{T}(\mathbf{b})}(\mathbf{t},r)\cap\cup_{k\ge n+2}\mathcal{W}_k=\emptyset$, we have such path has length at least $r$, which implies that $\mathbf{d}_{\mathcal{L}(\mathbf{g},\mathbf{b})}([(\mathbf{v},\mathbf{s})],[(\mathbf{u},\mathbf{t})])\ge r$. Hence
\begin{align*}
&\mathcal{Q}^{-1}(B_{\mathcal{L}(\mathbf{g},\mathbf{b})}([(\mathbf{u},\mathbf{t})],r))\subseteq \overline{B}_{\mathcal{U}(\mathbf{g})}(\mathbf{u},2^{n+1})\times B_{\mathcal{T}(\mathbf{b})}(\mathbf{t},r)\\
&\subseteq B_{\mathcal{U}(\mathbf{g})}(\mathbf{u},8r)\times B_{\mathcal{T}(\mathbf{b})}(\mathbf{t},r)\subseteq B_{\mathcal{P}(\mathbf{g},\mathbf{b})}((\mathbf{u},\mathbf{t}),8r).
\end{align*}
Moreover, for any $\mathbf{v}\in B_{\mathcal{U}(\mathbf{g})}(\mathbf{u},8r)$, $\mathbf{s}\in B_{\mathcal{T}(\mathbf{b})}(\mathbf{t},r)$, we have
\begin{align*}
&\mathbf{d}_{\mathcal{L}(\mathbf{g},\mathbf{b})}([(\mathbf{v},\mathbf{s})],[(\mathbf{u},\mathbf{t})])\le\mathbf{d}_{\mathcal{L}(\mathbf{g},\mathbf{b})}([(\mathbf{v},\mathbf{s})],[(\mathbf{u},\mathbf{s})])+\mathbf{d}_{\mathcal{L}(\mathbf{g},\mathbf{b})}([(\mathbf{u},\mathbf{s})],[(\mathbf{u},\mathbf{t})])\\
&\le3\mathbf{d}_{\mathcal{P}(\mathbf{g},\mathbf{b})}((\mathbf{v},\mathbf{s}),(\mathbf{u},\mathbf{s}))+\mathbf{d}_{\mathcal{T}(\mathbf{b})}(\mathbf{s},\mathbf{t})=3\mathbf{d}_{\mathcal{U}(\mathbf{g})}(\mathbf{v},\mathbf{u})+\mathbf{d}_{\mathcal{T}(\mathbf{b})}(\mathbf{s},\mathbf{t})<3\cdot 8r+r=25r,
\end{align*}
which gives
$$B_{\mathcal{U}(\mathbf{g})}(\mathbf{u},8r)\times B_{\mathcal{T}(\mathbf{b})}(\mathbf{t},r)\subseteq \mathcal{Q}^{-1}(B_{{\mathcal{L}(\mathbf{g},\mathbf{b})}}([(\mathbf{u},\mathbf{t})],25r)).$$

If $B_{\mathcal{T}(\mathbf{b})}(\mathbf{t},r)\cap\cup_{k\ge n+2}\mathcal{W}_k=\{\mathbf{t}^{(1)}\}$ is a one-point set, where $\mathbf{t}^{(1)}\in\mathcal{W}_m$ for some $m\ge n+2$. Take $(\mathbf{u},\mathbf{t})$ in $[(\mathbf{u},\mathbf{t})]$. For any $k=1,\ldots,\mathbf{g}(m)$, let $\mathbf{u}^{(k)}\in\mathcal{U}(\mathbf{g})$ be given by $\mathbf{u}^{(k)}|_{\mathbb{Z}\backslash\{m\}}=\mathbf{u}|_{\mathbb{Z}\backslash\{m\}}$ and $\mathbf{u}^{(k)}(m)=k-1$, then $[(\mathbf{u}^{(k)},\mathbf{t}^{(1)})]=[(\mathbf{u},\mathbf{t}^{(1)})]$. We claim that
\begin{align*}
&\mathcal{Q}^{-1}(B_{\mathcal{L}(\mathbf{g},\mathbf{b})}([(\mathbf{u},\mathbf{t})],r))\\
&\subseteq\left(\bigcup_{k=1}^{\mathbf{g}(m)}\left\{\mathbf{v}\in\mathcal{U}(\mathbf{g}):\mathbf{v}|_{\mathbb{Z}\cap[n+2,+\infty)}=\mathbf{u}^{(k)}|_{\mathbb{Z}\cap[n+2,+\infty)}\right\}\right)\times B_{\mathcal{T}(\mathbf{b})}(\mathbf{t},r).
\end{align*}
Indeed, we prove that if $(\mathbf{v},\mathbf{s})\in\mathcal{P}(\mathbf{g},\mathbf{b})$ satisfies $\mathbf{d}_{\mathcal{T}(\mathbf{b})}(\mathbf{s},\mathbf{t})\ge r$ or $\mathbf{v}|_{\mathbb{Z}\cap[n+2,+\infty)}\ne\mathbf{u}^{(k)}|_{\mathbb{Z}\cap[n+2,+\infty)}$ for any $k=1,\ldots,\mathbf{g}(m)$, then $\mathbf{d}_{\mathcal{L}(\mathbf{g},\mathbf{b})}([(\mathbf{v},\mathbf{s})],[(\mathbf{u},\mathbf{t})])\ge r$. For the first case, it is obvious that $\mathbf{d}_{\mathcal{L}(\mathbf{g},\mathbf{b})}([(\mathbf{v},\mathbf{s})],[(\mathbf{u},\mathbf{t})])\ge\mathbf{d}_{\mathcal{T}(\mathbf{b})}(\mathbf{s},\mathbf{t})\ge r$. For the second case, we have $\mathbf{v}|_{(\mathbb{Z}\cap[n+2,+\infty))\backslash\{m\}}\ne\mathbf{u}|_{(\mathbb{Z}\cap[n+2,+\infty))\backslash\{m\}}$, otherwise, $\mathbf{v}|_{\mathbb{Z}\cap[n+2,+\infty)}=\mathbf{u}^{(k)}|_{\mathbb{Z}\cap[n+2,+\infty)}$ with $k=\mathbf{v}(m)+1$. Then there exists $p\ge n+2$ with $p\ne m$ such that $\mathbf{v}(p)\ne\mathbf{u}(p)$, any path connecting $[(\mathbf{v},\mathbf{s})],[(\mathbf{u},\mathbf{t})]$ needs to jump through at least one level-$p$ wormhole, however, since $B_{\mathcal{T}(\mathbf{b})}(\mathbf{t},r)\cap\cup_{k\ge n+2}\mathcal{W}_k=\{\mathbf{t}^{(1)}\}\subseteq\mathcal{W}_m$, we have $B_{\mathcal{T}(\mathbf{b})}(\mathbf{t},r)\cap\mathcal{W}_p=\emptyset$, such path has length at least $r$, which implies that $\mathbf{d}_{\mathcal{L}(\mathbf{g},\mathbf{b})}([(\mathbf{v},\mathbf{s})],[(\mathbf{u},\mathbf{t})])\ge r$. Hence
\begin{align*}
&\mathcal{Q}^{-1}(B_{\mathcal{L}(\mathbf{g},\mathbf{b})}([(\mathbf{u},\mathbf{t})],r))\subseteq \left(\bigcup_{k=1}^{\mathbf{g}(m)} \overline{B}_{\mathcal{U}(\mathbf{g})}(\mathbf{u}^{(k)},2^{n+1})\right)\times B_{\mathcal{T}(\mathbf{b})}(\mathbf{t},r)\\
&\subseteq\left(\bigcup_{k=1}^{\mathbf{g}(m)}B_{\mathcal{U}(\mathbf{g})}(\mathbf{u}^{(k)},8r)\right)\times B_{\mathcal{T}(\mathbf{b})}(\mathbf{t},r)\subseteq\bigcup_{k=1}^{\mathbf{g}(m)}B_{\mathcal{P}(\mathbf{g},\mathbf{b})}((\mathbf{u}^{(k)},\mathbf{t}),8r).
\end{align*}
Moreover, for any $\mathbf{v}\in B_{\mathcal{U}(\mathbf{g})}(\mathbf{u}^{(k)},8r)$, $\mathbf{s}\in B_{\mathcal{T}(\mathbf{b})}(\mathbf{t},r)$, since $[(\mathbf{u}^{(k)},\mathbf{t}^{(1)})]=[(\mathbf{u},\mathbf{t}^{(1)})]$ and $\mathbf{t}^{(1)}\in B_{\mathcal{T}(\mathbf{b})}(\mathbf{t},r)$, we have
\begin{align*}
&\mathbf{d}_{\mathcal{L}(\mathbf{g},\mathbf{b})}([(\mathbf{v},\mathbf{s})],[(\mathbf{u},\mathbf{t})])\\
&\le\mathbf{d}_{\mathcal{L}(\mathbf{g},\mathbf{b})}([(\mathbf{v},\mathbf{s})],[(\mathbf{u}^{(k)},\mathbf{s})])\\
&\hspace{10pt}+\mathbf{d}_{\mathcal{L}(\mathbf{g},\mathbf{b})}([(\mathbf{u}^{(k)},\mathbf{s})],[(\mathbf{u}^{(k)},\mathbf{t}^{(1)})])+\mathbf{d}_{\mathcal{L}(\mathbf{g},\mathbf{b})}([(\mathbf{u}^{(k)},\mathbf{t}^{(1)})],[(\mathbf{u},\mathbf{t})])\\
&\le3\mathbf{d}_{\mathcal{P}(\mathbf{g},\mathbf{b})}((\mathbf{v},\mathbf{s}),(\mathbf{u}^{(k)},\mathbf{s}))+\mathbf{d}_{\mathcal{T}(\mathbf{b})}(\mathbf{s},\mathbf{t}^{(1)})+\mathbf{d}_{\mathcal{L}(\mathbf{g},\mathbf{b})}([(\mathbf{u},\mathbf{t}^{(1)})],[(\mathbf{u},\mathbf{t})])\\
&=3\mathbf{d}_{\mathcal{U}(\mathbf{g})}(\mathbf{v},\mathbf{u}^{(k)})+\mathbf{d}_{\mathcal{T}(\mathbf{b})}(\mathbf{s},\mathbf{t}^{(1)})+\mathbf{d}_{\mathcal{T}(\mathbf{b})}(\mathbf{t}^{(1)},\mathbf{t})<3\cdot8r+2r+r=27r,
\end{align*}
which gives
$$\left(\bigcup_{k=1}^{\mathbf{g}(m)}B_{\mathcal{U}(\mathbf{g})}(\mathbf{u}^{(k)},8r)\right)\times B_{\mathcal{T}(\mathbf{b})}(\mathbf{t},r)\subseteq\mathcal{Q}^{-1}\left(B_{\mathcal{L}(\mathbf{g},\mathbf{b})}([(\mathbf{u},\mathbf{t})],27r)\right).$$

Therefore, in both of the aforementioned cases, $\mathcal{Q}^{-1}(B_{\mathcal{L}(\mathbf{g},\mathbf{b})}([(\mathbf{u},\mathbf{t})],r))$ can always be contained in the union of $\sup_\mathbb{Z}\mathbf{g}$ balls with radius $8r$.

In summary, $\mathcal{Q}:(\mathcal{P}(\mathbf{g},\mathbf{b}),\mathbf{d}_{\mathcal{P}(\mathbf{g},\mathbf{b})})\to(\mathcal{L}(\mathbf{g},\mathbf{b}),\mathbf{d}_{\mathcal{L}(\mathbf{g},\mathbf{b})})$ is David-Semmes regular.

Finally, for any Cauchy sequence $\{[\mathbf{q}_n]\}_{n\ge1}$ in $\mathcal{L}(\mathbf{g},\mathbf{b})$, we have $\{[\mathbf{q}_n]\}_{n\ge1}$ is bounded in $\mathcal{L}(\mathbf{g},\mathbf{b})$, hence $\mathcal{Q}^{-1}(\{[\mathbf{q}_n]\}_{n\ge1})$ is bounded in $\mathcal{P}(\mathbf{g},\mathbf{b})$. Since $\mathcal{P}(\mathbf{g},\mathbf{b})$ is complete and locally compact, $\mathcal{Q}^{-1}(\{[\mathbf{q}_n]\}_{n\ge1})$ has a limit point $\mathbf{q}\in\mathcal{P}(\mathbf{g},\mathbf{b})$. By the Lipschitz continuity of $\mathcal{Q}$, $[\mathbf{q}]=\mathcal{Q}(\mathbf{q})\in\mathcal{L}(\mathbf{g},\mathbf{b})$ is the limit of $\{[\mathbf{q}_n]\}_{n\ge1}$ in $\mathcal{L}(\mathbf{g},\mathbf{b})$, hence $\mathcal{L}(\mathbf{g},\mathbf{b})$ is complete.

For any bounded closed set $K$ in $\mathcal{L}(\mathbf{g},\mathbf{b})$, we have $\mathcal{Q}^{-1}(K)$ is bounded closed in $\mathcal{P}(\mathbf{g},\mathbf{b})$. Since $\mathcal{P}(\mathbf{g},\mathbf{b})$ is locally compact, we have $\mathcal{Q}^{-1}(K)$ is compact in $\mathcal{P}(\mathbf{g},\mathbf{b})$. By the continuity of $\mathcal{Q}$, we have $K=\mathcal{Q}(\mathcal{Q}^{-1}(K))$ is compact in $\mathcal{L}(\mathbf{g},\mathbf{b})$, hence $\mathcal{L}(\mathbf{g},\mathbf{b})$ is locally compact.

Since $\mathcal{P}(\mathbf{g},\mathbf{b})$ is separable, let $C$ be a countable dense subset of $\mathcal{P}(\mathbf{g},\mathbf{b})$, then $\mathcal{Q}(C)$ is a countable dense subset of $\mathcal{L}(\mathbf{g},\mathbf{b})$. Hence $\mathcal{L}(\mathbf{g},\mathbf{b})$ is separable.
\end{proof}

Since $\mathcal{Q}$ is David-Semmes regular, we push forward measures on $\mathcal{P}(\mathbf{g},\mathbf{b})$ to measures on $\mathcal{L}(\mathbf{g},\mathbf{b})$ as follows.
\begin{align*}
\mathbf{m}_{\mathcal{L}(\mathbf{g},\mathbf{b})}&=\mathcal{Q}_*(\mathbf{m}_{\mathcal{U}(\mathbf{g})}\times\mathbf{m}_{\mathcal{T}(\mathbf{b})}),\\
\lambda_{\mathcal{L}(\mathbf{g},\mathbf{b})}&=\mathcal{Q}_*(\mathbf{m}_{\mathcal{U}(\mathbf{g})}\times\bf\lambda_{\mathcal{T}(\mathbf{b})}).
\end{align*}
There exists some positive constant $C$ depending only on $\sup_\mathbb{Z}\mathbf{g}$ and $\sup_\mathbb{Z}\mathbf{b}$ such that
\begin{equation}\label{eq_Laa_vol}
\frac{1}{C}V_\mathbf{g}(r)V_\mathbf{b}(r)\le\mathbf{m}_{\mathcal{L}(\mathbf{g},\mathbf{b})}(B_{\mathcal{L}(\mathbf{g},\mathbf{b})}([\mathbf{p}],r))\le{C}V_\mathbf{g}(r)V_\mathbf{b}(r)\text{ for any }[\mathbf{p}]\in\mathcal{L}(\mathbf{g},\mathbf{b}), r>0.
\end{equation}

\section{\texorpdfstring{$p$}{p}-energy on Laakso-type spaces}\label{sec_Laa_energy}

In this section, we introduce a $p$-energy with a $p$-energy measure on Laakso-type spaces, and prove the Poincar\'e inequality, the capacity upper bound and the cutoff Sobolev inequality.

Firstly, let us introduce a $p$-energy on $\mathcal{L}(\mathbf{g},\mathbf{b})$. For any function $f$ on $\mathcal{L}(\mathbf{g},\mathbf{b})$, let $\widehat{f}=f\circ\mathcal{Q}$ be a function on $\mathcal{P}(\mathbf{g},\mathbf{b})$. If $f\in C(\mathcal{L}(\mathbf{g},\mathbf{b}))$, since $\mathcal{Q}$ is continuous, $\widehat{f}\in C(\mathcal{P}(\mathbf{g},\mathbf{b}))$. Moreover, if $f\in C_c(\mathcal{L}(\mathbf{g},\mathbf{b}))$, since $\mathrm{supp}(f)$ is compact in $\mathcal{L}(\mathbf{g},\mathbf{b})$ and $\mathcal{Q}$ is David-Semmes regular, $\mathrm{supp}(\widehat{f})=\mathcal{Q}^{-1}(\mathrm{supp}(f))$ is compact in $\mathcal{P}(\mathbf{g},\mathbf{b})$, $\widehat{f}\in C_c(\mathcal{P}(\mathbf{g},\mathbf{b}))$.

For any $n\in\mathbb{Z}$, let
\begin{align*}
\mathcal{C}^\mathcal{L}_n&=\left\{f\in C_c(\mathcal{L}(\mathbf{g},\mathbf{b})):\widehat{f}=f\circ\mathcal{Q}\in C_c(\mathcal{P}(\mathbf{g},\mathbf{b}))\right.\\
&\left.\text{for any }\mathbf{u}\in\mathcal{U}(\mathbf{g}), \widehat{f}(\mathbf{u},\cdot)\in\mathcal{F}^\mathcal{T},\right.\\
&\left.\text{for any }\mathbf{t}\in\mathcal{T}(\mathbf{b}), \widehat{f}(\cdot,\mathbf{t})\text{ is constant on any closed ball with radius }2^n\right\}.
\end{align*}
Let $\mathcal{C}^\mathcal{L}=\cup_{n\in\mathbb{Z}}\mathcal{C}^\mathcal{L}_n$. For any $f\in\mathcal{C}^\mathcal{L}$ and any $[(\mathbf{u},\mathbf{t})]\in\mathcal{L}(\mathbf{g},\mathbf{b})$, let
$$|\nabla_\mathcal{L} f([(\mathbf{u},\mathbf{t})])|=|\nabla_\mathcal{T}\widehat{f}(\mathbf{u},\cdot)(\mathbf{t})|.$$
Let
$$\mathcal{E}^\mathcal{L}(f)=\int_{\mathcal{L}(\mathbf{g},\mathbf{b})}|\nabla_\mathcal{L} f|^p\mathrm{d}\lambda_{\mathcal{L}(\mathbf{g},\mathbf{b})}=\int_{\mathcal{U}(\mathbf{g})}\left(\int_{\mathcal{T}(\mathbf{b})}|\nabla_\mathcal{T}\widehat{f}(\mathbf{u},\cdot)(\mathbf{t})|^p\lambda_{\mathcal{T}(\mathbf{b})}(\mathrm{d}\mathbf{t})\right)\mathbf{m}_{\mathcal{U}(\mathbf{g})}(\mathrm{d}\mathbf{u}),$$
and
$$\mathcal{F}^\mathcal{L}=\text{the }(\mathcal{E}^\mathcal{L}(\cdot)^{1/p}+\lVert \cdot\rVert_{L^p(\mathcal{L}(\mathbf{g},\mathbf{b});\mathbf{m}_{\mathcal{L}(\mathbf{g},\mathbf{b})})})\text{-closure of }\mathcal{C}^\mathcal{L}.$$
Using approximating Cauchy sequences in the corresponding $L^p$-spaces, the function $|\nabla_\mathcal{L} f|\in L^p(\mathcal{L}(\mathbf{g},\mathbf{b});\lambda_{\mathcal{L}(\mathbf{g},\mathbf{b})})$ is well-defined for any $f\in\mathcal{F}^\mathcal{L}$.

We have the following result.

\begin{lemma}\label{lem_Laa_energy}
\begin{enumerate}[label=(\arabic*)]
\item $\mathcal{C}^\mathcal{L}$ is uniformly dense in $C_c(\mathcal{L}(\mathbf{g},\mathbf{b}))$.
\item $(\mathcal{E}^\mathcal{L},\mathcal{F}^\mathcal{L})$ is a $p$-energy on $(\mathcal{L}(\mathbf{g},\mathbf{b}),\mathbf{d}_{\mathcal{L}(\mathbf{g},\mathbf{b})},\mathbf{m}_{\mathcal{L}(\mathbf{g},\mathbf{b})})$ with a $p$-energy measure $\Gamma^\mathcal{L}$ given by 
$$\Gamma^\mathcal{L}(f)(A)=\int_A|\nabla_\mathcal{L} f|^p\mathrm{d}\lambda_{\mathcal{L}(\mathbf{g},\mathbf{b})}$$
for any $f\in\mathcal{F}^\mathcal{L}$, $A\in\mathcal{B}(\mathcal{L}(\mathbf{g},\mathbf{b}))$.
\end{enumerate}
\end{lemma}

\begin{proof}
We only give the proof of (1). The proof of (2) is easy and also similar to the proof of  $(\int_{\mathbb{R}}|f'(x)|^p\mathrm{d} x,W^{1,p}(\mathbb{R}))$ is a $p$-energy on $\mathbb{R}$ with a $p$-energy measure given by $(f,A)\mapsto\int_A|f'(x)|^p\mathrm{d} x$, see also \cite[Sections 4 and 5]{Ste08} or \cite[Subsections 4.4 and 4.5]{Ste10} for the construction of a Dirichlet form on the original Laakso space in \cite{Laa00}.

Indeed, since $\mathcal{C}^\mathcal{L}$ is a sub-algebra of $C_c(\mathcal{L}(\mathbf{g},\mathbf{b}))$, by the Stone-Weierstra{\ss} theorem, we only need to show that $\mathcal{C}^\mathcal{L}$ separates points and vanishes nowhere. For any distinct $x,y\in\mathcal{L}(\mathbf{g},\mathbf{b})$, let $D=\mathbf{d}_{\mathcal{L}(\mathbf{g},\mathbf{b})}(x,y)>0$, by the following Proposition \ref{prop_Laa_CS}, there exists $\phi_{\mathcal{L}}\in\mathcal{C}^{\mathcal{L}}$ such that $\phi_\mathcal{L}=1$ in $B_{\mathcal{L}(\mathbf{g},\mathbf{b})}(x,D/256)$ and $\phi_\mathcal{L}=0$ on $\mathcal{L}(\mathbf{g},\mathbf{b})\backslash B_{\mathcal{L}(\mathbf{g},\mathbf{b})}(x,D)$, in particular, we have $\phi_\mathcal{L}(x)=1\ne0=\phi_\mathcal{L}(y)$.
\end{proof}

Secondly, we prove the Poincar\'e inequality. We now present the main trick, known as the technique of the pencil of curves, as follows.

\begin{proposition}\label{prop_pencil}
There exists $C>0$ such that for any distinct $x, y\in\mathcal{L}(\mathbf{g},\mathbf{b})$, there exist a family $\Gamma_{x,y}$ of geodesics connecting $x,y$ and a probability measure $\mathbb{P}^{x,y}$ on $\Gamma_{x,y}$ such that for any non-negative measurable function $h$ on $\mathcal{L}(\mathbf{g},\mathbf{b})$, we have
\begin{align*}
&\int_{\Gamma_{x,y}}\left(\int_0^{\mathbf{d}_{\mathcal{L}(\mathbf{g},\mathbf{b})}(x,y)}h(\gamma(t))\mathrm{d} t\right)\mathrm{d}\mathbb{P}^{x,y}\\
&\le C\int_{B_{\mathcal{L}(\mathbf{g},\mathbf{b})}(x,\mathbf{d}_{\mathcal{L}(\mathbf{g},\mathbf{b})}(x,y))}\frac{h(z)}{V_\mathbf{g}(\mathbf{d}_{\mathcal{L}(\mathbf{g},\mathbf{b})}(x,z))\wedge V_\mathbf{g}(\mathbf{d}_{\mathcal{L}(\mathbf{g},\mathbf{b})}(y,z))}\lambda_{\mathcal{L}(\mathbf{g},\mathbf{b})}(\mathrm{d} z).
\end{align*}
\end{proposition}

\begin{remark}
We do \emph{NOT} require the family $\Gamma_{x,y}$ to be the family of \emph{ALL} geodesics connecting $x,y$. Roughly speaking, we will add some ``unnecessary" jumps to a geodesic $\gamma_0$, which was given by Proposition \ref{prop_Laa_geo}, to obtain a family of geodesics $\Gamma_{x,y}$, which would be rich enough for our purpose.
\end{remark}

\begin{proof}
Let $\gamma_0$ be a geodesic connecting $x,y$ given by Proposition \ref{prop_Laa_geo}. Let $D=\mathbf{d}_{\mathcal{L}(\mathbf{g},\mathbf{b})}(x,y)$ and $k_0$ the integer satisfying $2^{k_0+3}\le D<2^{k_0+4}$. Let $z=\gamma_0(D/2)$, then $\gamma_0|_{[0,D/2]}$ is a geodesic connecting $x, z$ given by Proposition \ref{prop_Laa_geo}, and $\gamma_0|_{[D/2,D]}$ is a geodesic connecting $z, y$ given by Proposition \ref{prop_Laa_geo}. Moreover, $l(\gamma_0|_{[0,D/2]})=l(\gamma_0|_{[D/2,D]})=D/2$. Since $2^{k_0+2}\le D/2$, by Lemma \ref{lem_onept} \ref{item_minonept}, for any $k\le k_0$, $\pi^\mathcal{T}(\gamma_0|_{[0,D/2]})\cap\mathcal{W}_k\ne\emptyset$, $\pi^\mathcal{T}(\gamma_0|_{[D/2,D]})\cap\mathcal{W}_k\ne\emptyset$.

For any $k\le k_0$, there exist only finitely many $t\in[0,D]$ such that $\pi^\mathcal{T}(\gamma_0(t))\in\mathcal{W}_k$, let
\begin{align*}
t_{k,1}&=\min\{t\in[0,D]:\pi^\mathcal{T}(\gamma_0(t))\in\mathcal{W}_k\},\\
t_{k,2}&=\max\{t\in[0,D]:\pi^\mathcal{T}(\gamma_0(t))\in\mathcal{W}_k\},
\end{align*}
then $t_{k,1}\in[0,D/2]$, $t_{k,2}\in[D/2,D]$. Moreover, $t_{k,1}\in[0,2^{k+2}]$, $t_{k,2}\in[D-2^{k+2},D]$. Indeed, for any $t\in[0,D]$, $\gamma_0|_{[0,t]}$ is a geodesic given by Proposition \ref{prop_Laa_geo} with length $t$, by Lemma \ref{lem_onept} \ref{item_minonept}, for any $t\ge2^{k+2}$, $\pi^{\mathcal{T}}(\gamma_0|_{[0,t]})\cap\mathcal{W}_k\ne\emptyset$, hence $t_{k,1}\in[0,2^{k+2}]$. Similarly, $t_{k,2}\in[D-2^{k+2},D]$. Since $\{\mathcal{W}_k\}_{k\in\mathbb{Z}}$ are disjoint, we have $\{\{t_{k,1},t_{k,2}\}\}_{k\le k_0}$ are also disjoint.

For any $t\in(0,D)$, let
$$k(t)=\min\{\lfloor \log_2t\rfloor,\lfloor{\log_2(D-t)} \rfloor\}-3,$$
then for any $k\le k(t)$, $t_{k,1}\le 2^{k+2}\le2^{k(t)+2}\le t/2<t$, $t_{k,2}\ge D-2^{k+2}\ge D-2^{k(t)+2}\ge D-(D-t)/2=(D+t)/2>t$, hence $t\in(t_{k,1},t_{k,2})$.

Let $(\Omega,\mathcal{F},\mathbb{P})$ be a probability space on which there exist independent random variables $\{\xi_{k}:k\le k_0\}$ satisfying that $\xi_{k}$ has the uniform distribution $\mathrm{Unif}(\{0,1,\ldots,\mathbf{g}(k)-1\})$. Write $\gamma_0(t)=[(\mathbf{v}_0(t),\mathbf{s}_0(t))]$ for $t\in[0,D]$. For any $\omega\in\Omega$, let $\gamma(\omega)=[(\mathbf{v}(\omega),\mathbf{s}(\omega))]:[0,D]\to\mathcal{L}(\mathbf{g},\mathbf{b})$, where $\mathbf{v}(\omega):[0,D]\to\mathcal{U}(\mathbf{g})$, $\mathbf{s}(\omega):[0,D]\to\mathcal{T}(\mathbf{b})$, be given as follows. Let $\mathbf{s}(\omega)=\mathbf{s}_0$. Recall that
$$\mathcal{U}(\mathbf{g})=\left\{\mathbf{u}:\mathbb{Z}\to\mathbb{Z}|\mathbf{u}(k)\in\{0,1,\ldots,\mathbf{g}(k)-1\}\text{ for any }k\in\mathbb{Z}\text{ and }\lim_{k\to+\infty}\mathbf{u}(k)=0\right\}.$$
We define the value of $\mathbf{v}(\omega)(t)(n)$ for $t\in[0,D]$ and $n\in\mathbb{Z}$ as follows. For any $t\in[0,D]$, let $\mathbf{v}(\omega)(t)|_{\mathbb{Z}\cap[k_0+1,+\infty)}=\mathbf{v}_0(t)|_{\mathbb{Z}\cap[k_0+1,+\infty)}$. For any $k\le k_0$, let
$$
\mathbf{v}(\omega)(t)(k)=
\begin{cases}
\mathbf{v}_0(t)(k)&\text{if }t\in[0,t_{k,1}]\cup(t_{k,2},D],\\
\xi_{k}(\omega)&\text{if }t\in(t_{k,1},t_{k,2}].
\end{cases}
$$

We claim that $\gamma(\omega)$ is a geodesic connecting $x,y$. It is obvious that $\gamma(\omega)(0)=x$, $\gamma(\omega)(D)=y$. We show that $\gamma(\omega):[0,D]\to(\mathcal{L}(\mathbf{g},\mathbf{b}),\mathbf{d}_{\mathcal{L}(\mathbf{g},\mathbf{b})})$ is continuous. For any $t\in(0,D)$, there exists an open interval $U\subset(0,D)$ containing $t$ such that $k(\cdot)\ge k(t)-1$ in $U$, then for any $k\le k(t)-1$, we have $U\subseteq(t_{k,1},t_{k,2})$, which implies $\mathbf{v}(\omega)(\cdot)|_{\mathbb{Z}\cap(-\infty,k(t)-1]}$ is constant in $U$. By definition, $U$ can be written as a \emph{finite} union of intervals such that $\mathbf{v}(\omega)(\cdot)|_{\mathbb{Z}\cap[k(t),+\infty)}$ is constant on each interval,\footnote{Notice that high level wormholes can only be used \emph{finitely} many times for a geodesic.} hence $\mathbf{v}(\omega)(\cdot)$ is constant on each interval, where the discontinuity between adjacent intervals comes from the jumps through wormholes, however, by the definition of $R_\mathcal{L}$, the function of the equivalence classes $\gamma(\omega)(\cdot)=[(\mathbf{v}(\omega)(\cdot),\mathbf{s}(\omega)(\cdot))]=[(\mathbf{v}(\omega)(\cdot),\mathbf{s}_0(\cdot))]$ is indeed continuous in $U$, hence $\gamma(\omega)$ is continuous at $t$ for any $t\in(0,D)$. At $0$, we have
\begin{align*}
&\mathbf{d}_{\mathcal{L}(\mathbf{g},\mathbf{b})}(\gamma(\omega)(t),\gamma(\omega)(0))=\mathbf{d}_{\mathcal{L}(\mathbf{g},\mathbf{b})}([(\mathbf{v}(\omega)(t),\mathbf{s}_0(t))],[(\mathbf{v}(\omega)(0),\mathbf{s}_0(0))])\\
&\le 3\mathbf{d}_{\mathcal{P}(\mathbf{g},\mathbf{b})}((\mathbf{v}(\omega)(t),\mathbf{s}_0(t)),(\mathbf{v}_0(0),\mathbf{s}_0(0)))\\
&=3\max\left\{\mathbf{d}_{\mathcal{U}(\mathbf{g})}(\mathbf{v}(\omega)(t),\mathbf{v}_0(0)),\mathbf{d}_{\mathcal{T}(\mathbf{b})}(\mathbf{s}_0(t),\mathbf{s}_0(0))\right\},
\end{align*}
where the above inequality follows from Lemma \ref{lem_DS}. Since $\lim_{t\downarrow0}\mathbf{d}_{\mathcal{T}(\mathbf{b})}(\mathbf{s}_0(t),\mathbf{s}_0(0))=0$, we only need to consider the term $\mathbf{v}(\omega)$. For any $\varepsilon>0$, let $k_1\in\mathbb{Z}$ satisfy $2^{k_1}<\varepsilon$. By definition, $[0,D]$ can be written as a \emph{finite} union of intervals such that $\mathbf{v}(\omega)(\cdot)|_{\mathbb{Z}\cap[k_1+1,+\infty)}$ is constant on each interval, take $(0,\delta)$ on which $\mathbf{v}(\omega)(\cdot)|_{\mathbb{Z}\cap[k_1+1,+\infty)}$ is constant, by choosing another element from the equivalence class $[(\mathbf{v}_0(0),\mathbf{s}_0(0))]$, we may assume that $\mathbf{v}_0(0)|_{\mathbb{Z}\cap[k_1+1,+\infty)}$ is this constant, then $\mathbf{d}_{\mathcal{U}(\mathbf{g})}(\mathbf{v}(\omega)(t),\mathbf{v}_0(0))\le2^{k_1}<\varepsilon$ for any $t\in(0,\delta)$, hence $\gamma(\omega)$ is continuous at $0$. Similarly, $\gamma(\omega)$ is continuous at $D$. Therefore, $\gamma(\omega):[0,D]\to(\mathcal{L}(\mathbf{g},\mathbf{b}),\mathbf{d}_{\mathcal{L}(\mathbf{g},\mathbf{b})})$ is continuous, $\gamma(\omega)$ is a path connecting $x, y$. By construction, $l(\gamma(\omega))=l(\gamma_0)=D$, hence $\gamma(\omega)$ is a geodesic connecting $x, y$.

Let $\Gamma_{x,y}=\{\gamma(\omega):\omega\in\Omega\}$ and $\mathbb{P}^{x,y}=\mathbb{P}$. We show $\Gamma_{x,y}$ and $\mathbb{P}^{x,y}$ give the desired result. Except for the end point $y$, all geodesics in $\Gamma_{x,y}$ are contained within the ball $B=B_{\mathcal{L}(\mathbf{g},\mathbf{b})}(x,D)$. For any non-negative measurable function $h$ on $\mathcal{L}(\mathbf{g},\mathbf{b})$, we have
\begin{align*}
&\int_{\Gamma_{x,y}}\left(\int_0^Dh(\gamma(t))\mathrm{d} t\right)\mathrm{d}\mathbb{P}^{x,y}=\int_{\Omega}\left(\int_0^Dh(\gamma(\omega)(t))1_{B}(\gamma(\omega)(t))\mathrm{d} t\right)\mathbb{P}(\mathrm{d}\omega)\\
&=\int_0^D\left(\int_\Omega (h1_B)\left([(\mathbf{v}(\omega)(t),\mathbf{s}(\omega)(t))]\right)\mathbb{P}(\mathrm{d}\omega)\right)\mathrm{d} t\\
&=\int_0^D\left(\int_\Omega (h1_B)\left([(\mathbf{v}(\omega)(t),\mathbf{s}_0(t))]\right)\mathbb{P}(\mathrm{d}\omega)\right)\mathrm{d} t.
\end{align*}

For fixed $t\in(0,D)$, let us consider the distribution $\mathbb{P}_{\mathbf{v}(\cdot)(t)}$ of $\mathbf{v}(\cdot)(t)$ under $\mathbb{P}$. Recall that for any 
$k\le k(t)$, we have $t\in(t_{k,1},t_{k,2})$, then $\mathbf{v}(\cdot)(t)(k)=\xi_k$. For any $\mathbf{u}\in\mathcal{U}(\mathbf{g})$, $n\in\mathbb{Z}$, we have
$$\mathbb{P}\left[\omega:\mathbf{v}(\omega)(t)\in{B}_{\mathcal{U}(\mathbf{g})}(\mathbf{u},2^n)\right]=\mathbb{P}\left[\omega:\mathbf{v}(\omega)(t)|_{\mathbb{Z}\cap[n,+\infty)}=\mathbf{u}|_{\mathbb{Z}\cap[n,+\infty)}\right].$$
If $n<k(t)$, then
\begin{align*}
&\mathbb{P}\left[\omega:\mathbf{v}(\omega)(t)|_{\mathbb{Z}\cap[n,+\infty)}=\mathbf{u}|_{\mathbb{Z}\cap[n,+\infty)}\right]\\
&\le\mathbb{P}\left[\omega:\xi_k(\omega)=\mathbf{u}(k)\text{ for any }k=n,\ldots,k(t)\right]=\frac{1}{\prod_{k=n}^{k(t)}\mathbf{g}(k)}.
\end{align*}
If $n\ge1$, then
$$\frac{1}{\prod_{k=n}^{k(t)}\mathbf{g}(k)}=\frac{\prod_{k=1}^{n-1}\mathbf{g}(k)}{\prod_{k=1}^{k(t)}\mathbf{g}(k)}\asymp \frac{V_\mathbf{g}(2^n)}{V_\mathbf{g}(2^{k(t)})}.$$
If $n\le0$, then
\begin{align*}
&\frac{1}{\prod_{k=n}^{k(t)}\mathbf{g}(k)}=\frac{1}{\prod_{k=n}^{0}\mathbf{g}(k)}\frac{\prod_{k=n}^{0}\mathbf{g}(k)}{\prod_{k=n}^{k(t)}\mathbf{g}(k)}\\
&=
\begin{cases}
\frac{1}{\prod_{k=n}^{0}\mathbf{g}(k)}\frac{1}{\prod_{k=1}^{k(t)}\mathbf{g}(k)}&\text{if }k(t)\ge1,\\
\frac{1}{\prod_{k=n}^{0}\mathbf{g}(k)}\prod_{k=k(t)+1}^{0}\mathbf{g}(k)&\text{if }k(t)\le0,
\end{cases}\\
&\asymp
\begin{cases}
V_\mathbf{g}(2^n)\frac{1}{V_{\mathbf{g}}(2^{k(t)})}&\text{if }k(t)\ge1,\\
V_\mathbf{g}(2^n)\frac{1}{V_{\mathbf{g}}(2^{k(t)})}&\text{if }k(t)\le0,
\end{cases}=\frac{V_\mathbf{g}(2^n)}{V_{\mathbf{g}}(2^{k(t)})}.
\end{align*}
Hence
$$\mathbb{P}\left[\omega:\mathbf{v}(\omega)(t)|_{\mathbb{Z}\cap[n,+\infty)}=\mathbf{u}|_{\mathbb{Z}\cap[n,+\infty)}\right]\lesssim \frac{V_\mathbf{g}(2^n)}{V_{\mathbf{g}}(2^{k(t)})}.$$
If $n\ge k(t)$, then
$$\mathbb{P}\left[\omega:\mathbf{v}(\omega)(t)|_{\mathbb{Z}\cap[n,+\infty)}=\mathbf{u}|_{\mathbb{Z}\cap[n,+\infty)}\right]\le1=\frac{V_\mathbf{g}(2^n)}{V_\mathbf{g}(2^n)}\le \frac{V_\mathbf{g}(2^n)}{V_{\mathbf{g}}(2^{k(t)})}.$$
In summary, we have
$$\mathbb{P}\left[\omega:\mathbf{v}(\omega)(t)\in{B}_{\mathcal{U}(\mathbf{g})}(\mathbf{u},2^n)\right]\lesssim \frac{V_\mathbf{g}(2^n)}{V_\mathbf{g}(2^{k(t)})}\asymp\frac{1}{V_{\mathbf{g}}(2^{k(t)})}\mathbf{m}_{\mathcal{U}(\mathbf{g})}(B_{\mathcal{U}(\mathbf{g})}(\mathbf{u},2^n)),$$
which gives $\mathbb{P}_{\mathbf{v}(\cdot)(t)}$ is absolutely continuous with respect to $\mathbf{m}_{\mathcal{U}(\mathbf{g})}$ with Radon derivative $\lesssim \frac{1}{V_\mathbf{g}(2^{k(t)})}$. Since for any $\omega$,
$$2^{k(t)}\asymp t\wedge (D-t)=\mathbf{d}_{\mathcal{L}(\mathbf{g},\mathbf{b})}([(\mathbf{v}(\omega)(t),\mathbf{s}_0(t))],x)\wedge \mathbf{d}_{\mathcal{L}(\mathbf{g},\mathbf{b})}([(\mathbf{v}(\omega)(t),\mathbf{s}_0(t))],y),$$
we have $V_\mathbf{g}(\mathbf{d}_{\mathcal{L}(\mathbf{g},\mathbf{b})}([(\cdot,\mathbf{s}_0(t))],x)\wedge \mathbf{d}_{\mathcal{L}(\mathbf{g},\mathbf{b})}([(\cdot,\mathbf{s}_0(t))],y))\mathbb{P}_{\mathbf{v}(\cdot)(t)}$ is absolutely continuous with respect to $\mathbf{m}_{\mathcal{U}(\mathbf{g})}$ with Radon derivative $\lesssim1$. Therefore
\begin{align*}
&\int_0^D\left(\int_\Omega (h1_B)\left([(\mathbf{v}(\omega)(t),\mathbf{s}_0(t))]\right)\mathbb{P}(\mathrm{d}\omega)\right)\mathrm{d} t\\
&=\int_0^D\left(\int_{\mathcal{U}(\mathbf{g})} (h1_B)\left([(\mathbf{u},\mathbf{s}_0(t))]\right)\mathbb{P}_{\mathbf{v}(\cdot)(t)}(\mathrm{d}\mathbf{u})\right)\mathrm{d} t\\
&=\int_0^D\left(\int_{\mathcal{U}(\mathbf{g})} (h1_B)\left([(\mathbf{u},\mathbf{s}_0(t))]\right)\right.\\
&\hspace{40pt}\left.\frac{V_\mathbf{g}(\mathbf{d}_{\mathcal{L}(\mathbf{g},\mathbf{b})}([(\mathbf{u},\mathbf{s}_0(t))],x)\wedge \mathbf{d}_{\mathcal{L}(\mathbf{g},\mathbf{b})}([(\mathbf{u},\mathbf{s}_0(t))],y))}{V_\mathbf{g}(\mathbf{d}_{\mathcal{L}(\mathbf{g},\mathbf{b})}([(\mathbf{u},\mathbf{s}_0(t))],x)\wedge \mathbf{d}_{\mathcal{L}(\mathbf{g},\mathbf{b})}([(\mathbf{u},\mathbf{s}_0(t))],y))} \mathbb{P}_{\mathbf{v}(\cdot)(t)}(\mathrm{d}\mathbf{u})\right)\mathrm{d} t\\
&\lesssim \int_0^D\left(\int_{\mathcal{U}(\mathbf{g})}\frac{(h1_B)\left([(\mathbf{u},\mathbf{s}_0(t))]\right)}{V_\mathbf{g}(\mathbf{d}_{\mathcal{L}(\mathbf{g},\mathbf{b})}([(\mathbf{u},\mathbf{s}_0(t))],x)\wedge \mathbf{d}_{\mathcal{L}(\mathbf{g},\mathbf{b})}([(\mathbf{u},\mathbf{s}_0(t))],y))}\mathbf{m}_{\mathcal{U}(\mathbf{g})}(\mathrm{d}\mathbf{u})\right)\mathrm{d} t\\
&=\int_{\mathcal{U}(\mathbf{g})}\left(\int_0^D\frac{(h1_B)\left([(\mathbf{u},\mathbf{s}_0(t))]\right)}{V_\mathbf{g}(\mathbf{d}_{\mathcal{L}(\mathbf{g},\mathbf{b})}([(\mathbf{u},\mathbf{s}_0(t))],x)\wedge \mathbf{d}_{\mathcal{L}(\mathbf{g},\mathbf{b})}([(\mathbf{u},\mathbf{s}_0(t))],y))}\mathrm{d} t\right)\mathbf{m}_{\mathcal{U}(\mathbf{g})}(\mathrm{d}\mathbf{u})\\
&\le\int_{\mathcal{U}(\mathbf{g})}\left(\int_{\mathcal{T}(\mathbf{b})}\frac{(h1_B)\left([(\mathbf{u},\mathbf{t})]\right)}{V_\mathbf{g}(\mathbf{d}_{\mathcal{L}(\mathbf{g},\mathbf{b})}([(\mathbf{u},\mathbf{t})],x)\wedge \mathbf{d}_{\mathcal{L}(\mathbf{g},\mathbf{b})}([(\mathbf{u},\mathbf{t})],y))}\lambda_{\mathcal{T}(\mathbf{b})}(\mathrm{d}\mathbf{t})\right)\mathbf{m}_{\mathcal{U}(\mathbf{g})}(\mathrm{d}\mathbf{u})\\
&=\int_{B}\frac{h(z)}{V_\mathbf{g}(\mathbf{d}_{\mathcal{L}(\mathbf{g},\mathbf{b})}(x,z))\wedge V_\mathbf{g}(\mathbf{d}_{\mathcal{L}(\mathbf{g},\mathbf{b})}(y,z))}\lambda_{\mathcal{L}(\mathbf{g},\mathbf{b})}(\mathrm{d} z).
\end{align*}
\end{proof}

We have the Poincar\'e inequality on $\mathcal{L}(\mathbf{g},\mathbf{b})$ as follows.

\begin{proposition}\label{prop_Laa_PI}
There exists $C>0$ such that for any ball $B_{\mathcal{L}(\mathbf{g},\mathbf{b})}(x_0,r)$, for any $f\in\mathcal{F}^\mathcal{L}$, we have
\begin{align*}
&\int_{B_{\mathcal{L}(\mathbf{g},\mathbf{b})}(x_0,r)}\lvert f-\dashint_{B_{\mathcal{L}(\mathbf{g},\mathbf{b})}(x_0,r)}f\mathrm{d}\mathbf{m}_{\mathcal{L}(\mathbf{g},\mathbf{b})}\rvert^p\mathrm{d}\mathbf{m}_{\mathcal{L}(\mathbf{g},\mathbf{b})}\\
&\le C\left(r^{p-1}V_\mathbf{b}(r)\right)\int_{B_{\mathcal{L}(\mathbf{g},\mathbf{b})}(x_0,4r)}|\nabla_{\mathcal{L}}f|^p\mathrm{d} \lambda_{\mathcal{L}(\mathbf{g},\mathbf{b})}.
\end{align*}
\end{proposition}

\begin{proof}
For notational convenience, let us denote $B(\cdot,\cdot)=B_{\mathcal{L}(\mathbf{g},\mathbf{b})}(\cdot,\cdot)$, $d=\mathbf{d}_{\mathcal{L}(\mathbf{g},\mathbf{b})}$, $m=\mathbf{m}_{\mathcal{L}(\mathbf{g},\mathbf{b})}$, $\lambda=\lambda_{\mathcal{L}(\mathbf{g},\mathbf{b})}$, $\nabla=\nabla_\mathcal{L}$. Since $m(B(x_0,r))\asymp V_\mathbf{g}(r)V_\mathbf{b}(r)$, we have
\begin{align*}
&\int_{B(x_0,r)}\lvert f-\dashint_{B(x_0,r)}f\mathrm{d} m\rvert^p\mathrm{d} m\le\int_{B(x_0,r)}\dashint_{B(x_0,r)}|f(x)-f(y)|^pm(\mathrm{d} y)m(\mathrm{d} x)\\
&\asymp\frac{1}{V_\mathbf{g}(r)V_\mathbf{b}(r)}\int_{B(x_0,r)}\int_{B(x_0,r)}|f(x)-f(y)|^pm(\mathrm{d} y)m(\mathrm{d} x).
\end{align*}
For any distinct $x,y\in B(x_0,r)$, let $\Gamma_{x,y}$, $\mathbb{P}^{x,y}$ be given by Proposition \ref{prop_pencil}, then for any $\gamma\in\Gamma_{x,y}$, we have
\begin{align*}
&|f(x)-f(y)|^p\le\left(\int_0^{d(x,y)}|\nabla f(\gamma(t))|\mathrm{d} t\right)^p\\
&\le d(x,y)^{p-1}\int_0^{d(x,y)}|\nabla f(\gamma(t))|^p\mathrm{d} t\le (2r)^{p-1}\int_0^{d(x,y)}|\nabla f(\gamma(t))|^p\mathrm{d} t.
\end{align*}
Taking expectation over $\gamma\in\Gamma_{x,y}$ under $\mathbb{P}^{x,y}$, by Proposition \ref{prop_pencil}, we have
\begin{align*}
&|f(x)-f(y)|^p\le(2r)^{p-1}\int_{\Gamma_{x,y}}\left(\int_0^{d(x,y)}|\nabla f(\gamma(t))|^p\mathrm{d} t\right)\mathrm{d} \mathbb{P}^{x,y}\\
&\lesssim r^{p-1}\int_{B(x,d(x,y))}\frac{|\nabla f(z)|^p}{V_\mathbf{g}(d(x,z))\wedge V_\mathbf{g}(d(y,z))}\lambda(\mathrm{d} z).
\end{align*}
Hence
\begin{align*}
&\int_{B(x_0,r)}\lvert f-\dashint_{B(x_0,r)}f\mathrm{d} m\rvert^p\mathrm{d} m\\
&\lesssim\frac{r^{p-1}}{V_\mathbf{g}(r)V_\mathbf{b}(r)}\int_{B(x_0,r)}\int_{B(x_0,r)}\left(\int_{B(x,d(x,y))}\frac{|\nabla f(z)|^p}{V_\mathbf{g}(d(x,z))\wedge V_\mathbf{g}(d(y,z))}\lambda(\mathrm{d} z)\right)m(\mathrm{d} x)m(\mathrm{d} y)\\
&\le\frac{r^{p-1}}{V_\mathbf{g}(r)V_\mathbf{b}(r)}\int_{B(x_0,4r)}\int_{B(x_0,4r)}\left(\int_{B(x_0,4r)}\frac{|\nabla f(z)|^p}{V_\mathbf{g}(d(x,z))\wedge V_\mathbf{g}(d(y,z))}\lambda(\mathrm{d} z)\right)m(\mathrm{d} x)m(\mathrm{d} y)\\
&=\frac{r^{p-1}}{V_\mathbf{g}(r)V_\mathbf{b}(r)}\\
&\hspace{10pt}\cdot\int_{B(x_0,4r)}|\nabla f(z)|^p\left(\int_{B(x_0,4r)}\int_{B(x_0,4r)}\frac{1}{V_\mathbf{g}(d(x,z))\wedge V_\mathbf{g}(d(y,z))}m(\mathrm{d} x)m(\mathrm{d} y)\right)\lambda(\mathrm{d} z),
\end{align*}
where for any $z\in B(x_0,4r)$, we have
\begin{align*}
&\int_{B(x_0,4r)}\int_{B(x_0,4r)}\frac{1}{V_\mathbf{g}(d(x,z))\wedge V_\mathbf{g}(d(y,z))}m(\mathrm{d} x)m(\mathrm{d} y)\\
&\le\int_{B(x_0,4r)}\int_{B(x_0,4r)}\frac{1_{d(x,z)\le d(y,z)}+1_{{d(x,z)\ge d(y,z)}}}{V_\mathbf{g}(d(x,z))\wedge V_\mathbf{g}(d(y,z))}m(\mathrm{d} y)m(\mathrm{d} x)\\
&=2\int_{B(x_0,4r)}\frac{m(B(x_0,4r)\backslash B(z,d(x,z)))}{V_\mathbf{g}(d(x,z))}m(\mathrm{d} x)\le 2m(B(x_0,4r))\int_{B(z,8r)}\frac{m(\mathrm{d} x)}{V_\mathbf{g}(d(x,z))}\\
&=2m(B(x_0,4r))\sum_{n\le3}\int_{B(z,2^nr)\backslash B(z,2^{n-1}r)}\frac{m(\mathrm{d} x)}{V_\mathbf{g}(d(x,z))}\le2m(B(x_0,4r))\sum_{n\le3}\frac{m(B(z,2^nr))}{V_\mathbf{g}(2^{n-1}r)}\\
&\asymp V_\mathbf{g}(r)V_\mathbf{b}(r)\sum_{n\le3}\frac{V_\mathbf{g}(2^nr)V_\mathbf{b}(2^nr)}{V_\mathbf{g}(2^{n-1}r)}\asymp V_\mathbf{g}(r)V_\mathbf{b}(r)\sum_{n\le3}V_\mathbf{b}(2^nr)\asymp V_\mathbf{g}(r)V_\mathbf{b}(r)^2,
\end{align*}
which gives
\begin{align*}
&\int_{B(x_0,r)}\lvert f-\dashint_{B(x_0,r)}f\mathrm{d} m\rvert^p\mathrm{d} m\\
&\lesssim \frac{r^{p-1}}{V_\mathbf{g}(r)V_\mathbf{b}(r)}\int_{B(x_0,4r)}|\nabla f(z)|^pV_\mathbf{g}(r)V_\mathbf{b}(r)^2\lambda(\mathrm{d} z)\\
&=r^{p-1}V_\mathbf{b}(r)\int_{B(x_0,4r)}|\nabla f(z)|^p\lambda(\mathrm{d} z).
\end{align*}
\end{proof}

Finally, we have the capacity upper bound and the cutoff Sobolev inequality on $\mathcal{L}(\mathbf{g},\mathbf{b})$ as follows.

\begin{proposition}\label{prop_Laa_CS}
There exist $C_1, C_2, C_3>0$ such that for any ball $B_{\mathcal{L}(\mathbf{g},\mathbf{b})}$ with radius $r$, there exists $\phi_\mathcal{L}\in\mathcal{C}^\mathcal{L}\subseteq\mathcal{F}^\mathcal{L}\cap C_c(\mathcal{L}(\mathbf{g},\mathbf{b}))$ with $0\le \phi_\mathcal{L}\le1$ in $\mathcal{L}(\mathbf{g},\mathbf{b})$, $\phi_\mathcal{L}=1$ in $B_{\mathcal{L}(\mathbf{g},\mathbf{b})}$, $\phi_\mathcal{L}=0$ on ${\mathcal{L}(\mathbf{g},\mathbf{b})}\backslash (256B_{{\mathcal{L}(\mathbf{g},\mathbf{b})}})$ such that
$$\mathcal{E}^\mathcal{L}(\phi_\mathcal{L})\le C_1\frac{\mathbf{m}_{\mathcal{L}(\mathbf{g},\mathbf{b})}(B_{\mathcal{L}(\mathbf{g},\mathbf{b})})}{r^{p-1}V_\mathbf{b}(r)}.$$
Hence
$$\mathrm{cap}_{{\mathcal{L}(\mathbf{g},\mathbf{b})}}(B_{\mathcal{L}(\mathbf{g},\mathbf{b})},\mathcal{L}(\mathbf{g},\mathbf{b})\backslash(256B_{\mathcal{L}(\mathbf{g},\mathbf{b})}))\le C_1\frac{\mathbf{m}_{\mathcal{L}(\mathbf{g},\mathbf{b})}(B_{\mathcal{L}(\mathbf{g},\mathbf{b})})}{r^{p-1}V_\mathbf{b}(r)}.$$
Moreover, for any $f\in\mathcal{F}^\mathcal{L}$, we have
\begin{align*}
&\int_{256B_{\mathcal{L}(\mathbf{g},\mathbf{b})}}|\widetilde{f}|^p|\nabla_\mathcal{L}\phi_\mathcal{L}|^p\mathrm{d}\lambda_{\mathcal{L}(\mathbf{g},\mathbf{b})}\\
&\le C_2\int_{256B_{\mathcal{L}(\mathbf{g},\mathbf{b})}}|\nabla_\mathcal{L} f|^p\mathrm{d}\lambda_{\mathcal{L}(\mathbf{g},\mathbf{b})}+\frac{C_3}{r^{p-1}V_\mathbf{b}(r)}\int_{256B_{\mathcal{L}(\mathbf{g},\mathbf{b})}}|f|^p\mathrm{d}\mathbf{m}_{\mathcal{L}(\mathbf{g},\mathbf{b})},
\end{align*}
where $\widetilde{f}$ is a quasi-continuous modification of $f$, such that $\widetilde{f}$ is uniquely determined\\
\noindent $|\nabla_\mathcal{L}\phi_\mathcal{L}|^p\mathrm{d}\lambda_{\mathcal{L}(\mathbf{g},\mathbf{b})}$-a.e. in $X$.
\end{proposition}

\begin{proof}
Write $B_{\mathcal{L}(\mathbf{g},\mathbf{b})}=B_{\mathcal{L}(\mathbf{g},\mathbf{b})}([(\mathbf{u},\mathbf{t})],r)$, where $(\mathbf{u},\mathbf{t})\in\mathcal{P}(\mathbf{g},\mathbf{b})$. Let $n$ be the integer satisfying $2^{n-1}\le r< 2^n$. By Lemma \ref{lem_DS}, we have the following dichotomies. For $B_{\mathcal{L}(\mathbf{g},\mathbf{b})}([(\mathbf{u},\mathbf{t})],r)$, either $B_{\mathcal{T}(\mathbf{b})}(\mathbf{t},r)\cap\cup_{k\ge n+2}\mathcal{W}_k=\emptyset$, then
\begin{equation}\label{eq_dic1r1}
\mathcal{Q}^{-1}(B_{\mathcal{L}(\mathbf{g},\mathbf{b})}([(\mathbf{u},\mathbf{t})],r))\subseteq B_{\mathcal{U}(\mathbf{g})}(\mathbf{u},8r)\times B_{\mathcal{T}(\mathbf{b})}(\mathbf{t},r),
\end{equation}
or $B_{\mathcal{T}(\mathbf{b})}(\mathbf{t},r)\cap\cup_{k\ge n+2}\mathcal{W}_k=\{\mathbf{t}^{(1)}\}$, where $\mathbf{t}^{(1)}\in\mathcal{W}_m$ for some $m\ge n+2$, then
\begin{equation}\label{eq_dic1r2}
\mathcal{Q}^{-1}(B_{{\mathcal{L}(\mathbf{g},\mathbf{b})}}([(\mathbf{u},\mathbf{t})],r))\subseteq\bigcup_{k=1}^{\mathbf{g}(m)}B_{\mathcal{U}(\mathbf{g})}(\overline{\mathbf{u}}^{(k)},8r)\times B_{\mathcal{T}(\mathbf{b})}(\mathbf{t},r),
\end{equation}
where $\overline{\mathbf{u}}^{(1)}$, \ldots, $\overline{\mathbf{u}}^{(\mathbf{g}(m))}\in\mathcal{U}(\mathbf{g})$ are all the points satisfying that $[(\overline{\mathbf{u}}^{(1)},\mathbf{t}^{(1)})]=\ldots=[(\overline{\mathbf{u}}^{(\mathbf{g}(m))},\mathbf{t}^{(1)})]=[(\mathbf{u},\mathbf{t}^{(1)})]$.

For $B_{\mathcal{L}(\mathbf{g},\mathbf{b})}([(\mathbf{u},\mathbf{t})],8r)$, either $B_{\mathcal{T}(\mathbf{b})}(\mathbf{t},8r)\cap\cup_{k\ge n+5}\mathcal{W}_k=\emptyset$, then
$$\mathcal{Q}^{-1}(B_{\mathcal{L}(\mathbf{g},\mathbf{b})}([(\mathbf{u},\mathbf{t})],8r))\subseteq B_{\mathcal{U}(\mathbf{g})}(\mathbf{u},64r)\times B_{\mathcal{T}(\mathbf{b})}(\mathbf{t},8r)\subseteq\mathcal{Q}^{-1}(B_{{\mathcal{L}(\mathbf{g},\mathbf{b})}}([(\mathbf{u},\mathbf{t})],256r)),$$
or $B_{\mathcal{T}(\mathbf{b})}(\mathbf{t},8r)\cap\cup_{k\ge n+5}\mathcal{W}_k=\{\mathbf{s}^{(1)}\}$, where $\mathbf{s}^{(1)}\in\mathcal{W}_l$ for some $l\ge n+5$, then
\begin{align}
&\mathcal{Q}^{-1}(B_{{\mathcal{L}(\mathbf{g},\mathbf{b})}}([(\mathbf{u},\mathbf{t})],8r))\nonumber\\
&\subseteq\bigcup_{k=1}^{\mathbf{g}(l)}B_{\mathcal{U}(\mathbf{g})}\left(\widetilde{\mathbf{u}}^{(k)},64r\right)\times B_{\mathcal{T}(\mathbf{b})}(\mathbf{t},8r)\nonumber\\
&\subseteq\mathcal{Q}^{-1}(B_{{\mathcal{L}(\mathbf{g},\mathbf{b})}}([(\mathbf{u},\mathbf{t})],256r))\label{eq_dic2r2},
\end{align}
where $\widetilde{\mathbf{u}}^{(1)}$, \ldots, $\widetilde{\mathbf{u}}^{(\mathbf{g}(l))}\in\mathcal{U}(\mathbf{g})$ are all the points satisfying that $[(\widetilde{\mathbf{u}}^{(1)},\mathbf{s}^{(1)})]=\ldots=[(\widetilde{\mathbf{u}}^{(\mathbf{g}(l))},\mathbf{s}^{(1)})]=[(\mathbf{u},\mathbf{s}^{(1)})]$. To write in a unified way, denote $\mathbf{u}^{(1)}$, \ldots, $\mathbf{u}^{(\sup_\mathbb{Z}\mathbf{g})}\in\mathcal{U}(\mathbf{g})$, such that
\begin{align}
&\mathcal{Q}^{-1}(B_{\mathcal{L}(\mathbf{g},\mathbf{b})}([(\mathbf{u},\mathbf{t})],8r))\nonumber\\
&\subseteq\bigcup_{k=1}^{\sup_\mathbb{Z}\mathbf{g}}B_{\mathcal{U}(\mathbf{g})}\left(\mathbf{u}^{(k)},64r\right)\times B_{\mathcal{T}(\mathbf{b})}\left(\mathbf{t},8r\right)\nonumber\\
&\subseteq\mathcal{Q}^{-1}(B_{\mathcal{L}(\mathbf{g},\mathbf{b})}([(\mathbf{u},\mathbf{t})],256r))\label{eq_dic2r}.
\end{align}
Then $\{\mathbf{u}^{(1)},\ldots,\mathbf{u}^{(\sup_\mathbb{Z}\mathbf{g})}\}$ is either the one-point set $\{\mathbf{u}\}$ or the set $\{\widetilde{\mathbf{u}}^{(1)},\ldots,\widetilde{\mathbf{u}}^{(\mathbf{g}(l))}\}$. It is possible that $\{\overline{\mathbf{u}}^{(\cdot)}\}$, $\{\widetilde{\mathbf{u}}^{(\cdot)}\}$ are different, but $\{\overline{\mathbf{u}}^{(\cdot)}\}$, $\{\widetilde{\mathbf{u}}^{(\cdot)}\}$ both contain $\mathbf{u}$.

We claim that
\begin{equation}\label{eq_claim_1r}
\mathcal{Q}^{-1}(B_{\mathcal{L}(\mathbf{g},\mathbf{b})}([(\mathbf{u},\mathbf{t})],r))\subseteq\bigcup_{k=1}^{\sup_\mathbb{Z}\mathbf{g}}B_{\mathcal{U}(\mathbf{g})}\left(\mathbf{u}^{(k)},64r\right)\times B_{\mathcal{T}(\mathbf{b})}(\mathbf{t},r).
\end{equation}
Indeed, we consider the dichotomy for $B_{\mathcal{L}(\mathbf{g},\mathbf{b})}([(\mathbf{u},\mathbf{t})],r)$. For the first case, by Equation (\ref{eq_dic1r1}), we have
$$\mathcal{Q}^{-1}(B_{\mathcal{L}(\mathbf{g},\mathbf{b})}([(\mathbf{u},\mathbf{t})],r))\subseteq B_{\mathcal{U}(\mathbf{g})}(\mathbf{u},8r)\times B_{\mathcal{T}(\mathbf{b})}(\mathbf{t},r)\subseteq\bigcup_{k=1}^{\sup_\mathbb{Z}\mathbf{g}}B_{\mathcal{U}(\mathbf{g})}\left(\mathbf{u}^{(k)},64r\right)\times B_{\mathcal{T}(\mathbf{b})}(\mathbf{t},r),$$
because $\{\mathbf{u}^{(\cdot)}\}$ always contains $\mathbf{u}$. For the second case, if $m\ge n+5$, then
$$B_{\mathcal{T}(\mathbf{b})}(\mathbf{t},8r)\cap\cup_{k\ge n+5}\mathcal{W}_k\supseteq B_{\mathcal{T}(\mathbf{b})}(\mathbf{t},r)\cap\cup_{k\ge n+5}\mathcal{W}_k=\{\mathbf{t}^{(1)}\},$$
we are in the second case of the dichotomy for $B_{{\mathcal{L}(\mathbf{g},\mathbf{b})}}([(\mathbf{u},\mathbf{t})],8r))$, $\mathbf{t}^{(1)}=\mathbf{s}^{(1)}$, $m=l$, and $\{\overline{\mathbf{u}}^{(\cdot)}\}$, $\{\widetilde{\mathbf{u}}^{(\cdot)}\}$ coincide. Then by Equation (\ref{eq_dic1r2}), we have
\begin{align*}
&\mathcal{Q}^{-1}(B_{{\mathcal{L}(\mathbf{g},\mathbf{b})}}([(\mathbf{u},\mathbf{t})],r))\subseteq\bigcup_{k=1}^{\mathbf{g}(m)}B_{\mathcal{U}(\mathbf{g})}(\overline{\mathbf{u}}^{(k)},8r)\times B_{\mathcal{T}(\mathbf{b})}(\mathbf{t},r)\\
&=\bigcup_{k=1}^{\mathbf{g}(l)}B_{\mathcal{U}(\mathbf{g})}(\widetilde{\mathbf{u}}^{(k)},8r)\times B_{\mathcal{T}(\mathbf{b})}(\mathbf{t},r)\subseteq\bigcup_{k=1}^{\mathbf{g}(l)}B_{\mathcal{U}(\mathbf{g})}(\widetilde{\mathbf{u}}^{(k)},64r)\times B_{\mathcal{T}(\mathbf{b})}(\mathbf{t},r)\\
&=\bigcup_{k=1}^{\sup_\mathbb{Z}\mathbf{g}}B_{\mathcal{U}(\mathbf{g})}\left(\mathbf{u}^{(k)},64r\right)\times B_{\mathcal{T}(\mathbf{b})}(\mathbf{t},r).
\end{align*}
If $m\in\{n+2,n+3,n+4\}$, then for any $\overline{\mathbf{u}}^{(k)}$, we have $\mathbf{d}_{\mathcal{U}(\mathbf{g})}(\overline{\mathbf{u}}^{(k)},\mathbf{u})\le2^m\le2^{n+4}$, hence
\begin{align*}
&\bigcup_{k=1}^{\mathbf{g}(m)}B_{\mathcal{U}(\mathbf{g})}(\overline{\mathbf{u}}^{(k)},8r)\subseteq\bigcup_{k=1}^{\mathbf{g}(m)} B_{\mathcal{U}(\mathbf{g})}({\mathbf{u}},\mathbf{d}_{\mathcal{U}(\mathbf{g})}(\overline{\mathbf{u}}^{(k)},\mathbf{u})+8r)\\
&\subseteq B_{\mathcal{U}(\mathbf{g})}({\mathbf{u}},2^{n+4}+8r)\subseteq B_{\mathcal{U}(\mathbf{g})}({\mathbf{u}},40r)\subseteq\bigcup_{k=1}^{\sup_\mathbb{Z}\mathbf{g}}B_{\mathcal{U}(\mathbf{g})}\left(\mathbf{u}^{(k)},64r\right),
\end{align*}
also because $\{\mathbf{u}^{(\cdot)}\}$ always contains $\mathbf{u}$. Then by Equation (\ref{eq_dic1r2}), we have
$$\mathcal{Q}^{-1}(B_{\mathcal{L}(\mathbf{g},\mathbf{b})}([(\mathbf{u},\mathbf{t})],r))\subseteq\bigcup_{k=1}^{\sup_\mathbb{Z}\mathbf{g}}B_{\mathcal{U}(\mathbf{g})}\left(\mathbf{u}^{(k)},64r\right)\times B_{\mathcal{T}(\mathbf{b})}(\mathbf{t},r).$$
Therefore, we finish the proof of Equation (\ref{eq_claim_1r}).

By the capacity upper bound on $\mathcal{T}(\mathbf{b})$, Proposition \ref{prop_tree_CS}, there exists $\phi_\mathcal{T}\in\mathcal{F}^\mathcal{T}\cap C_c(\mathcal{T}(\mathbf{b}))$ with $0\le\phi_\mathcal{T}\le1$ in $\mathcal{T}(\mathbf{b})$, $\phi_\mathcal{T}=1$ in $B_{\mathcal{T}(\mathbf{b})}(\mathbf{t},r)$, $\phi_\mathcal{T}=0$ on $\mathcal{T}(\mathbf{b})\backslash B_{\mathcal{T}(\mathbf{b})}(\mathbf{t},8r)$ such that $\mathcal{E}^\mathcal{T}(\phi_\mathcal{T})\lesssim \frac{1}{r^{p-1}}$. Let $\phi_\mathcal{L}:\mathcal{L}(\mathbf{g},\mathbf{b})\to\mathbb{R}$ be given by
$$
\phi_\mathcal{L}([(\mathbf{v},\cdot)])=
\begin{cases}
\phi_\mathcal{T}&\text{if }\mathbf{v}\in \bigcup_{k=1}^{\sup_\mathbb{Z}\mathbf{g}}B_{\mathcal{U}(\mathbf{g})}\left(\mathbf{u}^{(k)},64r\right),\\
0&\text{if }\mathbf{v}\not\in \bigcup_{k=1}^{\sup_\mathbb{Z}\mathbf{g}}B_{\mathcal{U}(\mathbf{g})}\left(\mathbf{u}^{(k)},64r\right).
\end{cases}
$$

Firstly, we show that ${\phi_\mathcal{L}}$ is well-defined. We only need to show that if
$$\mathbf{v}^{(1)}\in\bigcup_{k=1}^{\sup_\mathbb{Z}\mathbf{g}}B_{\mathcal{U}(\mathbf{g})}\left(\mathbf{u}^{(k)},64r\right),\mathbf{v}^{(2)}\not\in\bigcup_{k=1}^{\sup_\mathbb{Z}\mathbf{g}}B_{\mathcal{U}(\mathbf{g})}\left(\mathbf{u}^{(k)},64r\right),$$
satisfy $[(\mathbf{v}^{(1)},\mathbf{s})]=[(\mathbf{v}^{(2)},\mathbf{s})]$, then $\phi_\mathcal{T}(\mathbf{s})=0$. By assumption, we have $\mathbf{d}_{\mathcal{U}(\mathbf{d})}(\mathbf{v}^{(1)},\mathbf{v}^{(2)})\ge64r\ge2^{n+5}$, hence $\mathbf{s}\in\cup_{k\ge n+5}\mathcal{W}_k$. Suppose $\mathbf{s}\in B_{\mathcal{T}(\mathbf{b})}(\mathbf{t},8r)$, then we are in the second case of the dichotomy for $B_{\mathcal{L}(\mathbf{g},\mathbf{b})}([(\mathbf{u},\mathbf{t})],8r)$, hence $\mathbf{s}=\mathbf{s}^{(1)}$, $\{\mathbf{v}^{(1)},\mathbf{v}^{(2)}\}\subseteq\{\widetilde{\mathbf{u}}^{(1)},\ldots,\widetilde{\mathbf{u}}^{(\mathbf{g}(l))}\}$, in particular, we have
$$\mathbf{v}^{(2)}\in\bigcup_{k=1}^{\mathbf{g}(l)}B_{\mathcal{U}(\mathbf{g})}\left(\widetilde{{\mathbf{u}}}^{(k)},64r\right)=\bigcup_{k=1}^{\sup_\mathbb{Z}\mathbf{g}}B_{\mathcal{U}(\mathbf{g})}\left(\mathbf{u}^{(k)},64r\right),$$
contradiction. Hence $\mathbf{s}\not\in B_{\mathcal{T}(\mathbf{b})}(\mathbf{t},8r)$, $\phi_\mathcal{T}(\mathbf{s})=0$. Therefore, $\phi_\mathcal{L}$ is well-defined.

Secondly, it is obvious that $\phi_\mathcal{L}\in C_c(\mathcal{L}(\mathbf{g},\mathbf{b}))$. We claim that $\phi_\mathcal{L}\in\mathcal{C}^\mathcal{L}_{n+4}\subseteq\mathcal{C}^\mathcal{L}$. Indeed, let $\widetilde{\phi}_\mathcal{L}=\phi_\mathcal{L}\circ\mathcal{Q}\in C_c(\mathcal{P}(\mathbf{g},\mathbf{b}))$. For any $\mathbf{v}\in\mathcal{U}(\mathbf{g})$, $\widetilde{\phi}_\mathcal{L}(\mathbf{v},\cdot)$ is either $\phi_\mathcal{T}$ or $0$, hence in $\mathcal{F}^\mathcal{T}$. For any $\mathbf{s}\in\mathcal{T}(\mathbf{b})$, for any $\mathbf{v}^{(1)}$, $\mathbf{v}^{(2)}\in\mathcal{U}(\mathbf{g})$ with $\mathbf{d}_{\mathcal{U}(\mathbf{g})}(\mathbf{v}^{(1)},\mathbf{v}^{(2)})\le2^{n+4}$, suppose
$$\mathbf{v}^{(1)}\in\bigcup_{k=1}^{\sup_\mathbb{Z}\mathbf{g}}B_{\mathcal{U}(\mathbf{g})}\left(\mathbf{u}^{(k)},64r\right),\mathbf{v}^{(2)}\not\in\bigcup_{k=1}^{\sup_\mathbb{Z}\mathbf{g}}B_{\mathcal{U}(\mathbf{g})}\left(\mathbf{u}^{(k)},64r\right),$$
then $\mathbf{d}_{\mathcal{U}(\mathbf{g})}(\mathbf{v}^{(1)},\mathbf{v}^{(2)})\ge64r\ge2^{n+5}>2^{n+4}\ge \mathbf{d}_{\mathcal{U}(\mathbf{g})}(\mathbf{v}^{(1)},\mathbf{v}^{(2)})$, contradiction. Hence $\widetilde{\phi}_\mathcal{L}(\mathbf{v}^{(1)},\mathbf{s})=\widetilde{\phi}_\mathcal{L}(\mathbf{v}^{(2)},\mathbf{s})$ is either $\phi_\mathcal{T}(\mathbf{s})$ or $0$, that is, $\widetilde{\phi}_\mathcal{L}(\cdot,\mathbf{s})$ is constant on any closed ball with radius $2^{n+4}$. Therefore $\phi_\mathcal{L}\in\mathcal{C}^\mathcal{L}_{n+4}\subseteq\mathcal{C}^\mathcal{L}$.

Finally, it is obvious that $0\le\phi_\mathcal{L}\le1$ in $\mathcal{L}(\mathbf{g},\mathbf{b})$. By Equation (\ref{eq_dic2r}), we have $\phi_\mathcal{L}=0$ on $\mathcal{L}(\mathbf{g},\mathbf{b})\backslash B_{\mathcal{L}(\mathbf{g},\mathbf{b})}([(\mathbf{u},\mathbf{t})],256r)$. By Equation (\ref{eq_claim_1r}), we have $\phi_\mathcal{L}=1$ in $B_{\mathcal{L}(\mathbf{g},\mathbf{b})}([(\mathbf{u},\mathbf{t})],r)$. Moreover
$$\mathcal{E}^\mathcal{L}(\phi_\mathcal{L})\asymp V_\mathbf{g}(r)\mathcal{E}^\mathcal{T}(\phi_\mathcal{T})\lesssim V_\mathbf{g}(r)\frac{1}{r^{p-1}}=\frac{V_\mathbf{g}(r)V_\mathbf{b}(r)}{r^{p-1}V_\mathbf{b}(r)}\asymp \frac{\mathbf{m}_{\mathcal{L}(\mathbf{g},\mathbf{b})}(B_{\mathcal{L}(\mathbf{g},\mathbf{b})}([(\mathbf{u},\mathbf{t})],r))}{r^{p-1}V_{\mathbf{b}}(r)},$$
hence
$$\mathrm{cap}_{\mathcal{L}(\mathbf{g},\mathbf{b})}(B_{\mathcal{L}(\mathbf{g},\mathbf{b})}([(\mathbf{u},\mathbf{t})],r),\mathcal{L}(\mathbf{g},\mathbf{b})\backslash B_{\mathcal{L}(\mathbf{g},\mathbf{b})}([(\mathbf{u},\mathbf{t})],256r))\lesssim\frac{\mathbf{m}_{\mathcal{L}(\mathbf{g},\mathbf{b})}(B_{\mathcal{L}(\mathbf{g},\mathbf{b})}([(\mathbf{u},\mathbf{t})],r))}{r^{p-1}V_{\mathbf{b}}(r)}.$$
Moreover, for any $f\in\mathcal{C}^\mathcal{L}$, we have $\widehat{f}=f\circ\mathcal{Q}\in C_c(\mathcal{P}(\mathbf{g},\mathbf{b}))$ satisfies $\widehat{f}(\mathbf{v},\cdot)\in\mathcal{F}^\mathcal{T}$ for any $\mathbf{v}\in\mathcal{U}(\mathbf{g})$. Let $U=\bigcup_{k=1}^{\sup_\mathbb{Z}\mathbf{g}}B_{\mathcal{U}(\mathbf{g})}\left(\mathbf{u}^{(k)},64r\right)$. By Proposition \ref{prop_tree_CS}, we have
\begin{align}
&\int_{B_{\mathcal{L}(\mathbf{g},\mathbf{b})}([(\mathbf{u},\mathbf{t})],256r)}|f|^p|\nabla_\mathcal{L}\phi_\mathcal{L}|^p\mathrm{d}\lambda_{\mathcal{L}(\mathbf{g},\mathbf{b})}\nonumber\\
&=\int_{U}\left(\int_{B_{\mathcal{T}(\mathbf{b})}(\mathbf{t},8r)}|\widehat{f}(\mathbf{v},\cdot)|^p|\nabla_\mathcal{T}\phi_\mathcal{T}|^p\mathrm{d}\lambda_{\mathcal{T}(\mathbf{b})}\right)\mathbf{m}_{\mathcal{U}(\mathbf{g})}(\mathrm{d}\mathbf{v})\nonumber\\
&\lesssim\int_U \left(\int_{B_{\mathcal{T}(\mathbf{b})}(\mathbf{t},8r)}|\nabla_\mathcal{T}\widehat{f}(\mathbf{v},\cdot)|^p\mathrm{d}\lambda_{\mathcal{T}(\mathbf{b})}\right.\nonumber\\
&\hspace{20pt}\left.+\frac{1}{r^{p-1}V_\mathbf{b}(r)}\int_{B_{\mathcal{T}(\mathbf{b})}(\mathbf{t},8r)}|\widehat{f}(\mathbf{v},\cdot)|^p\mathrm{d}\mathbf{m}_{\mathcal{T}(\mathbf{b})}\right)\mathbf{m}_{\mathcal{U}(\mathbf{g})}(\mathrm{d}\mathbf{v})\nonumber\\
&\le\int_{B_{\mathcal{L}(\mathbf{g},\mathbf{b})}([(\mathbf{u},\mathbf{t})],256r)}|\nabla_\mathcal{L} f|^p\mathrm{d}\lambda_{\mathcal{L}(\mathbf{g},\mathbf{b})}\nonumber\\
&\hspace{20pt}+\frac{1}{r^{p-1}V_\mathbf{b}(r)}\int_{B_{\mathcal{L}(\mathbf{g},\mathbf{b})}([(\mathbf{u},\mathbf{t})],256r)}|f|^p\mathrm{d}\mathbf{m}_{\mathcal{L}(\mathbf{g},\mathbf{b})}\label{eq_Laa_CScalC},
\end{align}
where the last inequality follows from Equation (\ref{eq_dic2r}). For any $f\in\mathcal{F}^\mathcal{L}$, since $\mathcal{C}^\mathcal{L}$ is $(\mathcal{E}^\mathcal{L}_1(\cdot)=\mathcal{E}^\mathcal{L}(\cdot)+\lVert \cdot\rVert^p_{L^p(\mathcal{L}(\mathbf{g},\mathbf{b});\mathbf{m}_{\mathcal{L}(\mathbf{g},\mathbf{b})})})$-dense in $\mathcal{F}^\mathcal{L}$, by Proposition \ref{prop_quasi_exist}, there exist a quasi-continuous modification $\widetilde{f}$ of $f$, $\{f_n\}\subseteq\mathcal{C}^\mathcal{L}$ which is $\mathcal{E}^\mathcal{L}_1$-convergent to $f$, and an increasing sequence of closed sets $\{F_k\}$ with $\mathrm{cap}_1\left(\mathcal{L}(\mathbf{g},\mathbf{b})\backslash F_k\right)\downarrow0$, such that $\{f_n\}$ is uniformly convergent to $\widetilde{f}$ on each $F_k$. Hence $\mathrm{cap}_1\left(\mathcal{L}(\mathbf{g},\mathbf{b})\backslash\cup_kF_k\right)=0$ and $\lim_{n\to+\infty}f_n(x)=\widetilde{f}(x)$ for any $x\in\cup_kF_k$. By Proposition \ref{prop_quasi_charge}, $\Gamma^\mathcal{L}(\phi_\mathcal{L})\left(\mathcal{L}(\mathbf{g},\mathbf{b})\backslash\cup_kF_k\right)=0$. Recall that $\mathrm{d}\Gamma^\mathcal{L}(\phi_\mathcal{L})=|\nabla_\mathcal{L}\phi_\mathcal{L}|^p\mathrm{d}\lambda_{\mathcal{L}(\mathbf{g},\mathbf{b})}$. By Fatou's lemma, we have
\begin{align*}
&\int_{B_{\mathcal{L}(\mathbf{g},\mathbf{b})}([(\mathbf{u},\mathbf{t})],256r)}|\widetilde{f}|^p|\nabla_\mathcal{L}\phi_\mathcal{L}|^p\mathrm{d}\lambda_{\mathcal{L}(\mathbf{g},\mathbf{b})}=\int_{B_{\mathcal{L}(\mathbf{g},\mathbf{b})}([(\mathbf{u},\mathbf{t})],256r)}\lim_{n\to+\infty}|f_n|^p|\nabla_\mathcal{L}\phi_\mathcal{L}|^p\mathrm{d}\lambda_{\mathcal{L}(\mathbf{g},\mathbf{b})}\\
&\le\varliminf_{n\to+\infty}\int_{B_{\mathcal{L}(\mathbf{g},\mathbf{b})}([(\mathbf{u},\mathbf{t})],256r)}|{f}_n|^p|\nabla_\mathcal{L}\phi_\mathcal{L}|^p\mathrm{d}\lambda_{\mathcal{L}(\mathbf{g},\mathbf{b})}\\
&\overset{(\star)}{\scalebox{2}[1]{$\lesssim$}}\varliminf_{n\to+\infty}\left(\int_{B_{\mathcal{L}(\mathbf{g},\mathbf{b})}([(\mathbf{u},\mathbf{t})],256r)}|\nabla_\mathcal{L} f_n|^p\mathrm{d}\lambda_{\mathcal{L}(\mathbf{g},\mathbf{b})}\right.\\
&\hspace{50pt}\left.+\frac{1}{r^{p-1}V_\mathbf{b}(r)}\int_{B_{\mathcal{L}(\mathbf{g},\mathbf{b})}([(\mathbf{u},\mathbf{t})],256r)}|f_n|^p\mathrm{d}\mathbf{m}_{\mathcal{L}(\mathbf{g},\mathbf{b})}\right)\\
&\overset{(\diamond)}{\scalebox{2}[1]{$=$}}\int_{B_{\mathcal{L}(\mathbf{g},\mathbf{b})}([(\mathbf{u},\mathbf{t})],256r)}|\nabla_\mathcal{L} f|^p\mathrm{d}\lambda_{\mathcal{L}(\mathbf{g},\mathbf{b})}+\frac{1}{r^{p-1}V_\mathbf{b}(r)}\int_{B_{\mathcal{L}(\mathbf{g},\mathbf{b})}([(\mathbf{u},\mathbf{t})],256r)}|f|^p\mathrm{d}\mathbf{m}_{\mathcal{L}(\mathbf{g},\mathbf{b})},
\end{align*}
where $(\star)$ follows from Equation (\ref{eq_Laa_CScalC}) for $\{f_n\}$, and $(\diamond)$ follows from the fact that $\{f_n\}$ is $\mathcal{E}^\mathcal{L}_1$-convergent to $f$.
\end{proof}

\section{Proof of Theorem \ref{thm_main}}\label{sec_proof_main}

The proof can be reduced to find $\mathbf{g},\mathbf{b}$ for given $\Phi,\Psi$, which is the following elementary result. 

\begin{lemma}\label{lem_findgb}
Let $\theta_1,\theta_2\ge0$ with $\theta_1\le\theta_2$, and let $C>0$. Assume that $\Phi:(0,+\infty)\to(0,+\infty)$ satisfies that for any $R,r>0$ with $r\le R$, we have
$$\frac{1}{C}\left(\frac{R}{r}\right)^{\theta_1}\le \frac{\Phi(R)}{\Phi(r)}\le C\left(\frac{R}{r}\right)^{\theta_2}.$$
Then there exists an integer-valued function $\mathbf{a}:\mathbb{Z}\to\{2^{\lfloor {\theta_1}\rfloor},2^{\lceil {\theta_2}\rceil}\}$ such that
$$\frac{1}{C^32^{\lceil {\theta_2}\rceil-\lfloor {\theta_1}\rfloor}}\prod_{k=0}^n\mathbf{a}(k)\le \frac{\Phi(2^n)}{\Phi(1)}\le{C^32^{\lceil {\theta_2}\rceil-\lfloor {\theta_1}\rfloor}}\prod_{k=0}^n\mathbf{a}(k)\text{ for any }n\ge0,$$
$$\frac{1}{C^32^{\lceil {\theta_2}\rceil-\lfloor {\theta_1}\rfloor}}\prod_{k=n}^0\mathbf{a}(k)\le \frac{\Phi(1)}{\Phi(2^n)}\le{C^32^{\lceil {\theta_2}\rceil-\lfloor {\theta_1}\rfloor}}\prod_{k=n}^0\mathbf{a}(k)\text{ for any }n\le0.$$
\end{lemma}

\begin{proof}
By replacing $\theta_1, \theta_2$ by $\lfloor{\theta_1} \rfloor, \lceil {\theta_2}\rceil$, respectively, we may assume that $\theta_1\le\theta_2$ are non-negative integers. We construct $\mathbf{a}(n)$ for $n\ge0$, a similar construction applies for $n\le0$.

Let $\mathbf{a}(0)=2^{\theta_1}$. Denote $\mathbf{A}(n)=\Pi_{k=0}^n\mathbf{a}(k)$ if $\mathbf{a}(0),\ldots,\mathbf{a}(n)$ have already been defined.

By assumption, for any $n\ge0$, we have $\frac{\Phi(2^n)}{\Phi(1)}\ge \frac{1}{C}2^{\theta_1 n}$. If $\frac{\Phi(2^n)}{\Phi(1)}\le C2^{\theta_1n}$ for any $n\ge0$, then let $\mathbf{a}(n)=2^{\theta_1}$ for any $n\ge1$, the result is obvious. Otherwise, let $n_1=\min\{n\ge0:\frac{\Phi(2^n)}{\Phi(1)}>C2^{\theta_1n}\}$, then $n_1>0$, let $\mathbf{a}(n)=2^{\theta_1}$ for any $n=1,\ldots,n_1$. Since $\frac{\Phi(2^{n_1-1})}{\Phi(1)}\le C2^{\theta_1(n_1-1)}$, we have
$$\frac{\Phi(2^{n_1})}{\Phi(1)}=\frac{\Phi(2^{n_1})}{\Phi(2^{n_1-1})}\frac{\Phi(2^{n_1-1})}{\Phi(1)}\le (C2^{\theta_2})(C2^{\theta_1(n_1-1)})=C^22^{\theta_2-\theta_1}\mathbf{A}(n_1).$$
By assumption, for any $k\ge0$, we have
$$\frac{\Phi(2^{n_1+k})}{\Phi(1)}=\frac{\Phi(2^{n_1+k})}{\Phi(2^{n_1})}\frac{\Phi(2^{n_1})}{\Phi(1)}\le(C2^{\theta_2k})\left(C^22^{\theta_2-\theta_1}\mathbf{A}(n_1)\right)=C^32^{\theta_2-\theta_1}\mathbf{A}(n_1)2^{\theta_2k}.$$
By definition, $\frac{\Phi(2^{n_1})}{\Phi(1)}>C\mathbf{A}(n_1)\ge \frac{1}{C}\mathbf{A}(n_1)$. If $\frac{\Phi(2^{n_1+k})}{\Phi(1)}\ge \frac{1}{C}\mathbf{A}(n_1)2^{\theta_2k}$ for any $k\ge0$, then let $\mathbf{a}(n)=2^{\theta_2}$ for any $n>n_1$, the result is obvious. Otherwise, let $m_1=\min\{n_1+k:k\ge0,\frac{\Phi(2^{n_1+k})}{\Phi(1)}<\frac{1}{C}\mathbf{A}(n_1)2^{\theta_2k}\}$, then $m_1>n_1$, let $\mathbf{a}(n)=2^{\theta_2}$ for any $n=n_1+1,\ldots,m_1$. Since $\frac{\Phi(2^{m_1-1})}{\Phi(1)}\ge\frac{1}{C}\mathbf{A}(n_1)2^{\theta_2(m_1-1-n_1)}$, we have
$$\frac{\Phi(2^{m_1})}{\Phi(1)}=\frac{\Phi(2^{m_1})}{\Phi(2^{m_1-1})}\frac{\Phi(2^{m_1-1})}{\Phi(1)}\ge\left(\frac{1}{C}2^{\theta_1}\right)\left(\frac{1}{C}\mathbf{A}(n_1)2^{\theta_2(m_1-1-n_1)}\right)=\frac{1}{C^22^{\theta_2-\theta_1}}\mathbf{A}(m_1).$$
By assumption, for any $k\ge0$, we have
$$\frac{\Phi(2^{m_1+k})}{\Phi(1)}=\frac{\Phi(2^{m_1+k})}{\Phi(2^{m_1})}\frac{\Phi(2^{m_1})}{\Phi(1)}\ge\left(\frac{1}{C}2^{\theta_1k}\right)\left(\frac{1}{C^22^{\theta_2-\theta_1}}\mathbf{A}(m_1)\right)=\frac{1}{C^32^{\theta_2-\theta_1}}\mathbf{A}(m_1)2^{\theta_1k}.$$
By definition, $\frac{\Phi(2^{m_1})}{\Phi(1)}<\frac{1}{C}\mathbf{A}(n_1)2^{\theta_2(m_1-n_1)}=\frac{1}{C}\mathbf{A}(m_1)\le C\mathbf{A}(m_1)$.

Assume that we have constructed $n_l, m_l$. If $\frac{\Phi(2^{m_l+k})}{\Phi(1)}\le C\mathbf{A}(m_l)2^{\theta_1k}$ for any $k\ge0$, then let $\mathbf{a}(n)=2^{\theta_1}$ for any $n>m_l$, the result is obvious. Otherwise, let $n_{l+1}=\min\{m_l+k:k\ge0,\frac{\Phi(2^{m_l+k})}{\Phi(1)}>C\mathbf{A}(m_l)2^{\theta_1k}\}$, then $n_{l+1}>m_l$, let $\mathbf{a}(n)=2^{\theta_1}$ for any $n=m_l+1,\ldots,n_{l+1}$. Since $\frac{\Phi(2^{n_{l+1}-1})}{\Phi(1)}\le C\mathbf{A}(m_l)2^{\theta_1(n_{l+1}-1-m_l)}$, we have
$$\frac{\Phi(2^{n_{l+1}})}{\Phi(1)}=\frac{\Phi(2^{n_{l+1}})}{\Phi(2^{n_{l+1}-1})}\frac{\Phi(2^{n_{l+1}-1})}{\Phi(1)}\le (C2^{\theta_2})\left(C\mathbf{A}(m_l)2^{\theta_1(n_{l+1}-1-m_l)}\right)=C^22^{\theta_2-\theta_1}\mathbf{A}(n_{l+1}).$$
By assumption, for any $k\ge0$, we have
\begin{align*}
&\frac{\Phi(2^{n_{l+1}+k})}{\Phi(1)}=\frac{\Phi(2^{n_{l+1}+k})}{\Phi(2^{n_{l+1}})}\frac{\Phi(2^{n_{l+1}})}{\Phi(1)}\\
&\le\left(C2^{\theta_2k}\right)\left(C^22^{\theta_2-\theta_1}\mathbf{A}(n_{l+1})\right)=C^32^{\theta_2-\theta_1}\mathbf{A}(n_{l+1})2^{\theta_2k}.
\end{align*}
By definition, $\frac{\Phi(2^{n_{l+1}})}{\Phi(1)}>C\mathbf{A}(m_l)2^{\theta_1(n_{l+1}-m_l)}=C\mathbf{A}(n_{l+1})\ge \frac{1}{C}\mathbf{A}(n_{l+1})$. If $\frac{\Phi(2^{n_{l+1}+k})}{\Phi(1)}\ge \frac{1}{C}\mathbf{A}(n_{l+1})2^{\theta_2k}$ for any $k\ge0$, then let $\mathbf{a}(n)=2^{\theta_2}$ for any $n>n_{l+1}$, the result is obvious. Otherwise, let $m_{l+1}=\min\{n_{l+1}+k:k\ge0,\frac{\Phi(2^{n_{l+1}+k})}{\Phi(1)}<\frac{1}{C}\mathbf{A}(n_{l+1})2^{\theta_2k}\}$, then $m_{l+1}>n_{l+1}$, let $\mathbf{a}(n)=2^{\theta_2}$ for any $n=n_{l+1}+1,\ldots,m_{l+1}$. Since $\frac{\Phi(2^{m_{l+1}-1})}{\Phi(1)}\ge\frac{1}{C}\mathbf{A}(n_{l+1})2^{\theta_2(m_{l+1}-1-n_{l+1})}$, we have
\begin{align*}
&\frac{\Phi(2^{m_{l+1}})}{\Phi(1)}=\frac{\Phi(2^{m_{l+1}})}{\Phi(2^{m_{l+1}-1})}\frac{\Phi(2^{m_{l+1}-1})}{\Phi(1)}\\
&\ge\left(\frac{1}{C}2^{\theta_1}\right)\left(\frac{1}{C}\mathbf{A}(n_{l+1})2^{\theta_2(m_{l+1}-1-n_{l+1})}\right)=\frac{1}{C^22^{\theta_2-\theta_1}}\mathbf{A}(m_{l+1}).
\end{align*}
By assumption, for any $k\ge0$, we have
\begin{align*}
&\frac{\Phi(2^{m_{l+1}+k})}{\Phi(1)}=\frac{\Phi(2^{m_{l+1}+k})}{\Phi(2^{m_{l+1}})}\frac{\Phi(2^{m_{l+1}})}{\Phi(1)}\\
&\ge\left(\frac{1}{C}2^{\theta_1k}\right)\left(\frac{1}{C^22^{\theta_2-\theta_1}}\mathbf{A}(m_{l+1})\right)=\frac{1}{C^32^{\theta_2-\theta_1}}\mathbf{A}(m_{l+1})2^{\theta_1k}.
\end{align*}
By definition, $\frac{\Phi(2^{m_{l+1}})}{\Phi(1)}<\frac{1}{C}\mathbf{A}(n_{l+1})2^{\theta_2(m_{l+1}-n_{l+1})}=\frac{1}{C}\mathbf{A}(m_{l+1})\le C\mathbf{A}(m_{l+1})$.

By the principle of induction, we obtain $\mathbf{a}(n)$ for $n\ge0$ satisfying that for any $n\ge0$
$$\frac{1}{C^32^{\theta_2-\theta_1}}\prod_{k=0}^n\mathbf{a}(k)\le \frac{\Phi(2^n)}{\Phi(1)}\le{C^32^{\theta_2-\theta_1}}\prod_{k=0}^n\mathbf{a}(k).$$
\end{proof}

Now we give the proof of Theorem \ref{thm_main} as follows.

\begin{proof}[Proof of Theorem \ref{thm_main}]
By assumption, for any $R,r>0$ with $r\le R$, we have
$$\frac{1}{C}\frac{R}{r}\le \frac{\frac{\Psi(R)}{R^{p-1}}}{\frac{\Psi(r)}{r^{p-1}}}\le C \frac{\Phi(R)}{\Phi(r)}\le CC_\Phi\left(\frac{R}{r}\right)^{\log_2C_\Phi}.$$
By Lemma \ref{lem_findgb}, there exists $\mathbf{b}:\mathbb{Z}\to\{2^1,2^{\lceil{\log_2C_\Phi} \rceil}\}$ such that
$$\frac{\frac{\Psi(2^n)}{(2^n)^{p-1}}}{\frac{\Psi(1)}{1^{p-1}}}\asymp\prod_{k=0}^n\mathbf{b}(k)\text{ for any }n\ge0,$$
$$\frac{\frac{\Psi(1)}{1^{p-1}}}{\frac{\Psi(2^n)}{(2^n)^{p-1}}}\asymp \prod_{k=n}^0\mathbf{b}(k)\text{ for any }n\le0,$$
hence $\frac{\Psi(r)}{r^{p-1}}\asymp V_\mathbf{b}(r)$, or equivalently, $\Psi(r)\asymp r^{p-1}V_\mathbf{b}(r)$.

Moreover, by assumption, for any $R,r>0$ with $r\le R$, we have
$$\frac{1}{C}\le \frac{\frac{R^{p-1}\Phi(R)}{\Psi(R)}}{\frac{r^{p-1}\Phi(r)}{\Psi(r)}}=\left(\frac{R}{r}\right)^{p-1}\frac{\Phi(R)}{\Phi(r)}\frac{\Psi(r)}{\Psi(R)}\le C_\Phi \left(\frac{R}{r}\right)^{\log_2C_\Phi+p-1}.$$
By Lemma \ref{lem_findgb} again, there exists $\mathbf{g}:\mathbb{Z}\to\{2^0,2^{\lceil {\log_2C_\Phi+p}\rceil-1}\}$ such that $\frac{r^{p-1}\Phi(r)}{\Psi(r)}\asymp V_\mathbf{g}(r)$, hence $\Phi(r)\asymp r^{1-p}\Psi(r)V_\mathbf{g}(r)\asymp V_\mathbf{g}(r)V_\mathbf{b}(r)$.

Let $(\mathcal{L}(\mathbf{g},\mathbf{b}),\mathbf{d}_{\mathcal{L}(\mathbf{g},\mathbf{b})},\mathbf{m}_{\mathcal{L}(\mathbf{g},\mathbf{b})})$ be the corresponding Laakso-type space, and $(\mathcal{E}^\mathcal{L},\mathcal{F}^\mathcal{L})$ the corresponding $p$-energy with a $p$-energy measure $\Gamma^\mathcal{L}$. Then by Proposition \ref{prop_Laa_geo}, Equation (\ref{eq_Laa_vol}), Proposition \ref{prop_Laa_PI}, and Proposition \ref{prop_Laa_CS}, we have $(\mathcal{L}(\mathbf{g},\mathbf{b}),\mathbf{d}_{\mathcal{L}(\mathbf{g},\mathbf{b})})$ is a geodesic space, and \ref{eq_VPhi}, \ref{eq_PI}, \ref{eq_CS} hold. In particular, for any $d_h,\beta_p>0$ satisfying $p\le\beta_p\le d_h+(p-1)$, let $\Phi:r\mapsto r^{d_h}$, $\Psi:r\mapsto r^{\beta_p}$, then for any $R,r>0$ with $r\le R$, we have
$$\left(\frac{R}{r}\right)^p\le\left(\frac{R}{r}\right)^{\beta_p}=\frac{\Psi(R)}{\Psi(r)}=\left(\frac{R}{r}\right)^{\beta_p-d_h}\frac{\Phi(R)}{\Phi(r)}\le\left(\frac{R}{r}\right)^{p-1}\frac{\Phi(R)}{\Phi(r)},$$
hence \hyperlink{eq_Valpha}{V($d_h$)}, \hyperlink{eq_PIbeta}{PI($\beta_p$)}, \hyperlink{eq_CSbeta}{$\text{CS}(\beta_p)$} hold.
\end{proof}

\section{Some results from potential theory}\label{sec_quasi}

The main results of this section are Proposition \ref{prop_quasi_exist} and Proposition \ref{prop_quasi_charge}, which are the only results from this section used elsewhere in this paper.

The potential theory for $p$-energies parallels the theory for Dirichlet forms in \cite{FOT11} as outlined below, see also \cite{Hin04,BV05,Cla23,SZ25} for more general settings. Let $(X,d,m)$ be a metric measure space. Let $(\mathcal{E},\mathcal{F})$ be a $p$-energy with a $p$-energy measure $\Gamma$. Let $\mathcal{E}_1(u)=\mathcal{E}(u)+\lVert u\rVert_{L^p(X;m)}^p$ for any $u\in\mathcal{F}$.

We list some results related to the convexity originated from \ref{eq_Cla} as follows.

\begin{lemma}\label{lem_quasi_KS}
We have the following results.
\begin{enumerate}[label=(\arabic*)]
\item (\cite[Proposition 3.13]{KS24a}) $(\mathcal{F},\mathcal{E}_1^{1/p})$ is a uniformly convex reflexive Banach space.
\item (\cite[Theorem 3.7, Corollary 3.25]{KS24a})
For any $f,g\in\mathcal{F}$, the derivative
$$\mathcal{E}(f;g)=\frac{1}{p}\frac{\mathrm{d}}{\mathrm{d} t}\mathcal{E}(f+tg)|_{t=0}\in\mathbb{R}$$
exists, the map $\mathcal{E}(f;\cdot):\mathcal{F}\to\mathbb{R}$ is linear, $\mathcal{E}(f;f)=\mathcal{E}(f)$. Moreover, for any $f, g\in\mathcal{F}$ and any $a\in\mathbb{R}$, we have
\begin{equation}\label{eq_quasi_strict}
\mathbb{R}\ni t\mapsto\mathcal{E}(f+tg;g)\in\mathbb{R}\text{ is strictly increasing if and only if }\mathcal{E}(g)>0,
\end{equation}
$$\mathcal{E}(af;g)=\mathrm{sgn}(a)|a|^{p-1}\mathcal{E}(f;g),$$
$$|\mathcal{E}(f;g)|\le\mathcal{E}(f)^{(p-1)/p}\mathcal{E}(g)^{1/p}.$$
Moreover, all of the above results remain valid with $\mathcal{E}$ replaced by $\mathcal{E}_1$, and for any $f,g\in\mathcal{F}$, we have
$$\mathcal{E}_1(f;g)=\mathcal{E}(f;g)+\int_X\mathrm{sgn}(f)|f|^{p-1}g\mathrm{d} m.$$
\end{enumerate}
\end{lemma}

Let $\mathcal{O}$ be the family of all open subsets of $X$. For any $A\in\mathcal{O}$, let
\begin{equation}\label{eq_quasi_LAopen}
\mathcal{L}_A=\{u\in\mathcal{F}:u\ge 1\text{ }m\text{-a.e. on }A\},
\end{equation}
and
$$\mathrm{cap}_1(A)=\inf \left\{\mathcal{E}_1(u):u\in\mathcal{L}_A\right\},$$
where $\inf\emptyset=+\infty$. Let $\mathcal{K}$ be the family of all compact subsets of $X$.

\begin{lemma}\label{lem_quasi_open}
For any $A\in\mathcal{O}$ with $\mathcal{L}_A\ne\emptyset$, there exists a unique element $e_A\in\mathcal{L}_A$ such that $\mathrm{cap}_1(A)=\mathcal{E}_1(e_A)$. Moreover, we have $0\le e_A\le 1$ $m$-a.e. in $X$, $e_A=1$ $m$-a.e. on $A$, $e_A$ is the unique element $u\in\mathcal{F}$ satisfying $u=1$ $m$-a.e. on $A$ and $\mathcal{E}_1(u;v)\ge0$ for any $v\in\mathcal{F}$ with $v\ge0$ $m$-a.e. on $A$.
\end{lemma}

\begin{proof}
It is obvious that $\mathcal{L}_A$ is a non-empty closed convex subset of $\mathcal{F}$, by \cite[Theorem 8 in Chapter 5]{Lax02}, there exists a unique element $e_A\in\mathcal{L}_A$ such that $\mathcal{E}_1(e_A)=\inf_{u\in\mathcal{L}_A}\mathcal{E}_1(u)=\mathrm{cap}_1(A)$. Since $(e_A\vee 0)\wedge1\in\mathcal{L}_A$ and $\mathcal{E}_1((e_A\vee 0)\wedge1)\le\mathcal{E}_1(e_A)$, which follows from the Markovian property, by the uniqueness of $e_A$ in $\mathcal{L}_A$, we have $e_A=(e_A\vee 0)\wedge1$, which gives $0\le e_A\le1$ $m$-a.e. in $X$ and $e_A=1$ $m$-a.e. on $A$.

For any $v\in\mathcal{F}$ with $v\ge0$ $m$-a.e. on $A$, for any $t>0$, we have $e_A+tv\in\mathcal{L}_A$, hence $\mathcal{E}_1(e_A+tv)\ge\mathcal{E}_1(e_A)$, which gives $\mathcal{E}_1(e_A;v)=\frac{1}{p}\lim_{t\downarrow0}\frac{1}{t}\left(\mathcal{E}_1(e_A+tv)-\mathcal{E}_1(e_A)\right)\ge0$.

Assume that $u\in\mathcal{F}$ satisfies $u=1$ $m$-a.e. on $A$ and $\mathcal{E}_1(u;v)\ge0$ for any $v\in\mathcal{F}$ with $v\ge0$ $m$-a.e. on $A$. Since $u-e_A=e_A-u=0$ $m$-a.e. on $A$, by assumption we have $\mathcal{E}_1(u;u-e_A)\ge0$, $\mathcal{E}_1(u;e_A-u)\ge0$, which gives $\mathcal{E}_1(u;u-e_A)=0$, similarly, we have $\mathcal{E}_1(e_A;u-e_A)=0$. Suppose that $u\ne e_A$, consider the function $\varphi(t)=\mathcal{E}_1 \left((1-t)e_A+tu\right)=\mathcal{E}_1(e_A+t(u-e_A))$, $t\in\mathbb{R}$, then $\varphi'(t)=p\mathcal{E}_1(e_A+t(u-e_A);u-e_A)$, hence $\varphi'(0)=p\mathcal{E}_1(e_A;u-e_A)=0$, $\varphi'(1)=p\mathcal{E}_1(u;u-e_A)=0$. However, since $\mathcal{E}_1(u-e_A)>0$, by (\ref{eq_quasi_strict}) for $\mathcal{E}_1$, we have $\varphi'$ is strictly increasing, contradicting to $\varphi'(0)=\varphi'(1)=0$. Therefore, $u=e_A$.
\end{proof}

\begin{lemma}\label{lem_quasi_Choquet}
\begin{enumerate}[label=(\arabic*)]
\item For any $A, B\in\mathcal{O}$ with $A\subseteq B$, we have $\mathrm{cap}_1(A)\le \mathrm{cap}_1(B)$.
\item For any $A, B\in\mathcal{O}$, we have
$$\mathrm{cap}_1(A\cup B)+\mathrm{cap}_1(A\cap B)\le\mathrm{cap}_1(A)+\mathrm{cap}_1(B).$$
\item For any $\{A_n\}\subseteq\mathcal{O}$ with $A_n\subseteq A_{n+1}$ for any $n$, we have
$$\mathrm{cap}_1 \left(\cup_nA_n\right)=\sup_n\mathrm{cap}_1(A_n).$$
\end{enumerate}
For any subset $A\subseteq X$, let
\begin{equation}\label{eq_quasi_Choquet}
\mathrm{cap}_1(A)=\inf_{B\in\mathcal{O}:A\subseteq B}\mathrm{cap}_1(B).
\end{equation}
Then $\mathrm{cap}_1$ is a Choquet capacity. Moreover, for any Borel subset $A\subseteq X$, we have
\begin{equation}
\label{eq_quasi_regular}
\mathrm{cap}_1(A)=\sup_{K\in\mathcal{K}:K\subseteq A}\mathrm{cap}_1(K).
\end{equation}
\end{lemma}

\begin{proof}
For (1), it is obvious by definition. For (2), we may assume that $\mathrm{cap}_1(A)$, $\mathrm{cap}_1(B)$ are finite. By Lemma \ref{lem_quasi_open}, there exist unique $e_A\in\mathcal{L}_A$, $e_B\in\mathcal{L}_B$ such that $\mathrm{cap}_1(A)=\mathcal{E}_1(e_A)$, $\mathrm{cap}_1(B)=\mathcal{E}_1(e_B)$. Then $e_A\vee e_B\in\mathcal{L}_{A\cup B}$, $e_A\wedge e_B\in\mathcal{L}_{A\cap B}$, by \ref{eq_SubAdd}, we have
$$\mathrm{cap}_1(A\cup B)+\mathrm{cap}_1(A\cap B)\le\mathcal{E}_1(e_A\vee e_B)+\mathcal{E}_1(e_A\wedge e_B)\le\mathcal{E}_1(e_A)+\mathcal{E}_1( e_B)=\mathrm{cap}_1(A)+\mathrm{cap}_1(B).$$
For (3), by (1), we have ``$\ge$". To show ``$\le$", we may assume that the RHS is finite. By Lemma \ref{lem_quasi_open}, for any $n$, there exists $e_n\in\mathcal{L}_{A_n}$ such that $\mathrm{cap}_1(A_n)=\mathcal{E}_1(e_n)$. Since $\{e_n\}$ is $\mathcal{E}_1$-bounded and $(\mathcal{F},\mathcal{E}_1^{1/p})$ is a reflexive Banach space, by the Banach–Alaoglu theorem (see \cite[Theorem 3 in Chapter 12]{Lax02}), there exists a subsequence, still denoted by $\{e_n\}$, which is $\mathcal{E}_1$-weakly-convergent to some element $e\in\mathcal{F}$. By the Mazur's lemma, here we refer to the version in \cite[Theorem 2 in Section V.1]{Yos95}, for any $n\ge1$, since $\{e_m\}_{m\ge n}$ is $\mathcal{E}_1$-weakly-convergent to $e$, there exist $I_n\ge n$, $\alpha^{(n)}_i\ge0$ for $i=n,\ldots, I_n$ with $\sum_{i=n}^{I_n}\alpha^{(n)}_i=1$ such that $\mathcal{E}_1(\sum_{i=n}^{I_n}\alpha^{(n)}_ie_i-e)<\frac{1}{n}$, let $f_n=\sum_{i=n}^{I_n}\alpha^{(n)}_ie_i$, then $\{f_n\}$ is $\mathcal{E}_1$-convergent to $e$, which implies the $L^p$-convergence. Since $f_n=1$ $m$-a.e. on $A_n$, by passing to a subsequence to have the a.e. convergence, we have $e=1$ $m$-a.e. on $A_n$ for any $n$, which implies $e=1$ $m$-a.e. on $\cup_nA_n$, $e\in\mathcal{L}_{\cup_nA_n}$. Hence
$$\mathrm{cap}_1(\cup_nA_n)\le\mathcal{E}_1(e)\overset{(\star)}{\scalebox{2}[1]{$\le$}}\varliminf_{n\to+\infty}\mathcal{E}_1(e_n)=\sup_n\mathrm{cap}_1(A_n),$$
where $(\star)$ follows from the fact that $\{e_n\}$ is $\mathcal{E}_1$-weakly-convergent to $e$.

By \cite[Theorem A.1.2]{FOT11}, a $[0,+\infty]$-valued set function defined on $\mathcal{O}$ satisfying the above three conditions gives a Choquet capacity by Equation (\ref{eq_quasi_Choquet}). Since $(X,d)$ is a locally compact separable metric space, by \cite[Corollary A.1.1]{FOT11}, we have Equation (\ref{eq_quasi_regular}).
\end{proof}

It is obvious that for any subset $A\subseteq X$, $\mathrm{cap}_1(A)=0$ implies $m(A)=0$. Indeed, for any $n\ge1$, there exists $B_n\in\mathcal{O}$ with $A\subseteq B_n$ such that $\mathrm{cap}_1(B_n)<\frac{1}{n}$, let $e_{B_n}\in\mathcal{L}_{B_n}$ be given by Lemma \ref{lem_quasi_open}, then $m(e_{B_n})\le\mathcal{E}_1(e_{B_n})=\mathrm{cap}_1(B_n)<\frac{1}{n}$, hence $A\subseteq\cap_{n\ge1}B_n$ and $m(\cap_{n\ge1}B_n)=0$, which implies $m(A)$=0.

Let $A\subseteq X$. We say that a statement depending on $x\in A$ holds quasi-everywhere (abbreviated q.e.) on $A$ if there exists $N\subseteq A$ with $\mathrm{cap}_1(N)=0$ such that the statement holds for any $x\in A\backslash N$.

Let $u$ be a $[-\infty,+\infty]$-valued function defined q.e. on $X$. We say that $u$ is quasi-continuous if for any $\varepsilon>0$, there exists an open subset $G\subseteq X$ with $\mathrm{cap}_1(G)<\varepsilon$ such that $u|_{X\backslash G}$ is continuous.

\begin{lemma}\label{lem_quasi_nonnegative}
Let $G\subseteq X$ be an open subset and $u$ quasi-continuous on $G$. If $u\ge0$ $m$-a.e. in $G$, then $u\ge0$ q.e. on $G$.
\end{lemma}

\begin{proof}
For simplicity, we may assume that $G=X$. For any $k\ge1$, there exists a closed set $F_k$ with $\mathrm{cap}_1(X\backslash F_k)<\frac{1}{k}$ such that $u|_{F_k}$ is continuous. Since $\mathcal{L}_{X\backslash F_k}=\mathcal{L}_{X\backslash \mathrm{supp}(1_{F_k}m)}$, we have $\mathrm{cap}_1(X\backslash F_k)=\mathrm{cap}_1(X\backslash \mathrm{supp}(1_{F_k}m))$. By replacing $F_k$ by $\mathrm{supp}(1_{F_k}m)$, we may assume that $F_k=\mathrm{supp}(1_{F_k}m)$. By replacing $F_k$ by $\cup_{i=1}^kF_i$, we may assume that $\{F_k\}$ is increasing, then $\mathrm{cap}_1(X\backslash\cup_kF_k)=0$. We claim that $u(x)\ge0$ for any $x\in\cup_kF_k$. Indeed, suppose that $u(x)<0$ for some $x\in F_k$, since $u|_{F_k}$ is continuous, there exists an open neighborhood $U$ of $x$ such that $u<0$ on $U\cap F_k$. Since $\mathrm{supp}(1_{F_k}m)=F_k$, we have $m(U\cap F_k)>0$, contradicting to $u\ge0$ $m$-a.e. in $X$.
\end{proof}

The first main result of this section is as follows.

\begin{proposition}\label{prop_quasi_exist}
Each $u\in\mathcal{F}$ admits a quasi-continuous modification $\widetilde{u}$, that is, $\widetilde{u}$ is quasi-continuous and $u=\widetilde{u}$ $m$-a.e. in $X$. Moreover, if $\mathcal{C}\subseteq\mathcal{F}\cap C_c(X)$ is $\mathcal{E}_1$-dense in $\mathcal{F}$, then $\widetilde{u}$ can be chosen such that there exist $\{u_n\}\subseteq\mathcal{C}$ which is $\mathcal{E}_1$-convergent to $u$, and an increasing sequence of closed sets $\{F_k\}$ with $\mathrm{cap}_1(X\backslash F_k)\downarrow0$, such that $\{u_n\}$ is uniformly convergent to $\widetilde{u}$ on each $F_k$.
\end{proposition}

\begin{proof}
We follow the same argument as in the proof of \cite[Theorem 2.1.3]{FOT11}. Firstly, for any $u\in\mathcal{F}\cap C(X)$ and any $\lambda>0$, we have $\{x\in X:|u(x)|>\lambda\}$ is an open set, on which $\frac{|u|}{\lambda}\ge1$, hence
\begin{equation}\label{eq_quasi_capbd}
\mathrm{cap}_1 \left(\left\{x\in X:|u(x)|>\lambda\right\}\right)\le\mathcal{E}_1 \left(\frac{|u|}{\lambda}\right)\le \frac{1}{\lambda^p}\mathcal{E}_1(u).
\end{equation}
Secondly, for any $u\in\mathcal{F}$, since $\mathcal{C}\subseteq\mathcal{F}\cap C_c(X)$ is $\mathcal{E}_1$-dense in $\mathcal{F}$, there exists $\{u_n\}\subseteq\mathcal{C}$ which is $\mathcal{E}_1$-convergent to $u$ and satisfies that $\mathcal{E}_1(u_n-u_{n+1})<\frac{1}{2^{(p+1)n}}$ for any $n$. Since $u_n-u_{n+1}\in\mathcal{F}\cap C(X)$, let
$$U_n=\left\{x\in X:|u_n(x)-u_{n+1}(x)|>\frac{1}{2^n}\right\},$$
then $U_n$ is an open set, and by above, we have
$$\mathrm{cap}_1(U_n)\le \frac{1}{\left(\frac{1}{2^n}\right)^p}\mathcal{E}_1(u_n-u_{n+1})\le \frac{1}{2^n}.$$
Let $F_k=X\backslash\cup_{n=k}^{+\infty}U_n$, then $\{F_k\}$ is an increasing sequence of closed sets and
$$\mathrm{cap}_1(X\backslash F_k)=\mathrm{cap}_1(\cup_{n=k}^{+\infty}U_n)\le\sum_{n=k}^{+\infty}\mathrm{cap}_1(U_n)\le\sum_{n=k}^{+\infty}\frac{1}{2^n}=2^{1-k}\to0.$$
For any $k$, we have
$$F_k=\bigcap_{n=k}^{+\infty}\left\{x\in X:|u_n(x)-u_{n+1}(x)|\le \frac{1}{2^n}\right\},$$
hence there exists a real-valued function $v$ defined on $\cup_kF_k$ such that $\{u_n\}$ is uniformly convergent to $v$ on each $F_k$, which implies $\lim_{n\to+\infty}u_n(x)=v(x)$ for any $x\in\cup_kF_k$. Since $u_n|_{F_k}$ is continuous, by the standard $3\varepsilon$-argument, we have $v|_{F_k}$ is continuous for any $k$, which gives $v$ is quasi-continuous. Since $\{u_n\}$ is $\mathcal{E}_1$-convergent to $u$, in particular, also $L^p$-convergent to $u$, by passing to a subsequence to have the a.e. convergence, we have $u=v$ $m$-a.e. in $X$, that is, $v$ is a modification of $u$.
\end{proof}

For any $u\in\mathcal{F}$, we always use $\widetilde{u}$ to denote a quasi-continuous modification of $u$. Let
$$\widetilde{\mathcal{F}}=\left\{u\in\mathcal{F}:u\text{ is quasi-continuous}\right\},$$
then by Lemma \ref{lem_quasi_nonnegative} and Proposition \ref{prop_quasi_exist}, the equivalence classes of $\widetilde{\mathcal{F}}$ with respect to q.e. on $X$ are equal to the equivalence classes of $\mathcal{F}$ with respect to $m$-a.e. in $X$. We have a generalization of Equation (\ref{eq_quasi_capbd}) as follows.

\begin{lemma}\label{lem_quasi_capbd}
For any $u\in\widetilde{\mathcal{F}}$, for any $\lambda>0$, we have
$$\mathrm{cap}_1\left(\left\{x\in X:|u(x)|>\lambda\right\}\right)\le\frac{1}{\lambda^p}\mathcal{E}_1(u).$$
\end{lemma}

\begin{proof}
By Proposition \ref{prop_quasi_exist}, there exists $\{u_n\}\subseteq\mathcal{F}\cap C_c(X)$ satisfying that $\{u_n\}$ is $\mathcal{E}_1$-convergent to $u$, and for any $\varepsilon>0$, there exists an open set $G$ with $\mathrm{cap}_1(G)<\varepsilon$ such that $\left\{u_n\right\}$ converges uniformly to $u$ on $X\backslash G$. For any $\lambda_1\in(0,\lambda)$, there exists $N\ge1$ such that for any $n>N$, we have
$$\left\{x\in X:|u(x)|>\lambda\right\}\subseteq\left\{x\in X\backslash G:|u_n(x)|>\lambda_1\right\}\cup G\subseteq\left\{x\in X:|u_n(x)|>\lambda_1\right\}\cup G.$$
For any $n$, by Equation (\ref{eq_quasi_capbd}), we have
$$\mathrm{cap}_1\left(\left\{x\in X:|u_n(x)|>\lambda_1\right\}\right)\le\frac{1}{\lambda_1^p}\mathcal{E}_1(u_n),$$
hence
$$\mathrm{cap}_1\left(\left\{x\in X:|u(x)|>\lambda\right\}\right)\le\mathrm{cap}_1\left(\left\{x\in X:|u_n(x)|>\lambda_1\right\}\right)+\mathrm{cap}_1(G)\le\frac{1}{\lambda_1^p}\mathcal{E}_1(u_n)+\varepsilon.$$
Firstly, letting $n\to+\infty$, secondly, letting $\lambda_1\uparrow\lambda$, finally, letting $\varepsilon\downarrow0$, we have
$$\mathrm{cap}_1\left(\left\{x\in X:|u(x)|>\lambda\right\}\right)\le\frac{1}{\lambda^p}\mathcal{E}_1(u).$$
\end{proof}

A direct consequence of the above result is as follows.

\begin{corollary}\label{cor_quasi_Cauchy}
Let $\{u_n\}\subseteq\widetilde{\mathcal{F}}$ be an $\mathcal{E}_1$-Cauchy sequence. Then there exist a subsequence $\{u_{n_k}\}$ and a function $u\in\widetilde{\mathcal{F}}$ such that $\{{u}_{n_k}\}$ converges to ${u}$ q.e. on $X$ and $\left\{u_n\right\}$ is $\mathcal{E}_1$-convergent to $u$.
\end{corollary}

\begin{proof}
Take a subsequence $\{u_{n_k}\}$ satisfying $\mathcal{E}_1(u_{n_k}-u_{n_{k+1}})<\frac{1}{2^{(p+1)k}}$, then  by Lemma \ref{lem_quasi_capbd}, we have
$$\mathrm{cap}_1\left(\left\{x\in X:|u_{n_k}(x)-u_{n_{k+1}}(x)|>\frac{1}{2^k}\right\}\right)\le\frac{1}{\left(\frac{1}{2^k}\right)^p}\mathcal{E}_1(u_{n_k}-u_{n_{k+1}})\le\frac{1}{2^k}.$$
Let
$$G_k=\left\{x\in X:|u_{n_k}(x)-u_{n_{k+1}}(x)|>\frac{1}{2^k}\right\},$$
then for any $n$, we have
$$\mathrm{cap}_1\left(\bigcup_{k\ge n}G_k\right)\le\sum_{k\ge n}\mathrm{cap}_1\left(G_k\right)\le\sum_{k\ge n}\frac{1}{2^k}=2^{1-n}\to0,$$
and
$$\bigcap_{k\ge n}G_k^c=\left\{x\in X:|u_{n_k}(x)-u_{n_{k+1}}(x)|\le\frac{1}{2^k}\right\},$$
hence there exists a real-valued function $u$ defined on $\bigcup_{n\ge1}\bigcap_{k\ge n}G_k^c$ such that $\{u_{n_k}\}$ converges uniformly to $u$ on $\bigcap_{k\ge n}G_k^c$ for any $n$, which implies $\lim_{k\to+\infty}u_{n_k}(x)=u(x)$ for any $x\in\bigcup_{n\ge1}\bigcap_{k\ge n}G_k^c$, $\left\{u_{n_k}\right\}$ converges to $u$ q.e. on $X$, $\left\{u_{n_k}\right\}$ converges to $u$ $m$-a.e. in $X$. Since $\left\{u_n\right\}$ is an $\mathcal{E}_1$-Cauchy sequence, there exists $v\in\mathcal{F}$ such that $\left\{u_n\right\}$ is $\mathcal{E}_1$-convergent to $v$, there exists a subsequence of $\left\{u_{n_k}\right\}$ that converges to $v$ $m$-a.e. in $X$. Hence $u=v$ $m$-a.e. in $X$, $u\in\mathcal{F}$ and $\{u_n\}$ is $\mathcal{E}_1$-convergent to $u$.

It remains to show that $u$ is quasi-continuous. Indeed, since $\{u_n\}$ is a sequence of quasi-continuous functions, there exists an increasing sequence of closed sets $\left\{F_k\right\}$ with $\mathrm{cap}_1(F_k)\downarrow0$ such that $u_n|_{F_k}$ is continuous for any $n, k$. For any $\varepsilon>0$, take $n, l\ge1$ such that $2^{1-n}<\varepsilon/2$ and $\mathrm{cap}_1(X\backslash F_l)<\varepsilon/2$. Then $\mathrm{cap}_1\left(\cup_{k\ge n}G_k\right)\le{2^{1-n}}<{\varepsilon}/{2}$, there exists an open set $U\supseteq\cup_{k\ge n}G_k$ such that $\mathrm{cap}_1(U)<\varepsilon/2$, then $F_l\cap U^c$ is a closed set and
$$\mathrm{cap}_1(X\backslash(F_l\cap U^c))=\mathrm{cap}_1((X\backslash F_l)\cup U)\le\mathrm{cap}_1(X\backslash F_l)+\mathrm{cap}_1(U)<\frac{\varepsilon}{2}+\frac{\varepsilon}{2}=\varepsilon.$$
Since $u_{n_k}|_{F_l}$ is continuous, we have $u_{n_k}|_{F_l\cap U^c}$ is continuous. Since $\left\{u_{n_k}\right\}$ converges uniformly to $u$ on $\cap_{k\ge n}G_k^c\supseteq U^c$, we have $\left\{u_{n_k}\right\}$ converges uniformly to $u$ on $F_l\cap U^c$. Hence $u|_{F_l\cap U^c}$ is continuous. Therefore, $u$ is quasi-continuous.
\end{proof}

For any subset $A\subseteq X$, let
$$\mathcal{L}_A=\left\{u\in\mathcal{F}:\widetilde{u}\ge1\text{ q.e. on }A\right\}.$$
By Lemma \ref{lem_quasi_nonnegative}, if $A$ is open, then the above definition coincides with the one given by Equation (\ref{eq_quasi_LAopen}). Similar to Lemma \ref{lem_quasi_open}, we have the following result.

\begin{lemma}\label{lem_quasi_subset}
For any subset $A$ with $\mathcal{L}_A\ne\emptyset$, there exists a unique element $e_A\in\mathcal{L}_A$ such that $\mathrm{cap}_1(A)=\mathcal{E}_1(e_A)$. Moreover, we have $0\le e_A\le 1$ $m$-a.e. in $X$, $\widetilde{e}_A=1$ q.e. on $A$, $e_A$ is the unique element $u\in\mathcal{F}$ satisfying $\widetilde{u}=1$ q.e. on $A$ and $\mathcal{E}_1(u;v)\ge0$ for any $v\in\mathcal{F}$ with $\widetilde{v}\ge0$ q.e. on $A$.
\end{lemma}

\begin{proof}
It is obvious that $\mathcal{L}_A$ is a non-empty convex subset of $\mathcal{F}$. We claim that $\mathcal{L}_A$ is closed in $(\mathcal{F},\mathcal{E}_1^{1/p})$. Indeed, for any $\mathcal{E}_1$-Cauchy sequence $\{u_n\}\subseteq\mathcal{L}_A$, by Corollary \ref{cor_quasi_Cauchy}, there exist a subsequence $\{u_{n_k}\}$ and a function $u\in\widetilde{\mathcal{F}}$ such that $\{\widetilde{u}_{n_k}\}$ converges to $u$ q.e. on $X$ and $\{u_n\}$ is $\mathcal{E}_1$-convergent to $u$. Since $\widetilde{u}_{n_k}\ge1$ q.e. on $A$, we have $u\ge1$ q.e. on $A$, which gives $u\in\mathcal{L}_A$, $\mathcal{L}_A$ is closed. By the same argument as in the proof of Lemma \ref{lem_quasi_open}, there exists a unique element $e_A\in\mathcal{L}_A$ such that $\mathcal{E}_1(e_A)=\inf_{u\in\mathcal{L}_A}\mathcal{E}_1(u)$. Moreover, $0\le e_A\le 1$ $m$-a.e. in $X$, $\widetilde{e}_A=1$ q.e. on $A$, $e_A$ is the unique element $u\in\mathcal{F}$ satisfying $\widetilde{u}=1$ q.e. on $A$ and $\mathcal{E}_1(u;v)\ge0$ for any $v\in\mathcal{F}$ with $\widetilde{v}\ge0$ q.e. on $A$.

It remains to show that $\mathcal{E}_1(e_A)=\mathrm{cap}_1(A)$, or equivalently,
$$\mathcal{E}_1(e_A)=\inf_{u\in\mathcal{L}_A}\mathcal{E}_1(u)=\inf_{B\in\mathcal{O}:A\subseteq B}\mathrm{cap}_1(B)=\mathrm{cap}_1(A).$$
``$\le$": For any $B\in\mathcal{O}$ with $A\subseteq B$, be Lemma \ref{lem_quasi_open}, there exists $e_B\in\mathcal{L}_B$ with $e_B=1$ $m$-a.e. on $B$ such that $\mathrm{cap}_1(B)=\mathcal{E}_1(e_B)$. By Lemma \ref{lem_quasi_nonnegative}, $\widetilde{e}_B=1$ q.e. on $B\supseteq A$, hence $e_B\in\mathcal{L}_A$, which gives $\inf_{u\in\mathcal{L}_A}\mathcal{E}_1(u)\le\mathcal{E}_1(e_B)=\mathrm{cap}_1(B)$. Taking the infimum over all such $B$, we have the desired result.

``$\ge$": Since $\widetilde{e}_A$ is quasi-continuous and $\widetilde{e}_A=1$ q.e. on $A$, for any $\varepsilon>0$, there exists an open set $G$ with $\mathrm{cap}_1(G)<\varepsilon$ such that $\widetilde{e}_A|_{X\backslash G}$ is continuous and $\widetilde{e}_A=1$ on $A\backslash G$. For any $\delta\in(0,1)$, there exists an open set $U$ with $U\backslash G\supseteq A\backslash G$ such that $\widetilde{e}_A>1-\delta$ on $U\backslash G$, then $U\cup G$ is an open set satisfying $U\cup G\supseteq A$ and $U\cup G\subseteq\left\{x\in X:\widetilde{e}_A(x)>1-\delta\right\}\cup G$. By Lemma \ref{lem_quasi_capbd}, we have
$$\mathrm{cap}_1(U\cup G)\le\mathrm{cap}_1(\left\{x\in X:\widetilde{e}_A(x)>1-\delta\right\})+\mathrm{cap}_1(G)\le\frac{1}{(1-\delta)^p}\mathcal{E}_1(e_A)+\varepsilon,$$
hence
$$\mathrm{cap}_1(A)\le\mathrm{cap}_1(U\cup G)\le\frac{1}{(1-\delta)^2}\mathcal{E}_1(e_A)+\varepsilon.$$
Letting $\varepsilon\downarrow0$ and $\delta\downarrow0$, we have the desired result.
\end{proof}

Following \cite[Section 2.2]{FOT11}, we consider measures of finite energy integral as follows.

\begin{proposition}\label{prop_quasi_eneint}
Let $u\in\mathcal{F}$. The followings are equivalent.
\begin{enumerate}[label=(\arabic*),ref=(\arabic*)]
\item\label{cond_quasi_Radon} There exists a positive Radon measure $\mu$ on $X$ such that
$$\mathcal{E}_1(u;v)=\int_Xv\mathrm{d}\mu\text{ for any }v\in\mathcal{F}\cap C_c(X).$$
\item\label{cond_quasi_posFX} $\mathcal{E}_1(u;v)\ge0$ for any $v\in\mathcal{F}$ with $v\ge0$ $m$-a.e. in $X$.
\item\label{cond_quasi_posFCX} $\mathcal{E}_1(u;v)\ge0$ for any $v\in\mathcal{F}\cap C_c(X)$ with $v\ge0$ on $X$.
%\item $\mathcal{E}_1(u;v)\ge0$ for any $v\in\mathcal{C}$ with $v\ge0$ on $X$, for any (some) $\mathcal{C}\subseteq\mathcal{F}\cap C_c(X)$, which is uniformly dense in $C_c(X)$ and $\mathcal{E}_1$-dense in $\mathcal{F}$.
\end{enumerate}
Moreover, if the above conditions hold, then for any Borel subset $A\subseteq X$, we have
\begin{equation}\label{eq_quasi_mucap}
\mu(A)\le\mathcal{E}_1(u)^{(p-1)/p}\mathrm{cap}_1(A)^{1/p}.
\end{equation}
In particular, $\mu$ charges no set of zero capacity, that is, for any subset $A\subseteq X$, if $\mathrm{cap}_1(A)=0$, then $\mu(A)=0$. Hence $\widetilde{\mathcal{F}}\subseteq L^1(X;\mu)$ and
\begin{equation}\label{eq_quasi_L1mu}
\mathcal{E}_1(u;v)=\int_X\widetilde{v}\mathrm{d}\mu\text{ for any }v\in\mathcal{F}.
\end{equation}
\end{proposition}

\begin{proof}
``\ref{cond_quasi_Radon}$\Rightarrow$\ref{cond_quasi_posFCX}" and ``\ref{cond_quasi_posFX}$\Rightarrow$\ref{cond_quasi_posFCX}": Obvious.

``\ref{cond_quasi_posFCX}$\Rightarrow$\ref{cond_quasi_posFX}": For any $v\in\mathcal{F}$ with $v\ge0$ $m$-a.e. in $X$, there exists $\{v_n\}\subseteq\mathcal{F}\cap C_c(X)$ which is $\mathcal{E}_1$-convergent to $v$. By \cite[Corollary 3.19 (b)]{KS24a}, we have $\left\{v_n^+\right\}$ is $\mathcal{E}_1$-convergent to $v^+=v$. Since $v_n^+\in\mathcal{F}\cap C_c(X)$ satisfies $v_n^+\ge0$ on $X$, by assumption, we have $\mathcal{E}_1(u;v_n^+)\ge0$, which gives $\mathcal{E}_1(u;v)=\lim_{n\to+\infty}\mathcal{E}_1(u;v_n^+)\ge0$.

``\ref{cond_quasi_posFCX}$\Rightarrow$\ref{cond_quasi_Radon}": The idea is to construct a positive linear functional on $C_c(X)$ to apply the Riesz representation theorem. For any $v\in C_c(X)$, there exist $\left\{v_n\right\}\subseteq\mathcal{F}\cap C_c(X)$ that converges uniformly to $v$, and $\varphi\in\mathcal{F}\cap C_c(X)$ with $0\le\varphi\le1$ on $X$ and $\varphi=1$ on $\mathrm{supp}(v)$, then $\left\{\varphi v_n\right\}\subseteq\mathcal{F}\cap C_c(X)$ converges uniformly to $\varphi v=v$ and $\mathrm{supp}(\varphi v_n)\subseteq\mathrm{supp}(\varphi)$ for any $n$. For any $n, m$, we have $\lVert v_n-v_m\rVert_{L^\infty(X;m)}\varphi\pm\varphi(v_n-v_m)\in\mathcal{F}\cap C_c(X)$ and $\lVert v_n-v_m\rVert_{L^\infty(X;m)}\varphi\pm\varphi(v_n-v_m)\ge0$ on $X$. By assumption, we have
$$\mathcal{E}_1(u;\lVert v_n-v_m\rVert_{L^\infty(X;m)}\varphi\pm\varphi(v_n-v_m))\ge0,$$
hence
$$|\mathcal{E}_1(u;\varphi v_n)-\mathcal{E}_1(u;\varphi v_m)|\le \lVert v_n-v_m\rVert_{L^\infty(X;m)}\mathcal{E}_1(u;\varphi)\to0$$
as $n,m\to+\infty$, that is, $\{\mathcal{E}_1(u;\varphi v_n)\}$ is a Cauchy sequence. Define
$$l(v)=\lim_{n\to+\infty}\mathcal{E}_1(u;\varphi v_n).$$
It is obvious that $l(v)$ is well-defined, that is, the limit is independent of the choice of $\left\{v_n\right\}$ and $\varphi$, and that $l$ is a linear functional on $C_c(X)$. Moreover, $l(v)=\mathcal{E}_1(u;v)$ for any $v\in\mathcal{F}\cap C_c(X)$.

We claim that $l(v)\ge0$ if $v\in C_c(X)$ satisfies $v\ge0$ on $X$. Since $\{v_n\}\subseteq\mathcal{F}\cap C_c(X)$ converges uniformly to $v$ and $v\ge0$ on $X$, we have $\{v_n^+\}\subseteq\mathcal{F}\cap C_c(X)$ also converges uniformly to $v$, for any $\varepsilon>0$, there exists $N$ such that $\varphi v_n^+\ge-\varepsilon\varphi$ on $X$ for any $n>N$, then $\mathcal{E}_1(u;\varphi v_n^+)\ge-\varepsilon\mathcal{E}_1(u;\varphi)$. Hence
$$l(v)=\lim_{n\to+\infty}\mathcal{E}_1(u;\varphi v_n^+)\ge-\varepsilon\mathcal{E}_1(u;\varphi).$$
Letting $\varepsilon\downarrow0$, we have $l(v)\ge0$. Therefore, $l$ is a positive linear functional on $C_c(X)$. By the Riesz representation theorem, there exists a positive Radon measure $\mu$ on $X$ such that
$$l(v)=\int_Xv\mathrm{d}\mu\text{ for any }v\in C_c(X).$$
In particular, for any $v\in\mathcal{F}\cap C_c(X)$, we have
$$\mathcal{E}_1(u;v)=l(v)=\int_Xv\mathrm{d}\mu.$$

To prove Equation (\ref{eq_quasi_mucap}) for any Borel set $A$, by the regular property of $\mu$ and $\mathrm{cap}_1$, we only need to prove for any compact set $K$. Indeed, since $K$ is compact, we have $\mathcal{L}_K\ne\emptyset$, by Lemma \ref{lem_quasi_subset}, we have $\mathrm{cap}_1(K)<+\infty$, which gives
$$\mathrm{cap}_1(K)=\inf_{U\in \mathcal{O}:K\subseteq U}\mathrm{cap}_1(U)=\inf_{U\in \mathcal{O}:\mathcal{L}_U\ne\emptyset,K\subseteq U}\mathrm{cap}_1(U).$$
For any $U\in \mathcal{O}$ satisfying $\mathcal{L}_U\ne\emptyset$ and $K\subseteq U$, by Lemma \ref{lem_quasi_open}, there exists $e_U\in \mathcal{F}$ satisfying the conditions therein. Take $\{u_n\}\subseteq\mathcal{F}\cap C_c(X)$ which is $\mathcal{E}_1$-convergent to $e_U$. Since $0\le e_U\le 1$ $m$-a.e. in $X$, by replacing $u_n$ by $(u_n\vee0)\wedge1$ and applying \cite[Corollary 3.19 (b)]{KS24a}, we may assume that $0\le u_n\le 1$ on $X$. Moreover, there exists $\varphi\in\mathcal{F}\cap C_c(X)$ such that $0\le\varphi\le1$ on $X$, $\varphi=1$ on $K$, and $\mathrm{supp}(\varphi)\subseteq U$. Let $v_n=\varphi+u_n-\varphi u_n$, then $v_n\in\mathcal{F}\cap C_c(X)$, $v_n\ge0$ on $X$, and $v_n=1+u_n-u_n=1$ on $K$. By \ref{eq_Alg}, we have
\begin{align*}
&{\mathcal{E}(\varphi u_n)}^{1/p}\le C_p \max\{\lVert \varphi\rVert_{L^\infty(X;m)},\lVert u_n\rVert_{L^\infty(X;m)}\} \left(\mathcal{E}(\varphi)^{1/p}+\mathcal{E}(u_n)^{1/p}\right)\\
&\le C_p \left(\mathcal{E}(\varphi)^{1/p}+\sup_n\mathcal{E}(u_n)^{1/p}\right)<+\infty,
\end{align*}
hence $\sup_{n}\mathcal{E}(\varphi u_n)<+\infty$. Since $\{\varphi u_n\}$ is $L^p(X;m)$-convergent to $\varphi e_U$, by \cite[Lemma 3.17]{KS24a}, we have $\{\varphi u_n\}$ is $\mathcal{E}_1$-weakly-convergent to $\varphi e_U$. By the Mazur's lemma (see \cite[Theorem 2 in Section V.1]{Yos95}), for any $n$, there exist $I_n\ge n$, $\lambda^{(n)}_k\ge0$ for $k=n,\ldots,I_n$ with $\sum_{k=n}^{I_n}\lambda_k^{(n)}=1$, such that $\{\sum_{k=n}^{I_n}\lambda^{(n)}_k\varphi u_k\}_n$ is $\mathcal{E}_1$-convergent to $\varphi e_U$. Let $w_n=\sum_{k=n}^{I_n}\lambda^{(n)}_kv_k$, then $w_n\in\mathcal{F}\cap C_c(X)$, $w_n\ge0$ on $X$, $w_n=1$ on $K$, and $\{w_n\}$ is $\mathcal{E}_1$-convergent to $\varphi+e_U-\varphi e_U=e_U+\varphi(1-e_U)=e_U+\varphi(1-1)=e_U$, by noting that $e_U=1$ $m$-a.e. in $U\supseteq \mathrm{supp}(\varphi)$. Hence
\begin{align*}
&\mu(K)=\int_K w_n \mathrm{d}\mu\le \int_X w_n \mathrm{d}\mu\overset{(\star)}{\scalebox{2}[1]{$=$}}\mathcal{E}_1(u;w_n)\le \mathcal{E}_1(u)^{(p-1)/p}\mathcal{E}_1(w_n)^{1/p}\\
&\to \mathcal{E}_1(u)^{(p-1)/p}\mathcal{E}_1(e_U)^{1/p}=\mathcal{E}_1(u)^{(p-1)/p}\mathrm{cap}_1(U)^{1/p},
\end{align*}
where $(\star)$ follows from \ref{cond_quasi_Radon} for $w_n$. Taking the infimum over all such $U$, we have $\mu(K)\le \mathcal{E}_1(u)^{(p-1)/p}\mathrm{cap}_1(K)^{1/p}$.

In particular, for any subset $A\subseteq X$, if $\mathrm{cap}_1(A)=0$, then for any $n\ge1$, there exists an open set $G_n\supseteq A$ such that $\mathrm{cap}_1(G_n)<\frac{1}{n}$, then $A\subseteq\cap_{n\ge1}G_n$ and
$$\mu(\cap_{n\ge1}G_n)\le\mu(G_n)\le\mathcal{E}_1(u)^{(p-1)/p}\mathrm{cap}_1(G_n)^{1/p}\le\mathcal{E}_1(u)^{(p-1)/p}\left(\frac{1}{n}\right)^{1/p}\to0,$$
hence $\mu(\cap_{n\ge1}G_n)=0$, which gives $\mu(A)=0$.

Finally, for any $v\in\mathcal{F}$, by Proposition \ref{prop_quasi_exist}, there exists $\{v_n\}\subseteq\mathcal{F}\cap C_c(X)$ such that $\{v_n\}$ is $\mathcal{E}_1$-convergent to $v$ and $\{v_n\}$ converges to $\widetilde{v}$ q.e. on $X$. By above, the q.e. convergence implies the $\mu$-a.e. convergence, hence by Fatou's lemma, for any $n$, we have
\begin{align*}
&\int_X|v_n-\widetilde{v}|\mathrm{d}\mu=\int_X\lim_{m\to+\infty}|v_n-v_m|\mathrm{d}\mu\le\varliminf_{m\to+\infty}\int_X|v_n-v_m|\mathrm{d}\mu=\varliminf_{m\to+\infty}\mathcal{E}_1(u;|v_n-v_m|)\\
&\le\varliminf_{m\to+\infty}\mathcal{E}_1(u)^{(p-1)/p}\mathcal{E}_1(|v_n-v_m|)^{1/p}\le\varliminf_{m\to+\infty}\mathcal{E}_1(u)^{(p-1)/p}\mathcal{E}_1(v_n-v_m)^{1/p},
\end{align*}
which implies that $\widetilde{v}\in L^1(X;\mu)$ and $\{v_n\}$ is $L^1(X;\mu)$-convergent to $\widetilde{v}$, hence
$$\int_X\widetilde{v}\mathrm{d}\mu=\lim_{n\to+\infty}\int_Xv_n\mathrm{d}\mu=\lim_{n\to+\infty}\mathcal{E}_1(u;v_n)=\mathcal{E}_1(u;v).$$
\end{proof}

\begin{lemma}\label{lem_quasi_muonF}
Let $u\in\mathcal{F}$ and let $F\subseteq X$ be a closed subset. The followings are equivalent.
\begin{enumerate}[label=(\roman*),ref=(\roman*)]
\item\label{cond_quasi_RadonF} There exists a positive Radon measure $\mu$ on $X$ with $\mathrm{supp}(\mu)\subseteq F$ such that
$$\mathcal{E}_1(u;v)=\int_Xv\mathrm{d}\mu\text{ for any }v\in\mathcal{F}\cap C_c(X).$$
\item\label{cond_quasi_posFF} $\mathcal{E}_1(u;v)\ge0$ for any $v\in\mathcal{F}$ with $\widetilde{v}\ge0$ q.e. on $F$.
\item\label{cond_quasi_posFCF} $\mathcal{E}_1(u;v)\ge0$ for any $v\in\mathcal{F}\cap C_c(X)$ with $v\ge0$ on $F$.
%\item $\mathcal{E}_1(u;v)\ge0$ for any $v\in\mathcal{C}$ with $v\ge0$ on $F$, for any (some) $\mathcal{C}\subseteq\mathcal{F}\cap C_c(X)$, which is uniformly dense in $C_c(X)$ and $\mathcal{E}_1$-dense in $\mathcal{F}$.
\end{enumerate}
\end{lemma}

\begin{proof}
``\ref{cond_quasi_RadonF}$\Rightarrow$\ref{cond_quasi_posFCF}" and ``\ref{cond_quasi_posFF}$\Rightarrow$\ref{cond_quasi_posFCF}": Obvious.

``\ref{cond_quasi_posFCF}$\Rightarrow$\ref{cond_quasi_RadonF}": It is obvious that \ref{cond_quasi_posFCX} in Proposition \ref{prop_quasi_eneint} holds. By Proposition \ref{prop_quasi_eneint}, there exists a positive Radon measure $\mu$ on $X$ such that $\mathcal{E}(u;v)=\int_Xv\mathrm{d}\mu$ for any $v\in\mathcal{F}\cap C_c(X)$. It remains to show that $\mathrm{supp}(\mu)\subseteq F$. Suppose that $\mathrm{supp}(\mu)\not\subseteq F$. Since $F$ is closed, there exist a compact set $K$  and an open set $G$ satisfying $K\subseteq G\subseteq X\backslash F$ and $\mu(K)>0$. There exists $v\in\mathcal{F}\cap C_c(X)$ with $0\le v\le1$ on $X$, $v=1$ on $K$ and $\mathrm{supp}(v)\subseteq G$. Then $v=0$ on $F$, by assumption, we have $\mathcal{E}_1(u;v)=0$, which implies
$$0=\mathcal{E}_1(u;v)=\int_Xv\mathrm{d}\mu\ge\mu(K)>0,$$
contradiction. Hence $\mathrm{supp}(\mu)\subseteq F$.

``\ref{cond_quasi_RadonF}$\Rightarrow$\ref{cond_quasi_posFF}": It is obvious that \ref{cond_quasi_Radon} in Proposition \ref{prop_quasi_eneint} holds. For any $v\in\mathcal{F}$ with $\widetilde{v}\ge0$ q.e. on $F$, by Proposition \ref{prop_quasi_eneint}, we have $\widetilde{v}\in L^1(X;\mu)$ and $\widetilde{v}\ge0$ $\mu$-a.e. on $F$. By assumption, applying Equation (\ref{eq_quasi_L1mu}), we have
$$\mathcal{E}_1(u;v)=\int_X\widetilde{v}\mathrm{d}\mu=\int_F\widetilde{v}\mathrm{d}\mu\ge0.$$
\end{proof}

We have the following characterization of capacity for compact sets.

\begin{lemma}\label{lem_quasi_compact}
For any compact set $K$, we have
$$\mathrm{cap}_1(K)=\inf_{u\in\mathcal{C}^K}\mathcal{E}_1(u),$$
where
$$\mathcal{C}^K=\left\{u\in\mathcal{F}\cap C_c(X):u\ge1\text{ on }K\right\}.$$
Moreover, there exists $\{u_n\}\subseteq\mathcal{F}\cap C_c(X)$ with $0\le u_n\le1$ on $X$, $u_n=1$ on $K$ such that $\lim_{n\to+\infty}\mathcal{E}_1(u_n)=\mathrm{cap}_1(K)$.
\end{lemma}

\begin{proof}
Let $\overline{\mathcal{C}}^K$ be the $\mathcal{E}_1$-closure of $\mathcal{C}^K$ in $\mathcal{F}$, then $\overline{\mathcal{C}}^K$ is a non-empty closed convex subset of $\mathcal{F}$, by \cite[Theorem 8 in Chapter 5]{Lax02} and its proof, there exists a unique element $u\in\overline{\mathcal{C}}^K$ such that $\mathcal{E}_1(u)=\inf_{u\in\overline{\mathcal{C}}^K}\mathcal{E}_1(u)=\inf_{u\in\mathcal{C}^K}\mathcal{E}_1(u)$, and any minimizing sequence in $\mathcal{C}^K$ is an $\mathcal{E}_1$-Cauchy sequence which is $\mathcal{E}_1$-convergent to $u$.

To prove that $\mathcal{E}_1(u)=\mathrm{cap}_1(K)$, by Lemma \ref{lem_quasi_subset}, we only need to prove that $\widetilde{u}=1$ q.e. on $K$ and $\mathcal{E}_1(u;v)\ge0$ for any $v\in\mathcal{F}$ with $\widetilde{v}\ge0$ q.e. on $K$.

Firstly, we show that $\widetilde{u}=1$ q.e. on $K$. Let $\phi:\mathbb{R}\to\mathbb{R}$, $t\mapsto(t\vee0)\wedge1$. For any minimizing sequence $\{u_n\}$ in $\mathcal{C}^K$, we have $\{\phi(u_n)\}\subseteq\mathcal{C}^K$ and $\mathcal{E}_1(u)\le\mathcal{E}_1(\phi(u_n))\le\mathcal{E}_1(u_n)\to\mathcal{E}_1(u)$, hence $\{\phi(u_n)\}$ is also a minimizing sequence in $\mathcal{C}^K$, which is an $\mathcal{E}_1$-Cauchy sequence that $\mathcal{E}_1$-converges to $u$. Moreover, $0\le\phi(u_n)\le1$ on $X$ and $\phi(u_n)=1$ on $K$. By Corollary \ref{cor_quasi_Cauchy}, there exists a subsequence $\left\{\phi(u_{n_k})\right\}$ that converges to $\widetilde{u}$ q.e. on $X$. Since $\phi(u_{n_k})=1$ on $K$, we have $\widetilde{u}=1$ q.e. on $K$.

Secondly, we show that $\mathcal{E}_1(u;v)\ge0$ for any ${v}\in\mathcal{F}$ with $\widetilde{v}\ge0$ q.e. on $K$. By Lemma \ref{lem_quasi_muonF}, we only need to show that $\mathcal{E}_1(u;v)\ge0$ for any $v\in\mathcal{F}\cap C_c(X)$ with $v\ge0$ on $K$. For any minimizing sequence $\{u_n\}$ in $\mathcal{C}^K$, we have $\{u_n\}$ is $\mathcal{E}_1$-convergent to $u$, for any $\varepsilon>0$, we have $\{u_n+\varepsilon v\}\subseteq\mathcal{C}^K$, hence $\mathcal{E}_1(u)=\inf_{u\in\mathcal{C}^K}\mathcal{E}_1(u)\le\mathcal{E}_1(u_n+\varepsilon v)$. Letting $n\to+\infty$, we have $\mathcal{E}_1(u)\le\mathcal{E}_1(u+\varepsilon v)$, hence $\mathcal{E}_1(u;v)=\frac{1}{p}\lim_{\varepsilon\downarrow0}\frac{1}{\varepsilon}\left(\mathcal{E}_1(u+\varepsilon v)-\mathcal{E}_1(u)\right)\ge0$.
\end{proof}

The second main result of this section is as follows.

\begin{proposition}\label{prop_quasi_charge}
For any $u\in\mathcal{F}$, we have $\Gamma(u)$ charges no set of zero capacity, that is, for any subset $A\subseteq X$, $\mathrm{cap}_1(A)=0$ implies $\Gamma(u)(A)=0$.
\end{proposition}

\begin{proof}
Since $\mathcal{F}\cap C_c(X)$ is $\mathcal{E}_1$-dense in $\mathcal{F}$, we only need to show that for any $u\in\mathcal{F}\cap C_c(X)$, for any compact subset $K\subseteq X$ with $\mathrm{cap}_1(K)=0$, we have $\Gamma(u)(K)=0$. Indeed, by Lemma \ref{lem_quasi_compact}, there exists $\{\phi_n\}\subseteq\mathcal{F}\cap C_c(X)$ with $0\le\phi_n\le1$ on $X$, $\phi_n=1$ on $K$ such that $\lim_{n\to+\infty}\mathcal{E}_1(\phi_n)=\mathrm{cap}_1(K)=0$. By \cite[Proposition 4.16]{KS24a}, we have
$$0\le\Gamma(u)(K)\le\int_X\phi_n\mathrm{d}\Gamma(u)=\mathcal{E}(u;u\phi_n)-\left(\frac{p-1}{p}\right)^{p-1}\mathcal{E} \left(|u|^{\frac{p}{p-1}};\phi_n\right),$$
where
$$\lvert\mathcal{E}\left(|u|^{\frac{p}{p-1}};\phi_n\right)\rvert\le\mathcal{E} \left(|u|^{\frac{p}{p-1}}\right)^{{(p-1)}/{p}}\mathcal{E}(\phi_n)^{{1}/{p}}\to0.$$
Since $\{\phi_n\}$ is $L^p(X;m)$-convergent to $0$ and $u\in\mathcal{F}\cap C_c(X)$, we have $\{u\phi_n\}$ is $L^p(X;m)$-convergent to $0$. Since
\begin{align*}
&\mathcal{E}(u\phi_n)^{{1}/{p}}\le C_p\max\{\lVert u\rVert_{L^\infty(X;m)},\lVert \phi_n\rVert_{L^\infty(X;m)}\}\left(\mathcal{E}(\phi_n)^{{1}/{p}}+\mathcal{E}(u)^{{1}/{p}}\right)\\
&\le C_p\max\{\lVert u\rVert_{L^\infty(X;m)},1\}\left(\sup_n\mathcal{E}(\phi_n)^{{1}/{p}}+\mathcal{E}(u)^{{1}/{p}}\right)<+\infty,
\end{align*}
we have $\sup_{n}\mathcal{E}(u\phi_n)<+\infty$. By \cite[Lemma 3.17]{KS24a}, we have $\{u\phi_n\}$ is $\mathcal{E}_1$-weakly-convergent to 0, which gives $\lim_{n\to+\infty}\mathcal{E}(u;u\phi_n)=0$. Hence $0\le\Gamma(u)(K)\le\lim_{n\to+\infty}\int_X\phi_n\mathrm{d}\Gamma(u)=0$, that is, $\Gamma(u)(K)=0$.
\end{proof}

\bibliographystyle{plain}
%\bibliography{/Users/meng/Dropbox/myref}

\begin{thebibliography}{10}

\bibitem{AB25}
Patricia Alonso-Ruiz and Fabrice Baudoin.
\newblock Korevaar-{S}choen {$p$}-energies and their {$\Gamma $}-limits on
  {C}heeger spaces.
\newblock {\em Nonlinear Anal.}, 256:Paper No. 113779, 22, 2025.

\bibitem{ABCRST3}
Patricia Alonso-Ruiz, Fabrice Baudoin, Li~Chen, Luke Rogers, Nageswari
  Shanmugalingam, and Alexander Teplyaev.
\newblock Besov class via heat semigroup on {D}irichlet spaces {III}: {BV}
  functions and sub-{G}aussian heat kernel estimates.
\newblock {\em Calc. Var. Partial Differential Equations}, 60(5):Paper No. 170,
  38, 2021.

\bibitem{AES25a}
Riku {Anttila}, Sylvester {Eriksson-Bique}, and Ryosuke {Shimizu}.
\newblock {Construction of self-similar energy forms and singularity of Sobolev
  spaces on Laakso-type fractal spaces}.
\newblock {\em arXiv e-prints}, page arXiv:2503.13258, March 2025.

\bibitem{AEW13}
Siva Athreya, Michael Eckhoff, and Anita Winter.
\newblock Brownian motion on {$\mathbb{R}$}-trees.
\newblock {\em Trans. Amer. Math. Soc.}, 365(6):3115--3150, 2013.

\bibitem{Bar04}
Martin~T. Barlow.
\newblock Which values of the volume growth and escape time exponent are
  possible for a graph?
\newblock {\em Rev. Mat. Iberoamericana}, 20(1):1--31, 2004.

\bibitem{BB89}
Martin~T. Barlow and Richard~F. Bass.
\newblock The construction of {B}rownian motion on the {S}ierpi\'{n}ski carpet.
\newblock {\em Ann. Inst. H. Poincar\'{e} Probab. Statist.}, 25(3):225--257,
  1989.

\bibitem{BB90}
Martin~T. Barlow and Richard~F. Bass.
\newblock On the resistance of the {S}ierpi\'{n}ski carpet.
\newblock {\em Proc. Roy. Soc. London Ser. A}, 431(1882):345--360, 1990.

\bibitem{BB92}
Martin~T. Barlow and Richard~F. Bass.
\newblock Transition densities for {B}rownian motion on the {S}ierpi\'{n}ski
  carpet.
\newblock {\em Probab. Theory Related Fields}, 91(3-4):307--330, 1992.

\bibitem{BB00}
Martin~T. Barlow and Richard~F. Bass.
\newblock Divergence form operators on fractal-like domains.
\newblock {\em J. Funct. Anal.}, 175(1):214--247, 2000.

\bibitem{BB04}
Martin~T. Barlow and Richard~F. Bass.
\newblock Stability of parabolic {H}arnack inequalities.
\newblock {\em Trans. Amer. Math. Soc.}, 356(4):1501--1533, 2004.

\bibitem{BBK06}
Martin~T. Barlow, Richard~F. Bass, and Takashi Kumagai.
\newblock Stability of parabolic {H}arnack inequalities on metric measure
  spaces.
\newblock {\em J. Math. Soc. Japan}, 58(2):485--519, 2006.

\bibitem{BBS90}
Martin~T. Barlow, Richard~F. Bass, and John~D. Sherwood.
\newblock Resistance and spectral dimension of {S}ierpi\'{n}ski carpets.
\newblock {\em J. Phys. A}, 23(6):L253--L258, 1990.

\bibitem{BP88}
Martin~T. Barlow and Edwin~A. Perkins.
\newblock Brownian motion on the {S}ierpi\'{n}ski gasket.
\newblock {\em Probab. Theory Related Fields}, 79(4):543--623, 1988.

\bibitem{Bau24}
Fabrice Baudoin.
\newblock Korevaar-{S}choen-{S}obolev spaces and critical exponents in metric
  measure spaces.
\newblock {\em Ann. Fenn. Math.}, 49(2):487--527, 2024.

\bibitem{BC23}
Fabrice Baudoin and Li~Chen.
\newblock Sobolev spaces and {P}oincar\'e{} inequalities on the {V}icsek
  fractal.
\newblock {\em Ann. Fenn. Math.}, 48(1):3--26, 2023.

\bibitem{BC24}
Fabrice Baudoin and Li~Chen.
\newblock Heat kernel gradient estimates for the {V}icsek set.
\newblock {\em Math. Nachr.}, 297(12):4450--4477, 2024.

\bibitem{BV05}
Marco Biroli and Paola~G. Vernole.
\newblock Strongly local nonlinear {D}irichlet functionals and forms.
\newblock {\em Adv. Math. Sci. Appl.}, 15(2):655--682, 2005.

\bibitem{CGQ22}
Shiping Cao, Qingsong Gu, and Hua Qiu.
\newblock {$p$}-energies on p.c.f. self-similar sets.
\newblock {\em Adv. Math.}, 405:Paper No. 108517, 58, 2022.

\bibitem{Cap24}
Marco Capolli.
\newblock An overview on {L}aakso spaces.
\newblock {\em Note Mat.}, 44(2):53--75, 2024.

\bibitem{CGYZ26}
Aobo Chen, Jin Gao, Zhenyu Yu, and Junda Zhang.
\newblock Besov--{L}ipschitz norm and {$p$}-energy measure on scale-irregular
  {V}icsek sets.
\newblock {\em J. Fractal Geom.}, 13(1-2):37--85, 2026.

\bibitem{Cla23}
Burkhard Claus.
\newblock Energy spaces, {D}irichlet forms and capacities in a nonlinear
  setting.
\newblock {\em Potential Anal.}, 58(1):159--179, 2023.

\bibitem{Eri26}
Sylvester {Eriksson-Bique}.
\newblock {On the Resistance Conjecture}.
\newblock {\em arXiv e-prints}, page arXiv:2602.05477v2, February 2026.

\bibitem{Eva08}
Steven~N. Evans.
\newblock {\em Probability and real trees}, volume 1920 of {\em Lecture Notes
  in Mathematics}.
\newblock Springer, Berlin, 2008.
\newblock Lectures from the 35th Summer School on Probability Theory held in
  Saint-Flour, July 6--23, 2005.

\bibitem{EPW06}
Steven~N. Evans, Jim Pitman, and Anita Winter.
\newblock Rayleigh processes, real trees, and root growth with re-grafting.
\newblock {\em Probab. Theory Related Fields}, 134(1):81--126, 2006.

\bibitem{FOT11}
Masatoshi Fukushima, Yoichi Oshima, and Masayoshi Takeda.
\newblock {\em Dirichlet forms and symmetric {M}arkov processes}, volume~19 of
  {\em De Gruyter Studies in Mathematics}.
\newblock Walter de Gruyter \& Co., Berlin, extended edition, 2011.

\bibitem{Gri92}
Alexander Grigor'yan.
\newblock The heat equation on noncompact {R}iemannian manifolds.
\newblock {\em Mat. Sb.}, 182(1):55--87, 1991.

\bibitem{GHH18}
Alexander Grigor'yan, Eryan Hu, and Jiaxin Hu.
\newblock Two-sided estimates of heat kernels of jump type {D}irichlet forms.
\newblock {\em Adv. Math.}, 330:433--515, 2018.

\bibitem{GHL03}
Alexander Grigor'yan, Jiaxin Hu, and Ka-Sing Lau.
\newblock Heat kernels on metric measure spaces and an application to
  semilinear elliptic equations.
\newblock {\em Trans. Amer. Math. Soc.}, 355(5):2065--2095, 2003.

\bibitem{GHL14}
Alexander Grigor'yan, Jiaxin Hu, and Ka-Sing Lau.
\newblock Heat kernels on metric measure spaces.
\newblock In {\em Geometry and analysis of fractals}, volume~88 of {\em
  Springer Proc. Math. Stat.}, pages 147--207. Springer, Heidelberg, 2014.

\bibitem{GHL15}
Alexander Grigor'yan, Jiaxin Hu, and Ka-Sing Lau.
\newblock Generalized capacity, {H}arnack inequality and heat kernels of
  {D}irichlet forms on metric measure spaces.
\newblock {\em J. Math. Soc. Japan}, 67(4):1485--1549, 2015.

\bibitem{GT12}
Alexander Grigor'yan and Andras Telcs.
\newblock Two-sided estimates of heat kernels on metric measure spaces.
\newblock {\em Ann. Probab.}, 40(3):1212--1284, 2012.

\bibitem{HKKZ00}
Ben~M. Hambly, Takashi Kumagai, Shigeo Kusuoka, and Xian~Yin Zhou.
\newblock Transition density estimates for diffusion processes on homogeneous
  random {S}ierpinski carpets.
\newblock {\em J. Math. Soc. Japan}, 52(2):373--408, 2000.

\bibitem{HPS04}
Paul~Edward Herman, Roberto Peirone, and Robert~S. Strichartz.
\newblock {$p$}-energy and {$p$}-harmonic functions on {S}ierpinski gasket type
  fractals.
\newblock {\em Potential Anal.}, 20(2):125--148, 2004.

\bibitem{Hin04}
Masanori Hino.
\newblock Integral representation of linear functionals on vector lattices and
  its application to {BV} functions on {W}iener space.
\newblock In {\em Stochastic analysis and related topics in {K}yoto}, volume~41
  of {\em Adv. Stud. Pure Math.}, pages 121--140. Math. Soc. Japan, Tokyo,
  2004.

\bibitem{KS24a}
Naotaka {Kajino} and Ryosuke {Shimizu}.
\newblock {Contraction properties and differentiability of $p$-energy forms
  with applications to nonlinear potential theory on self-similar sets}.
\newblock {\em arXiv e-prints}, page arXiv:2404.13668v2, April 2024.

\bibitem{Kig89}
Jun Kigami.
\newblock A harmonic calculus on the {S}ierpi\'{n}ski spaces.
\newblock {\em Japan J. Appl. Math.}, 6(2):259--290, 1989.

\bibitem{Kig21book}
Jun Kigami.
\newblock {\em Conductive homogeneity of compact metric spaces and construction
  of {$p$}-energy}, volume~5 of {\em Memoirs of the European Mathematical
  Society}.
\newblock European Mathematical Society (EMS), Berlin, [2023] \copyright 2023.

\bibitem{KZ92}
Shigeo Kusuoka and Xianyin Zhou.
\newblock Dirichlet forms on fractals: {P}oincar\'{e} constant and resistance.
\newblock {\em Probab. Theory Related Fields}, 93(2):169--196, 1992.

\bibitem{Laa00}
Tomi~J. Laakso.
\newblock Ahlfors {$Q$}-regular spaces with arbitrary {$Q>1$} admitting weak
  {P}oincar\'e{} inequality.
\newblock {\em Geom. Funct. Anal.}, 10(1):111--123, 2000.

\bibitem{Lax02}
Peter~D. Lax.
\newblock {\em Functional analysis}.
\newblock Pure and Applied Mathematics (New York). Wiley-Interscience [John
  Wiley \& Sons], New York, 2002.

\bibitem{Mur20}
Mathav Murugan.
\newblock On the length of chains in a metric space.
\newblock {\em J. Funct. Anal.}, 279(6):108627, 18, 2020.

\bibitem{Mur24a}
Mathav Murugan.
\newblock Heat kernel for reflected diffusion and extension property on uniform
  domains.
\newblock {\em Probab. Theory Related Fields}, 190(1-2):543--599, 2024.

\bibitem{Mur24}
Mathav {Murugan}.
\newblock Diffusions and random walks with prescribed sub-gaussian heat kernel
  estimates.
\newblock {\em Ann. Probab.}, 2025.
\newblock To appear.

\bibitem{Mur26}
Mathav {Murugan}.
\newblock {A simplified characterization of stable-like heat kernel estimates}.
\newblock {\em arXiv e-prints}, page arXiv:2602.06388, February 2026.

\bibitem{MS25}
Mathav Murugan and Ryosuke Shimizu.
\newblock First-order {S}obolev spaces, self-similar energies and energy
  measures on the {S}ierpi\'nski carpet.
\newblock {\em Comm. Pure Appl. Math.}, 78(9):1523--1608, 2025.

\bibitem{Sal92}
Laurent Saloff-Coste.
\newblock A note on {P}oincar\'{e}, {S}obolev, and {H}arnack inequalities.
\newblock {\em Internat. Math. Res. Notices}, 1992(2):27--38, 1992.

\bibitem{Sal95}
Laurent Saloff-Coste.
\newblock Parabolic {H}arnack inequality for divergence-form second-order
  differential operators.
\newblock {\em Potential Anal.}, 4(4):429--467, 1995.
\newblock Potential theory and degenerate partial differential operators
  (Parma).

\bibitem{Sas26}
K{\^o}hei {Sasaya}.
\newblock {Construction of $p$-energy measures associated with strongly local
  $p$-energy forms}.
\newblock {\em J. Funct. Anal.}, 291(1):Paper No. 111489, 2026.

\bibitem{SZ25}
Marcel {Schmidt} and Ian {Zimmermann}.
\newblock {The extended Dirichlet space and criticality theory for nonlinear
  Dirichlet forms}.
\newblock {\em arXiv e-prints}, page arXiv:2501.18391v2, January 2025.

\bibitem{Shi24}
Ryosuke Shimizu.
\newblock Construction of {$p$}-energy and associated energy measures on
  {S}ierpi\'{n}ski carpets.
\newblock {\em Trans. Amer. Math. Soc.}, 377(2):951--1032, 2024.

\bibitem{Shi24a}
Ryosuke Shimizu.
\newblock Characterizations of {S}obolev functions via {B}esov-type energy
  functionals in fractals.
\newblock {\em Potential Anal.}, 63(4):2121--2156, 2025.

\bibitem{Ste08}
Benjamin {Steinhurst}.
\newblock {Dirichlet Forms on Laakso and Barlow-Evans Fractals of Arbitrary
  Dimension}.
\newblock {\em arXiv e-prints}, page arXiv:0811.1378v3, November 2008.

\bibitem{Ste10}
Benjamin Steinhurst.
\newblock {\em Diffusions and {L}aplacians on {L}aakso, {B}arlow-{E}vans, and
  other fractals}.
\newblock ProQuest LLC, Ann Arbor, MI, 2010.
\newblock Thesis (Ph.D.)--University of Connecticut.

\bibitem{Yos95}
K{\=o}saku Yosida.
\newblock {\em Functional analysis}.
\newblock Classics in Mathematics. Springer-Verlag, Berlin, 1995.
\newblock Reprint of the sixth (1980) edition.

\end{thebibliography}

\end{document}